\documentclass[10pt]{article}
\usepackage[]{authblk}
\usepackage{amsmath}
\usepackage{graphicx}
\usepackage{latexsym}
\usepackage{amssymb}
\usepackage{amsbsy}

\DeclareSymbolFont{msbm}{U}{msb}{m}{n}
\DeclareMathSymbol{\C}{\mathalpha}{msbm}{'103}
\DeclareMathSymbol{\R}{\mathalpha}{msbm}{'122}
\DeclareMathSymbol{\Z}{\mathalpha}{msbm}{'132}
\DeclareMathSymbol{\N}{\mathalpha}{msbm}{'116}

\def\RR{\mathbb R}

\def\be{\begin{equation}}
\def\ee{\end{equation}}
\def\bea{\begin{eqnarray}}
\def\ba{\begin{array}{l}\displaystyle}
\def\eea{\end{eqnarray}}
\def\ea{\end{array}}

\begin{document}

\title{The Moment Guided Monte Carlo method for the\\ Boltzmann equation
}

\author[1]{Giacomo Dimarco\footnote
{Corresponding author address: Institut de Math\'{e}matiques de
Toulouse, UMR 5219 Universit\'{e} Paul Sabatier, 118, route de
Narbonne 31062 TOULOUSE Cedex, FRANCE. \emph{E-mail}: giacomo.dimarco@math.univ-toulouse.fr}}

\affil[1]{Universit\'{e} de Toulouse; UPS, INSA, UT1, UTM;
CNRS, UMR 5219;

Institut de Math\'{e}matiques de Toulouse; F-31062
Toulouse, France.}
%
\maketitle

\begin{abstract}
In this work we propose a generalization of the Moment Guided Monte Carlo method
developed in \cite{dimarco1}. This approach permits to reduce the variance of the particle methods
through a matching with a set of suitable macroscopic moment
equations. In order to guarantee that the moment equations provide
the correct solutions, they are coupled to the kinetic equation
through a non equilibrium term. Here, at the contrary to the previous work in which we considered the
simplified BGK operator, we deal with the full Boltzmann operator. Moreover, we introduce
an hybrid setting which permits to entirely remove the resolution of the kinetic equation
in the limit of infinite number of collisions and to consider only the solution of the compressible Euler equation.
This modification additionally reduce the statistical error with respect to our previous work and permits to perform
simulations of non equilibrium gases using only a few number of particles. We show at the end of the paper several numerical
tests which prove the efficiency and the low level of numerical noise of the method.
\end{abstract}


\section{Introduction}
The kinetic description of a gas is based on the integro-differential Boltzmann equation \cite{bird, cercignani}. This equation is satisfied by the probability distribution function of the particles
depending on the time $t\geq 0$, the
space variable $x \in \mathbb{R}^{d}$ and the velocity $v \in \mathbb{R}^{d}$. In addition to the high dimensionality of the problem,
the Boltzmann collision term that characterizes the kinetic equation is a five fold integral very hard to treat
due to its nonlinear nature and to the physical properties which need to be preserved during its resolution.

For these reasons, the numerical simulation of the Boltzmann equation with
deterministic techniques appears to be very difficult and thus, in practice, probabilistic techniques such as Direct Simulation Monte Carlo
(DSMC) methods are largely used in real situations due to their
flexibility and low computational cost compared to finite
volume or spectral methods for kinetic equations
\cite{Babovsky, bird, Cf, CPima, sliu, Nanbu80}. However, DSMC solutions
are affected by large fluctuations. Moreover, in non stationary
situations the impossibility to use time averages to reduce
fluctuations leads to poorly accurate solutions or
computationally very expensive simulations. To overcome this problem,  several methods have been developed.
We quote \cite{Cf} for an overview on both efficient and low
variance Monte Carlo methods. For applications of variance reduction
techniques to kinetic equation we mention the works of Homolle and
Hadjiconstantinou \cite{Hadji} and \cite{Hadji1}. We mention also
the research papers by Pullin
\cite{Pullin78} which developed a low diffusion particle method for
simulating compressible inviscid flows and the work of Weinan E. and co-authors \cite{WE1} on the development of multiscale numerical methods. We quote also the work of Lemou and co-authors
who developed a low variance method for the Vlasov equation close to the fluid limit \cite{Lemou}.
We finally recall the results obtained by the author of the present paper \cite{dimarco4, dimarco5} on the construction
of efficient and low variance methods for kinetic equations in transitional regimes.

A second important drawback of kinetic approaches is that the collision term becomes stiff
close to the fluid regimes. A non dimensional measure of the importance of collision is given by
the Knudsen number which is large in the rarefied regions and small in the
fluid ones. Thus, standard computational approaches lose their efficiency due to the necessity of using very small time steps
in deterministic schemes or, equivalently, a large number of collisions in probabilistic approaches. One possibility which permits to avoid the severe time step restrictions caused by the small free path limit is to construct the so-called Asymptotic Preserving schemes \cite{Lemou, Filbet, Jin2} which permits to avoid time steps dependence of the mean free path. Another possibility to avoid the excessive computational cost and the stiff regimes is to construct numerical schemes which combine continuum models with microscopic kinetic models \cite{BLT, Boyd2, degond, degond3, degond4, dimarco3, HashHassan, Letallec, tiwari_JCP, tiwari_JCP1}. However, this approach has several drawbacks, first the identification of the different regions is not trivial and second match the two models at the interfaces is not simple. For these reasons alternatively approaches are often preferred.

The two main goals of this work are first, design a particle closure of the macroscopic fluid equations
which avoid the severe time step restrictions caused by the collisional scale and second to reduce the statistical
error caused by the use of DSMC methods. This approach has explained in detail in \cite{degond4} can also be the basis of very efficient
domain decomposition strategies. In our first paper \cite{dimarco1}, we developed a method for the simple case of the BGK equation which consists in forcing particles to match prescribed sets of moments given by the solution of deterministic
equations. Now, we generalize the approach to the case of the Boltzmann equation. Moreover, we construct the scheme
in such a way that the time step restriction due to the small Knudsen number are avoided, thanks to the introduction of exponential schemes for the
resolution of the Boltzmann collision term \cite{dimarco7}. In particular, in the limit
$\varepsilon\rightarrow 0$ the scheme becomes automatically an high order method for the compressible Euler equation in which the numerical noise
completely disappears. In this sense, the method belongs to the class of the so-called asymptotic preserving methods \cite{dimarco7, Filbet, Filbet3, Jin2}.

The method is based on the following idea. The resolution of the problem is done through
the kinetic equation which is solved by DSMC method and by a set of closed moment equations which are solved by means of finite volume techniques.
In order to provide the correct solution for all the regimes of the Knudsen number, the moment
equations are coupled to the Boltzmann equation through a kinetic correction term, which takes into account
departures from thermodynamical equilibrium. In this sense, we perform a closure to the infinite set of moment equations, using the
solution of the kinetic equation. Observe that the kinetic model will be sufficient to give the correct solution to the problem.
However, as already explained, the Boltzmann equation should be resolved by very fast methods because of its too high computational cost. Typically in real situations Monte Carlo is the only possible choice. The drawback is that the solution is only poorly accurate. On the other hand, the moment equations can be resolved by using high order finite volume techniques. The result is a method which has a computational cost comparable to the cost of DSMC method but which a smaller level of numerical noise. Moreover, the stiffness of the problem is removed thanks to the exponential resolution of the collision term. This gives a method which is faster than classical DSMC methods \cite{Babovsky, bird} in the limit of the Knudsen number which goes to zero i.e. in the fluid limit. Finally, in order to recover at each time step the same moments from the solution of the kinetic equation and from the solution of the moments equations and then advance in time, we constrain the DSMC method to match the moments obtained through the deterministic resolution of the macroscopic equations in such a way that the higher accuracy of the moments resolution improves the accuracy of the DSMC method. We experimentally show that this is indeed the case.

The remainder of the paper is organized as follows. In the next
section we recall some basic notions on the Boltzmann equations and
its fluid limit. The scheme is described in section 3. In particular, the numerical method used for the kinetic equation is described
in section 3.1, the numerical method used for the moment equation is described in section 3.2, while the matching procedure is described in section 3.3. In section 4 numerical examples which demonstrate the
capability of the method are presented. Finally some conclusions and remarks are drawn in the last section.

\section{The Boltzmann equation and its fluid limit}

We consider the Boltzmann equation of rarefied gas dynamics \cite{cercignani}
\be
\partial_t f + v\cdot\nabla_{x}f
=Q(f,f)\label{eq:1}\ee equipped of an initial data \be
f|_{t=0}=f_{0}.\label{eq:ini}\ee Here $f(x,v,t)$ is a non negative
function describing the time evolution of the distribution of
particles with velocity $v \in \R^{3}$ and position $x \in \Omega
\subset \R^{d_x}$ at time $ t
> 0$. In the sequel for notation simplicity we will omit the
dependence of $f$ from the independent variables $x,v,t$ unless
strictly necessary. The operator $Q(f,f)$ which describes particle
interactions has the form \be Q(f,f)=\int_{\RR^3\times S^2}
B(|v-v_*|,n) [f(v')f(v'_*)-f(v)f(v_*)]\,dv_*\,dn \label{eq:Q} \ee
where \be v'=v+\frac12(v-v_*)+\frac12|v-v_*|n,\quad
v'_*=v+\frac12(v-v_*)-\frac12|v-v_*|n, \ee and $B(|v-v_*|,n)$ is a
nonnegative collision kernel characterizing the microscopic details
of the collision given by \be
B(|v-v_*|,n)=\sigma\left(\frac{(v-v_*)}{|v-v_*|}\cdot n\right)
|v-v_*|^\gamma, \label{eq:coll}\ee with $\gamma \in [0,3)$ and
$\sigma$ the collision cross section which in general depends on the
collision angle. The case $\gamma=1$ is referred to as hard spheres
case, whereas the simplified situation $\gamma=0$, is referred to as
Maxwell case. Note that in most applications the angle dependence is
ignored and $\sigma$ is assumed constant \cite{bird}.

The operator $Q(f,f)$ is such that the local conservation
properties are satisfied \be\int_{\R^{3}} mQ(f,f) dv=:\langle
mQ(f,f)\rangle=0 \label{eq:QC}\ee where $m(v)=(1,v,\frac{|v|^2}{2})$
are the collision invariants. In addition it satisfies the entropy
inequality \be \frac{d}{dt}\int_{\R^{3}}f\log f\,dv = \int_{\R^{3}}
Q(f,f)\log f dv \leq 0. \label{eq:entropy} \ee Integrating
(\ref{eq:1}) against its invariants in the velocity space leads to
the following set of non closed conservations laws \be
\partial_t \langle mf\rangle+\nabla_x\cdot
\langle vmf\rangle=0.\label{eq:macr}\ee Equilibrium functions for
the operator $Q(f,f)$ (i.e. solutions of $Q(f,f)=0$) are local
Maxwellian of the form \be M_{f}(\rho,u,T)=\frac{\rho}{(2\pi
T)^{3/2}}\exp\left(\frac{-|u-v|^{2}}{2T}\right), \label{eq:M}\ee
where $\rho$, $u$, $T$ are the density, mean velocity and
temperature of the gas in the x-position and at time $t$ defined as
\be (\rho,\rho u,\rho e)=\langle m f \rangle, \qquad E=\rho e, \qquad
T=\frac1{3\rho}(E-\rho|u|^2). \ee We will denote by \be
U=(\rho,\rho u,E),\qquad M[U]=M_f, \ee for definition of $M[U]$ we have \be U=\langle
mM[U]\rangle.\ee When the number of collisions is large which means that the mean free path between particles is
very small, it is convenient to rescale the space
and time variables in (\ref{eq:1}) as \be x'=\varepsilon x, \ \
t'=\varepsilon t\ee which leads to \be
\partial_t f + v\cdot\nabla_{x}f
=\frac1{\varepsilon}Q(f,f)\label{eq:1b}\ee where $\varepsilon$ is the so-called Knudsen number
proportional to the mean free path while the primes have been omitted for simplicity.

If we pass formally to the limit for $\varepsilon\rightarrow 0$ we get
$f\rightarrow M[U]$. This leads to the well-known hyperbolic system
of compressible Euler equations for the macroscopic variables $U$
\be
\partial_t U+\nabla_x\cdot F(U)=0
\label{eq:Euler} \ee with
\[
F(U)=\langle vmM[U]\rangle=(\rho u, \varrho u \otimes u+pI,
Eu+pu),\quad p=\rho T,
\]
where $I$ is the identity matrix.

\section{The Monte Carlo method for the Boltzmann equation}

The method is based on the
decomposition \be f=M[U]+g.\ee In the above formula the function
$g$ represents the non-equilibrium part of the distribution
function. Since $f$ and $M[U]$ share the same first three
moments we get \be \langle mg\rangle=0,\ee which means that the non-equilibrium part $g$ has always non positive parts. With the above definitions $U$, $f$ and $g$ satisfy
the following system of equations \begin{eqnarray}
\partial_t U+\nabla_x\cdot
F(U)+\nabla_x\cdot\langle vmg\rangle&=&0,\label{eq:2}\\
\partial_t f+v\cdot\nabla_x
f&=&\frac{1}{\varepsilon}Q(f,f), \label{eq:3}\\
g=f-M[U]. \end{eqnarray}  The proof of the above statement can be found in \cite{degond3}, where we also remind for details.

As explained in the introduction, the goal of the method is to solve equation (\ref{eq:2}), which means to close the system of moments equations.
In order to accomplish that task, we need to know the time evolution of $g$, which is given by the solution of equation (\ref{eq:3}).
Of course solving equation (\ref{eq:3}) will be sufficient to compute the solution also of (\ref{eq:2}).
However, the numerical solution of equation (\ref{eq:3}) is very expensive and thus we seek for very fast solvers, as for instance
the Monte Carlo method. The counterpart of fast solvers is that the solution of (\ref{eq:3}) will be normally known with very low accuracy.
Thus, our statement is the following: if the departure from equilibrium $g$ is small, the low accuracy and the high level of numerical noise with which $g$ is known from (\ref{eq:3}) will still permit to compute the moments in (\ref{eq:2}) with an error which will be smaller than the error with which equation (\ref{eq:3}) is solved. We will numerically show that this is indeed the case. In particular, thanks to the way in which the collisional operator is solved, the coupling method will permit to make disappear $g$ from the computation of the moments in the limit $\varepsilon\rightarrow 0$.
The precise time marching procedure will be detailed next, the method can be summarized as following: at each
time step $t^n$
\begin{enumerate}
 \item Solve the kinetic equation (\ref{eq:3}) with a Monte Carlo scheme and obtain the
 set of moments $U^*=\langle mf^*\rangle$.
 \item Solve the fluid equation (\ref{eq:2}) with a finite volume/difference scheme using the particles solution to close the system and obtain a second set of moments $U^{n+1}$. This means evaluate the term $\nabla_x\cdot\langle vmg\rangle$ with particles.
 \item Match the moments of the kinetic solution with those of the fluid solution through a linear transformation of the
 samples velocity $f^{n+1}=T(f^*)$ so that $\langle mf^{n+1}\rangle=U^{n+1}$ and continue.
 \end{enumerate}
In the first step, the method can be generalized by substituting the Monte Carlo method with any low accurate but very fast solver. The other two steps are the key points of the method they involve the evaluation of
$\nabla_x\cdot\langle vmg\rangle$ and the matching procedure detailed next.

We observe that it is possible to improve
the method, adding to system (\ref{eq:2}) additional equations for
the time evolution of higher order moments and get \be
\partial_t \langle m_n f\rangle+\nabla_x \cdot\langle v
m_nf\rangle=\langle m_nQ(f,f)\rangle\label{eq:5}\ee with $m_n=v^n$
and $n\geq 3$. However, this will lead to more complex equations which discretiazion will not be easy. We will not discuss this possibility here.


\subsection{Solution of Boltzmann equation by exponential Runge-Kutta Monte Carlo methods}\label{sec:Bol}
The starting point in the solution of
kinetic equations by Monte Carlo methods is given by an operator splitting of
(\ref{eq:1}) or equivalently (\ref{eq:3}) in a time interval $[0,\Delta t]$ between relaxation
\be
\partial_t f=\frac{1}{\varepsilon}Q(f,f),\label{eq:9}\ee
and free transport \be
\partial_t f + v\cdot\nabla_{x}f
=0. \label{eq:tr} \ee Even if this splitting is limited to first order it permits to
treat separately the hyperbolic free transport from the stiff
relaxation step. Observe that high order splitting schemes are possible with Monte Carlo techniques. However, we did not consider this possibility in this work, we remind to the future for the exploration of high order time discretizations.

We introduce now a space discretization of mesh size $\Delta x$ and a time
discretization of time step $\Delta t$. The discretization of the
domain is not needed for the transport step which is solved exactly
by pushing the particles. However, it is necessary to solve the collision part of the
problem which acts locally in space and thus it turns to be necessary to
solve the full problem. In Monte Carlo simulations the distribution
function $f$ is discretized by a finite set of $N$ particles
\begin{eqnarray}
& & f = \frac{m_p}{N}\sum_{i=1}^N  \alpha_i(t) \, \delta(x-X_i(t))
\delta(v-V_i(t)), \label{particle}
\end{eqnarray}
where $X_i(t)$ represents the particle position in the three spatial
directions, $V_i(t)$ the particle velocities in the velocity space,
and $\alpha_i(t)$ the weight to associate to each
particle. The constant $m_p$ is defined at the beginning of the
computation in the following way \be
m_p=\frac{1}{N}\int_{\Omega}\int_{\R^{3}} f(t=0) dv dx.\ee
During the transport step (\ref{eq:tr}), the particles
move to their next positions according to \be X_i(t+\Delta
t)=X_i(t)+V_i(t)\Delta t\label{transport}\ee where (\ref{particle})
with (\ref{transport}) and a constant $\alpha_i(t)$ represents an exact solution
of equation (\ref{eq:tr}). The role of the function $\alpha$ will be clarified in the moment matching procedure.

We consider now the solution of the collision step (\ref{eq:9}).
The effect of this step is to change the shape of the velocity
distribution and to project it towards the equilibrium distribution
leaving unchanged mass, momentum and energy. The collision operator
acts locally in space which means that we solve it independently in
each spatial cell. The particles are assumed to have all the
same weight in one cell. The initial data of this
step is given by the solution of the transport step because of the splitting.

There exists many different ways to solve the collisions, see for instance \cite{Babovsky, bird, cercignani}, however,
we seek for a method which will be stable for all choices of $\Delta t$ independently of $\varepsilon$.
This will permit to avoid the stiffness of the equation and at the same time to reduce the statistical fluctuations.
Since we aim at developing unconditionally stable schemes, the most
natural choice would be to use implicit solvers applied to
(\ref{eq:9}). Unfortunately the use of fully implicit schemes for
(\ref{eq:9}) is unpracticable in the case of DSMC methods. To overcome these difficulties, we introduce a method based on exponential integrators. First we rewrite the homogeneous equation (\ref{eq:9}) in the form
\be
\partial_t f=\frac{1}{\varepsilon}(P(f,f)-\mu f),\label{eq:10}
\ee where $P(f,f)=Q(f,f)+\mu f$ and $\mu>0$ is a constant such that
$P(f,f)\geq 0$. Typically $\mu$ is an estimate
of the largest rate of the negative term in the Boltzmann operator. For example for
Maxwellian molecules we have \be P(f,f)= Q^+(f,f)(v) = \int_{\RR^3}
\int_{S^2} b_0(\cos\theta) f(v^{\prime})f(v_{\ast}^{\prime}) \,
d\omega\,dv_{\ast}, \label{eq:Qpiu} \ee and $\mu=\varrho$.
By construction the following property holds \be \frac1{\mu}\langle m
P(f,f)\rangle = \langle m f\rangle=U, \ee this means that $P(f,f)/\mu$ is a
density function. The homogeneous equation can be written now, adding and subtracting $M[U]$, in the form
\be
\partial_t
f=\frac{\mu}{\varepsilon}\left(\frac{P(f,f)}{\mu}-M[U]\right)+
\frac{\mu}{\varepsilon}(M[U]-f).\label{eq:11}\ee Note that even if
$M[U]$ is nonlinear in $f$ thanks to the conservation properties
(\ref{eq:QC}), it remains constant during the relaxation process.

We apply now to the reformulated equation (\ref{eq:11}) an approach based on the so-called exponential
integrators where the exact solution of the linear part is used
for the construction of the numerical method \cite{dimarco7}. In order to derive the methods it is useful to rewrite
(\ref{eq:11}) as \be \frac{\partial (f-M[U])e^{\mu
{t}/{\varepsilon}}}{\partial t}=\frac1{\varepsilon}(P(f,f)-\mu
M[U])e^{\mu {t}/{\varepsilon}}. \label{eq:15} \ee

The above form is readily obtained if one multiplies (\ref{eq:11})
by the integrating factor $\exp(\mu t/\varepsilon)$ and takes into
account the fact that $M[U]$ does not depend of time during the collisional process. A class of
explicit exponential Runge-Kutta schemes is then obtained by direct
application of an explicit Runge-Kutta method to (\ref{eq:15}). More
generally we can consider the family of methods characterized by
\begin{eqnarray}
\nonumber F^{(i)}&=&e^{-c_i\mu\Delta
t/\varepsilon}f^n+\frac{\mu\Delta t}{\varepsilon}
\sum_{j=1}^{i-1}A_{ij}(\mu\Delta
t/\varepsilon)\left(\frac{P(F^{(j)},F^{(j)})}{\mu}-M[U^n]\right)\\
\label{eq:rk1}
\\[-.25cm]
\nonumber &+& \left(1-e^{-c_i\mu\Delta t/\varepsilon}\right)M[U^n],
\qquad
i=1,\ldots,\nu\\
\nonumber f^{n+1}&=&e^{-\mu\Delta
t/\varepsilon}f^*+\frac{\mu\Delta
t}{\varepsilon}\sum_{i=1}^{\nu}W_i(\mu\Delta
t/\varepsilon)\left(\frac{P(F^{(i)},F^{(i)})}{\mu}-M[U^n]\right)\\
\label{eq:rk2}
\\[-.25cm]
\nonumber &+&\left(1-e^{-\mu\Delta t/\varepsilon}\right)M[U^n],
\end{eqnarray}
where $f^*$ is the distribution function value after the transport step, $F^{(i)}$ are called stages,
$c_i \geq 0$, while the coefficients $A_{ij}$ and the weights $W_i$
are such that
\[
A_{ij}(0)=a_{ij},\quad W_i(0)=w_i,\quad i,j=1,\ldots,\nu
\]
with coefficients $a_{ij}$, $c_i$ and weights $w_i$ given by a standard
explicit Runge-Kutta method called the underlying method. Various
schemes come from the different choices of the underlying method.
The most popular approach is the integrating factor (IF) method.
For the so-called Integrating Factor methods, which correspond to
a direct application of the underlying method to (\ref{eq:15}), we
have
\begin{eqnarray}
\nonumber
A_{ij}(\lambda)&=&a_{ij}e^{-(c_i-c_j)\lambda},\quad
i,j=1,\ldots,\nu,\quad i > j\\[-.1cm]
\label{eq:jin}\\
\nonumber W_i(\lambda)&=&w_i e^{-(1-c_i)\lambda},\quad
i=1,\ldots,\nu,
\end{eqnarray}
with $\lambda=\mu\Delta t/\varepsilon$. In the sequel, among all possible choices for discretizing (\ref{eq:10}), we will use the
first order in time integrating factor scheme which reads
\begin{equation}
f^{n+1}=e^{-\frac{\mu \Delta t}{\varepsilon}}
f^{*}+\frac{\mu\Delta t}{\varepsilon}e^{-\frac{\mu \Delta
t}{\varepsilon}}\left(\frac{P(f^*,f^*)}{\mu}-M[U^{n}]\right)+\left(1-e^{-\frac{\mu
\Delta t}{\varepsilon}}\right)M[U^{n}], \label{eq:collfirst}
\end{equation}
which is based on the simple firs order in time explicit Euler Runge-Kutta method.
Observe that this method permits to avoid the stiffness of the collisional operator. This is, in fact, unconditionally stable for all choices of time step $\Delta t$. At the same time, the method proposed is, in the limit $\varepsilon\rightarrow 0$, nothing else than a kinetic scheme for the numerical solution of the compressible Euler equation. Observe, in fact, that when $\varepsilon\rightarrow 0$ at each time step the distribution function is projected on the equilibrium distribution $M[U]$.

The Monte Carlo interpretation of equation (\ref{eq:collfirst}) is the following. If we define \be A=e^{-\frac{\mu \Delta t}{\varepsilon}}, \qquad B=\frac{\mu\Delta t}{\varepsilon}e^{-\frac{\mu \Delta
t}{\varepsilon}}, \qquad C=\left(1-\frac{\mu\Delta t}{\varepsilon}e^{-\frac{\mu \Delta
t}{\varepsilon}}-e^{-\frac{\mu
\Delta t}{\varepsilon}}\right), \quad C=1-A-B,\ee
with probability $A$ the velocity of the particle does not change during the collision step. With probability $B$ the particle occurs in one collision, the details of the collision are described by the operator $P$, this is, for instance, the gain part of the collisional process, described by Bird \cite{bird}. With probability $C$ the velocity of the particle is replaced by a particle sampled from the Maxwellian distribution $M[U]$.
We will discuss in the matching section the way in which we link the solution of the kinetic equation to the solution of the moment equations.

\subsection{Solution of the moments equations}

In this section we detail the way in which the moments
equations have been discretized. We suppose for the construction of the scheme that the function $g^{n}$ is known at time $n$.
The numerical scheme proposed will take advantage from the fact that part of these moment equations are in fact nothing else that the compressible Euler equations \be \underbrace{\partial_t U+\nabla_x\cdot
F(U)}_{\hbox{Euler equations}}+\nabla_x\cdot \langle vmg \rangle=0.
\label{eq:2b}\ee We first solve the set
of compressible Euler equations and then we consider
the discretization of the kinetic flux $\nabla_x\cdot\langle
vmg\rangle$. For the space discretization of the
equilibrium fluxes we use a second order MUSCL central
scheme. For simplicity, we indicate in the same way the numerical
flux in one or in more spatial dimensions. Thus, we have\be \frac{U^{*}_j-U^{n}_j}{\Delta t}
+\frac{\psi_{j+1/2}(U^n)-\psi_{j-1/2}(U^n)}{\Delta
x}=0.\label{eq:discmom}\ee The discrete flux reads \cite{JX, leveque:numerical-methods}\be
   \psi_{j+1/2}(U^n)=\frac{1}{2}(F(U^n_{j})+F(U^n_{j+1}))-\frac{1}{2}\alpha(U^n_{j+1}-U^n_{i})+\frac{1}{4}(\sigma^{n,+}_j-\sigma^{n,-}_{j+1})
\ee where \be \sigma^{n,\pm}_j=\left(F(U^n_{j+1})\pm \alpha
U^n_{j+1}-F(U^n_{j})\mp \alpha
U^n_{j}\right)\varphi(\chi^{n,\pm}_j)\ee with
$\varphi$ the Van Leer slope limiter\be
\varphi(\chi)=\frac{|\chi|+\chi}{1+\chi},\ee finally, the variable $\chi^{\pm}$ is defined as
following \be \chi^{n,\pm}_j=\frac{F(U^n_{j})\pm \alpha
U^n_{j}-F(U^n_{j-1})\mp \alpha U^n_{j-1}}{F(U^n_{j+1})\pm \alpha
U^n_{j+1}-F(U^n_{j})\mp \alpha U^n_{j}}\ee where the above ratio of
vectors is defined componentwise and
$\alpha$ is equal to the largest eigenvalue of the Euler system.

We now discuss how to discretize the non equilibrium term
$\nabla_x\cdot<vmg>$. To this aim, the first step is to introduce a filter to
eliminate some fluctuations. We choose the weighted moving average method, which is a convolution
of the pointwise value of $g$ with a fixed weighting function. Thus, given the function $g^{n}$, we define
\be \widetilde{g}_{j}^{n}=\frac{1}{2K+1}\sum_{k=-K}^{k=K}\omega_k g_{j-k}, \ \sum_{k=-K}^{k=K}\omega_k=1.\label{eq:movav}\ee
The smoothing permits to reduce the fluctuations but on the other hand it causes a degradation of the accuracy
in the solution. Thus, in practice, we use a weak filter in order to keep the solution as precise as possible. In the numerical simulations, we used $K=1$, $\omega_{-1}=\omega_{1}=1/6$ and $\omega_{0}=2/3$.
Observe that a different smoothing can be used at the particle level instead that at the moments level. This different smoothing has to be applied during the reconstruction of the moments of $g$ from the particles. We discuss this possibility in the next section. Once that the filer is applied, the non equilibrium term is discretized with the same second order MUSCL scheme used for computing the flux of the Euler equations where, however, the numerical diffusion is taken equal to zero. This is because this term is at the leading order a diffusion term, and thus it does not need numerical diffusion for stability to be assured. The complete scheme for the macroscopic equations reads\be \frac{U^{*}_j-U^{n}_j}{\Delta t}
+\frac{\psi_{j+1/2}(U^n)-\psi_{j-1/2}(U^n)}{\Delta
x}=0,\label{eq:discmom2}\ee
\be \frac{U^{n+1}_j-U^{*}_j}{\Delta t}
+\frac{\Psi_{j+1/2}(<vm\widetilde{g}^n>)-\Psi_{j-1/2}(<vm\widetilde{g}^n>)}{\Delta
x}=0.\label{eq:discmom3}\ee
where $\Psi_{j+1/2}(<vm\widetilde{g}^n>)$ is defined as
\be
   \Psi_{j+1/2}(<vm\widetilde{g}^n>)=\frac{1}{2}(<vm\widetilde{g}^n_{j}>+<vmg^n_{j+1}>)+\frac{1}{4}(\sigma^{n}_{g,j}-\sigma^{n}_{g,j+1})
\ee
with
\be \sigma^{n}_{g,j}=\left(<vm\widetilde{g}^n_{j+1}>-<vm\widetilde{g}^n_{j}>\right)\varphi(\chi^{n}_{g,j})\ee
and \be \chi^{n}_{g,j}=\frac{<vm\widetilde{g}^n_{j}>-<vm\widetilde{g}^n_{j-1}>}{<vm\widetilde{g}^n_{j+1}>-<vm\widetilde{g}^n_{j}>}.\ee
where again the above ratio of vectors is defined componentwise.


\subsection{The Moment Matching}
\label{sec:mm}
In this section we discuss the details of the moment matching and the coupling between the kinetic and the macroscopic model. Suppose, known $f^{n-1}$, $\widetilde{g}^{n-1}$ and $U^{n-1}$. The time marching procedure is the following\newpage
\begin{enumerate}
  \item Solve the moments equations from time $n-1$ to time $n$ using equations (\ref{eq:discmom2}-\ref{eq:discmom3}) this gives $U^{n}$.
  \item Solve the transport part of the kinetic equation from time $n-1$ to the intermediate stage $*$ pushing the particles as in equation (\ref{transport}).
  \item Compute the moments of the Boltzmann equation after the transport of the particles. This gives $\widetilde{U}^{n}$, the collision part, being conservative, does not alter the moments. The reconstruction of the moments from the particles will be discussed next.
  \item Match the moments of the kinetic and the macroscopic equations forcing the particles to have the moments $U^{n}$, i.e. $\widetilde{U}^{n}\rightarrow U^{n}$.
  \item Compute the collision step of the kinetic equation through (\ref{eq:collfirst}) this gives $f^{n}$.
  \item Use the computed value of $f^{n}$ to compute $g^{n}$.
  \item Compute $\widetilde{g}^{n}$ with the moving average technique (equation \ref{eq:movav}).
  \item Plug the value of $\widetilde{g}^{n}$ in the moments equations and compute $U^{n+1}$ using (\ref{eq:discmom2}-\ref{eq:discmom3}) and continue.
\end{enumerate}

We discuss now, how to reconstruct the moments from the particles (point $3$ of the above procedure), the matching of the different moments (point 4) and the computations of $g^{n}$ from $f^{n}$ (point 6).

We first consider the matching of the mass. To this aim, let consider
the set of particles $X_1,\ldots,X_{N_{\mathcal{I}_j}}$ inside the
cell $\mathcal{I}_j$ after the transport step. The corresponding mass is computed as
\be
\widetilde{\varrho}^{n}_{j}=\sum_{i\in
\mathcal{I}_{j}}^{N_{\mathcal{I}_j}} m_p \alpha^{n-1}_i, \qquad j=1,N_x\label{massrecon}\ee others possible reconstructions will be discussed next.
In our previous paper, \cite{dimarco1}, among the possible techniques that can be used to restore a prescribed density we choose to replicate or discard particles inside the cells. Here, in order to restore the mass we assign to each particle inside the cell $j$ the new value
\be \alpha_{i}^{n}=\frac{\varrho_{j}^{n}}{N_{\mathcal{I}_j}} \ \forall i=1,..,N_{\mathcal{I}_j},\label{match1}\ee where $\varrho_{j}^{n}$ is density computed from the solution of the moments equations.
This rescaling of the mass introduce an error in the evaluation of the distribution function $f$. However, we recall that the two models, the microscopic and macroscopic one, give the same solution in term of the moments apart from the numerical errors. In particular, the DSMC method gives solutions which oscillates around the exact solution of the problem. This implies that, the above renormalization of the mass only force the density to be closer to the value furnished by the exact solution of the problem. In fact, the moments equations furnish solutions which contain less fluctuations with respect to the DSMC method. The reason for this lower level of noise of the macroscopic equations is that only the perturbation from the equilibrium $g$ is computed statistically and not the full solution as in the original DSMC method.

We discuss now the matching of the momentum and of the energy, which are done after the matching of the density. In the Monte Carlo setting these moments can be obtained by the following piecewise constant reconstructions \be
\widetilde{u}^{n}_{j}=\frac{1}{\varrho^{n}_{j}N_{\mathcal{I}_j}}\sum_{i\in
\mathcal{I}_{j}}^{N_{\mathcal{I}_j}} V_i \qquad
\widetilde{e}_{j}^{n}=\frac{1}{\varrho^{n}_{j}N_{\mathcal{I}_j}}\sum_{i\in
\mathcal{I}_{j}}^{N_{\mathcal{I}_j}} \frac{1}{2}|V_i|^2, \label{momrecon}\ee where
$\widetilde{u}^{n}_{j}$ is a vector representing the mean velocity in the three
spatial directions. To match mean velocity and energy with those of equations (\ref{eq:discmom2}-\ref{eq:discmom3}) we then apply the transformation described in \cite{Cf} which permits to get
a new set of velocities $V_i^*$ defined by \be V_i^*=(V_i-\widetilde{u}^{n}_{j})/c+u^{n}_{j}\quad
c=\sqrt{\frac{\widetilde{e}^{n}_{j}-(\widetilde{u}^{n}_{j})^2}{e^{n}_{j}-(u^{n}_{j})^2}},\quad
i=1,\ldots,N_{\mathcal{I}_j} \label{match2}\ee which gives
\[
\frac{1}{\varrho^{n}_{j}N_{\mathcal{I}_j}}\sum_{i\in
\mathcal{I}_{j}}^{N_{\mathcal{I}_j}} V_i^*=u^{n}_{j},\qquad
\frac{1}{\varrho^{n}_{j}N_{\mathcal{I}_j}}\sum_{i\in
\mathcal{I}_{j}}^{N_{\mathcal{I}_j}} \frac{1}{2}|V_i^*|^2=e^{n}_{j}.
\]

After the matching at time $n$, the next point in the time marching procedure described before, regards the collisions. This step has been already described in section \ref{sec:Bol}. It uses equation (\ref{eq:collfirst}) to compute from $f^{*}$, which is now the distribution function value after the transport and the matching, the new distribution $f^n$.

We now discuss the sixth point of the method: the computation of $g^{n}$ from $f^{n}$. The perturbation $g^{n}$ is given by
\begin{eqnarray}
   g^{n}_{j}&=&f^{n}_{j}-M[U^{n}_{j}]=\nonumber\\&=& e^{-\frac{\mu \Delta t}{\varepsilon}}
f^{n-1}+\frac{\mu\Delta t}{\varepsilon}e^{-\frac{\mu \Delta
t}{\varepsilon}}\frac{P(f^{n-1},f^{n-1})}{\mu}+\left(1-e^{-\frac{\mu
\Delta t}{\varepsilon}}-\frac{\mu\Delta t}{\varepsilon}e^{-\frac{\mu \Delta
t}{\varepsilon}}\right)M[U^{n}_{j}]-M[U^{n}_{j}]\nonumber  \\
   &=&e^{-\frac{\mu \Delta t}{\varepsilon}}
f^{n-1}+\frac{\mu\Delta t}{\varepsilon}e^{-\frac{\mu \Delta
t}{\varepsilon}}\frac{P(f^{n-1},f^{n-1})}{\mu}-\left(e^{-\frac{\mu
\Delta t}{\varepsilon}}+\frac{\mu\Delta t}{\varepsilon}e^{-\frac{\mu \Delta
t}{\varepsilon}}\right)M[U^{n}_{j}],
\end{eqnarray}
where once again we used the
first order in time integrating factor scheme (\ref{eq:collfirst}).
The above expression tells us that the moments of $g^{n}_{j}$ can be obtained as a contribution of two terms. The first term is obtained by reconstructing from the particles the moments as in (\ref{massrecon}, \ref{momrecon}). The second term is obtained by integrating over the velocity space the analytic expression of the Maxwellian distribution. This implies directly that the moments related to the second part does not contain any statistical error. Finally, observe that in the limit $\varepsilon\rightarrow 0$ the contribution of the perturbation $g$ goes to zero. As a consequence, the method proposed resolves a macroscopic system which is as closer as the Knudsen number diminishes to the compressible Euler equation. In the limit of $\varepsilon=0$, the scheme exactly solve the compressible Euler equations without any source of statistical error.
Finally, the discretized moments equations (\ref{eq:discmom3}) can be rewritten as
\begin{eqnarray} \frac{U^{n+1}_j-U^{*}_j}{\Delta t}
&=&-\frac{\Psi_{j+1/2}(<vm\widetilde{g}^n>)-\Psi_{j-1/2}(<vm\widetilde{g}^n>)}{\Delta
x}= \nonumber\\ &=& -\frac{\Psi_{j+1/2}(<vm \left(e^{-\frac{\mu
\Delta t}{\varepsilon}}+\frac{\mu\Delta t}{\varepsilon}e^{-\frac{\mu \Delta
t}{\varepsilon}}\right)(\widetilde{g'}-\widetilde{M})^n>)}{\Delta
x}+ \nonumber\\ &+& \frac{\Psi_{j-1/2}(<vm \left(e^{-\frac{\mu
\Delta t}{\varepsilon}}+\frac{\mu\Delta t}{\varepsilon}e^{-\frac{\mu \Delta
t}{\varepsilon}}\right)(\widetilde{g'}-\widetilde{M})^n>)}{\Delta
x}
\end{eqnarray}
where \be (g')^{n}_{j}=\frac{e^{-\frac{\mu \Delta t}{\varepsilon}}
f^{n-1}+\frac{\mu\Delta t}{\varepsilon}e^{-\frac{\mu \Delta
t}{\varepsilon}}\frac{P(f^{n-1},f^{n-1})}{\mu}}{\left(e^{-\frac{\mu
\Delta t}{\varepsilon}}+\frac{\mu\Delta t}{\varepsilon}e^{-\frac{\mu \Delta
t}{\varepsilon}}\right)}\ee
As $\Delta t /\varepsilon$ grows, which means that the system
approaches the equilibrium, the contribution of the kinetic term
vanishes even though it is evaluated through particles. Observe that, this does
not happen if we just compute the kinetic term $\nabla_x\cdot \langle
vmg\rangle$ from the particles without considering the structure of
the distribution function $f$. This dramatically decreases
fluctuations when the Knudsen number is small.

As a conclusion for this section, we discuss some different reconstructions of the moments starting from the particles.
Instead of using the piecewise constant reconstructions (\ref{massrecon}-\ref{momrecon}), one can think to smoother reconstructions which can probably further diminish the statistical noise. This procedure will be alternative to the moving average algorithm used for computing $\widetilde{g}$ from $g$ (\ref{eq:movav}). In other words, we first smooth the contribution of the particles in the computation of the macroscopic quantities and we then after just apply the finite volume method to set of resulting moments equations without using average algorithms.

A general way to compute the integral over the velocity space is to sum over the particles with the use of weight functions. In this strategy, the value of the perturbation $g$ will be given by
\be \langle vmg \rangle=
\sum_{i}^{N}m_p\alpha_i v_i B(X_i-x_{j}) m_i \qquad j=1,..,N_x\ee where $B\geq 0$ is a suitable weight function s.t.
\[
\int_{\RR} B(x)\,dx=1.
\]
For example, $B(x)=1$ if $|x|\leq \Delta x/2$ and $B(x)=0$ elsewhere,
gives rise to the previous piecewise constant reconstruction, while  $B(x)=1/\Delta x$ if $|x|\leq \Delta x/2$ and $B(x)=0$ elsewhere,
is the so called Nearest Grid Point procedure in plasma
physics \cite{birsdall}. Smoother reconstructions can be recovered
by convolving the samples with a bell-shaped weight like a B-spline. Note that the value $B(x)$ has a strong influence on
the fluctuations in the reconstructed function, and in general
should be selected as a good compromise between fluctuations and
resolution as for the moving average method.

\begin{figure}
\begin{center}
\includegraphics[scale=0.39]{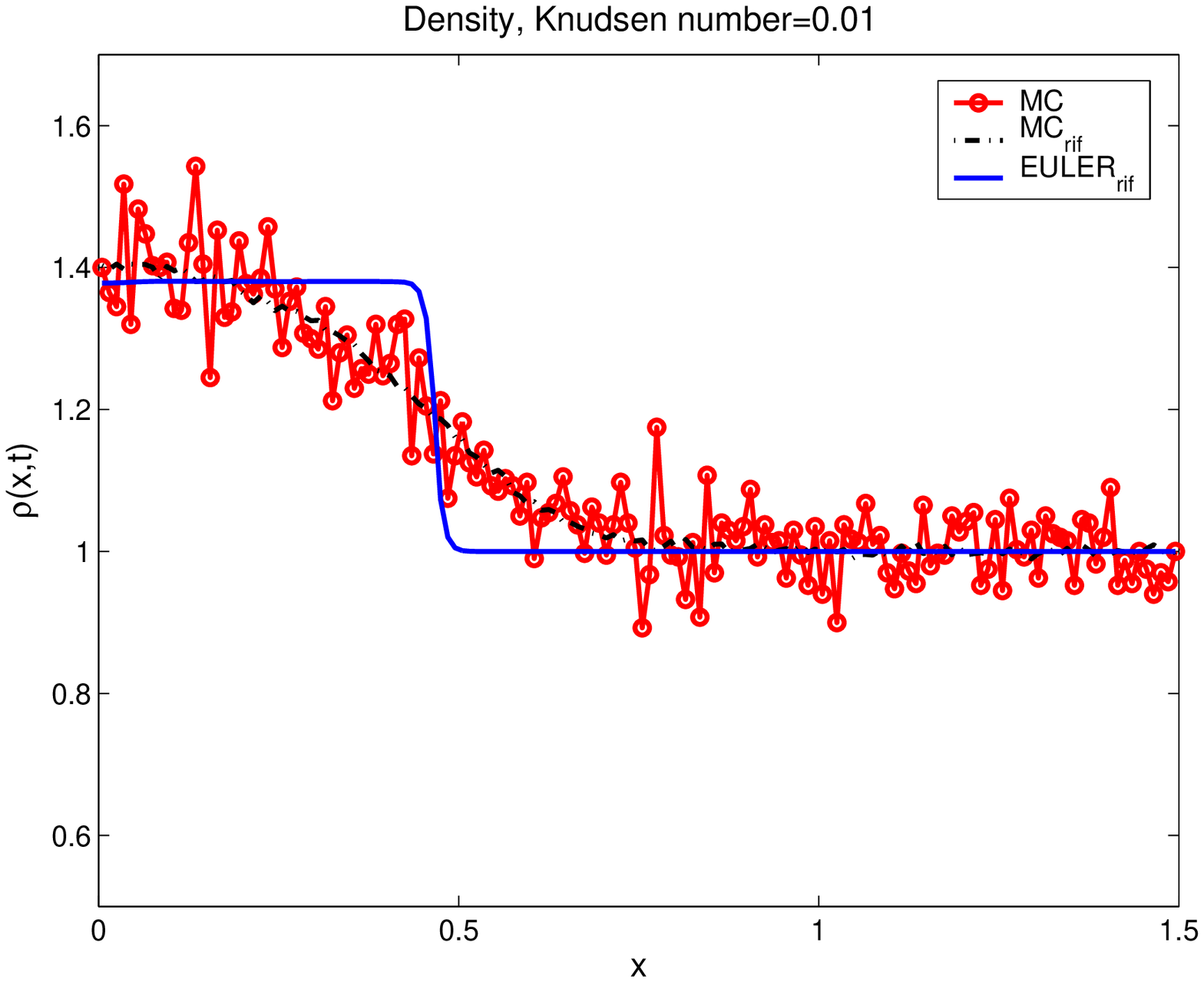}
\includegraphics[scale=0.39]{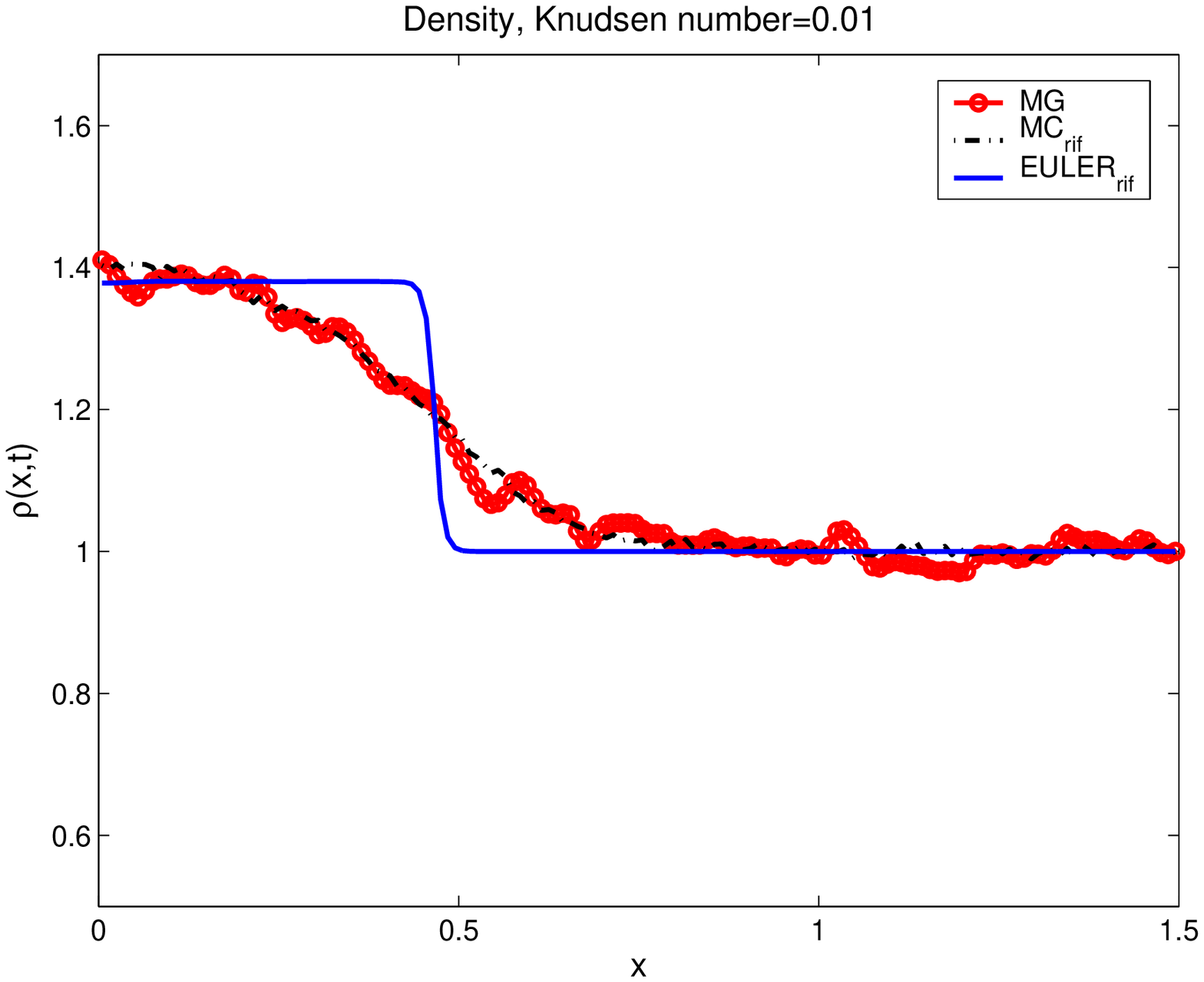}
\includegraphics[scale=0.39]{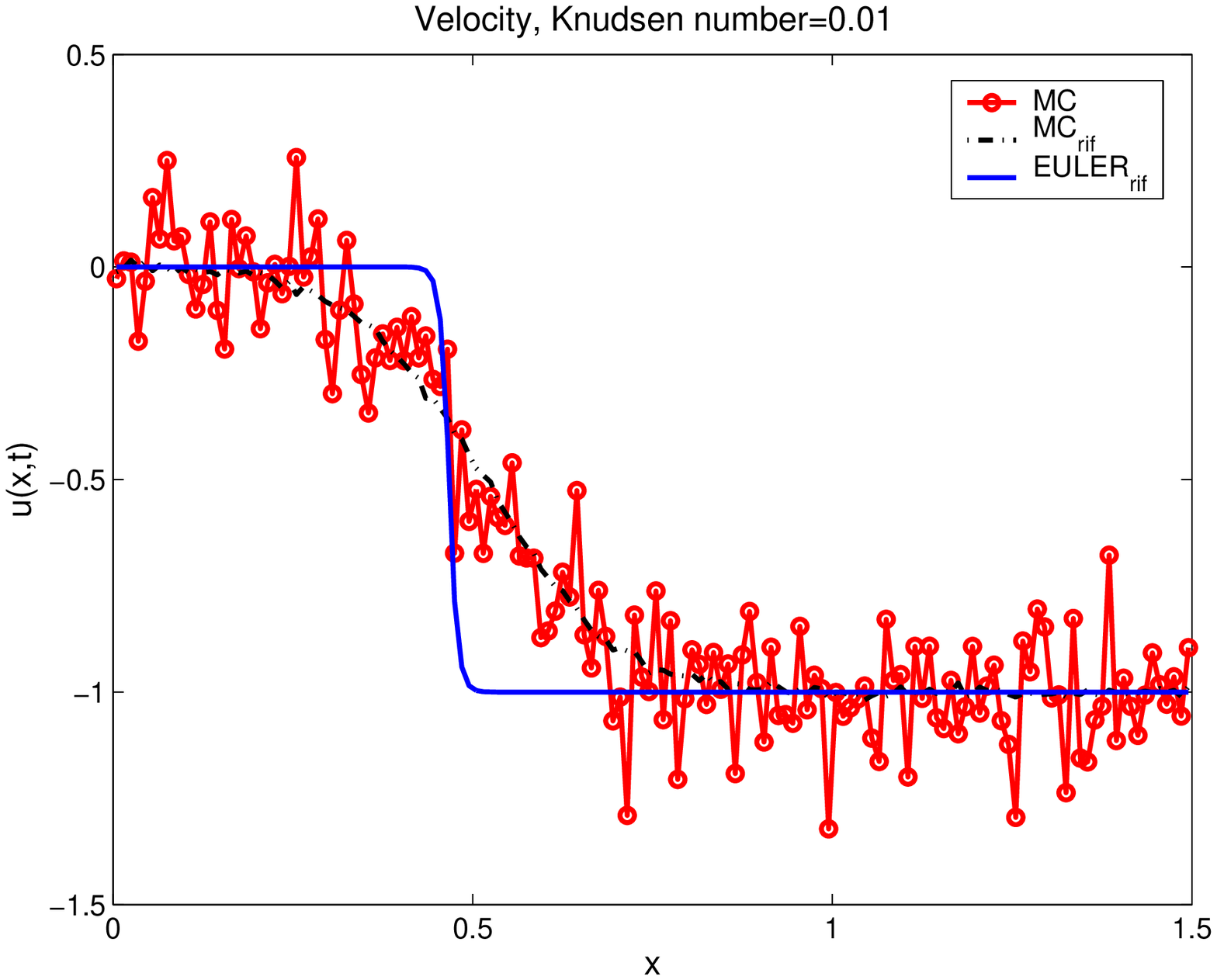}
\includegraphics[scale=0.39]{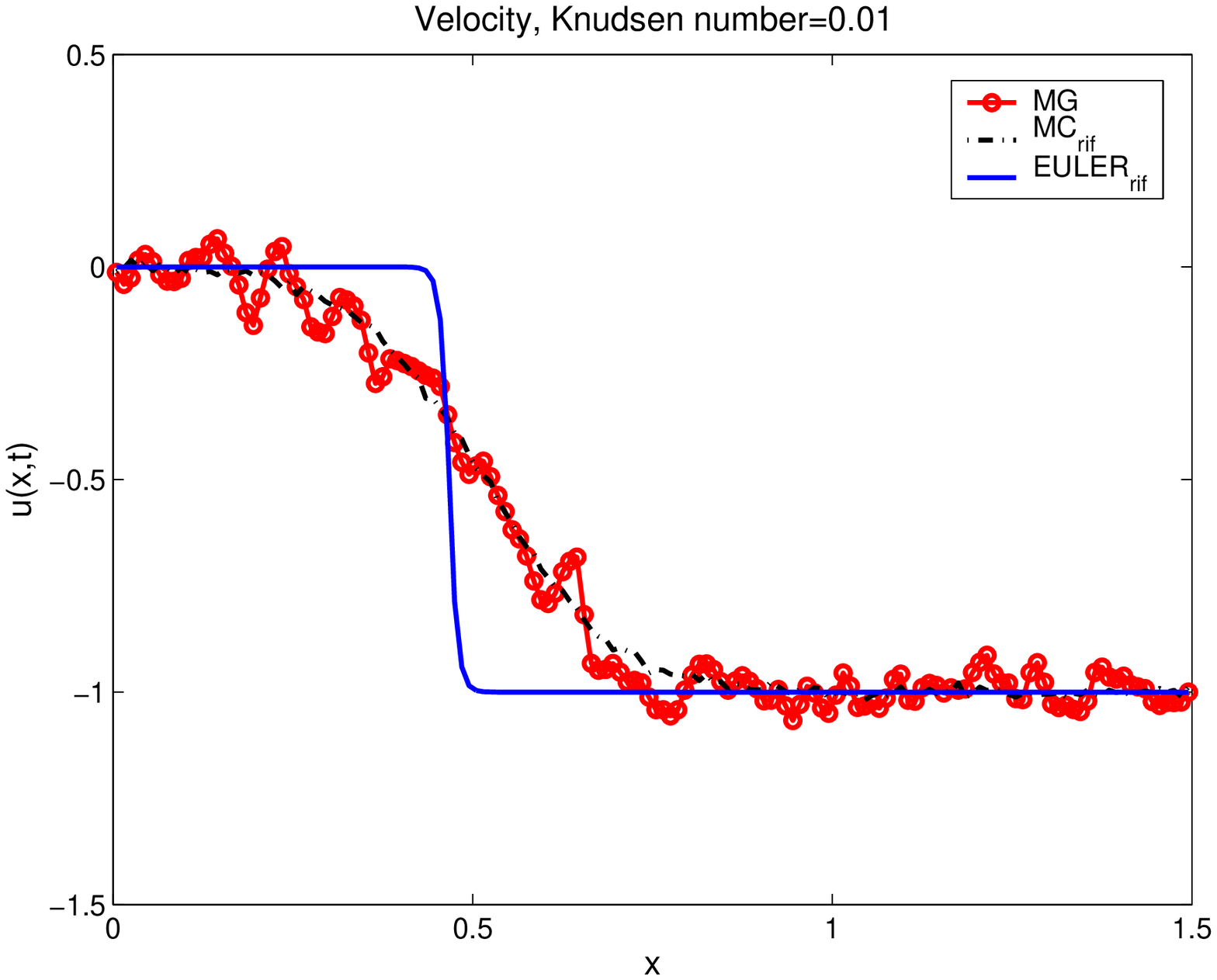}
\includegraphics[scale=0.39]{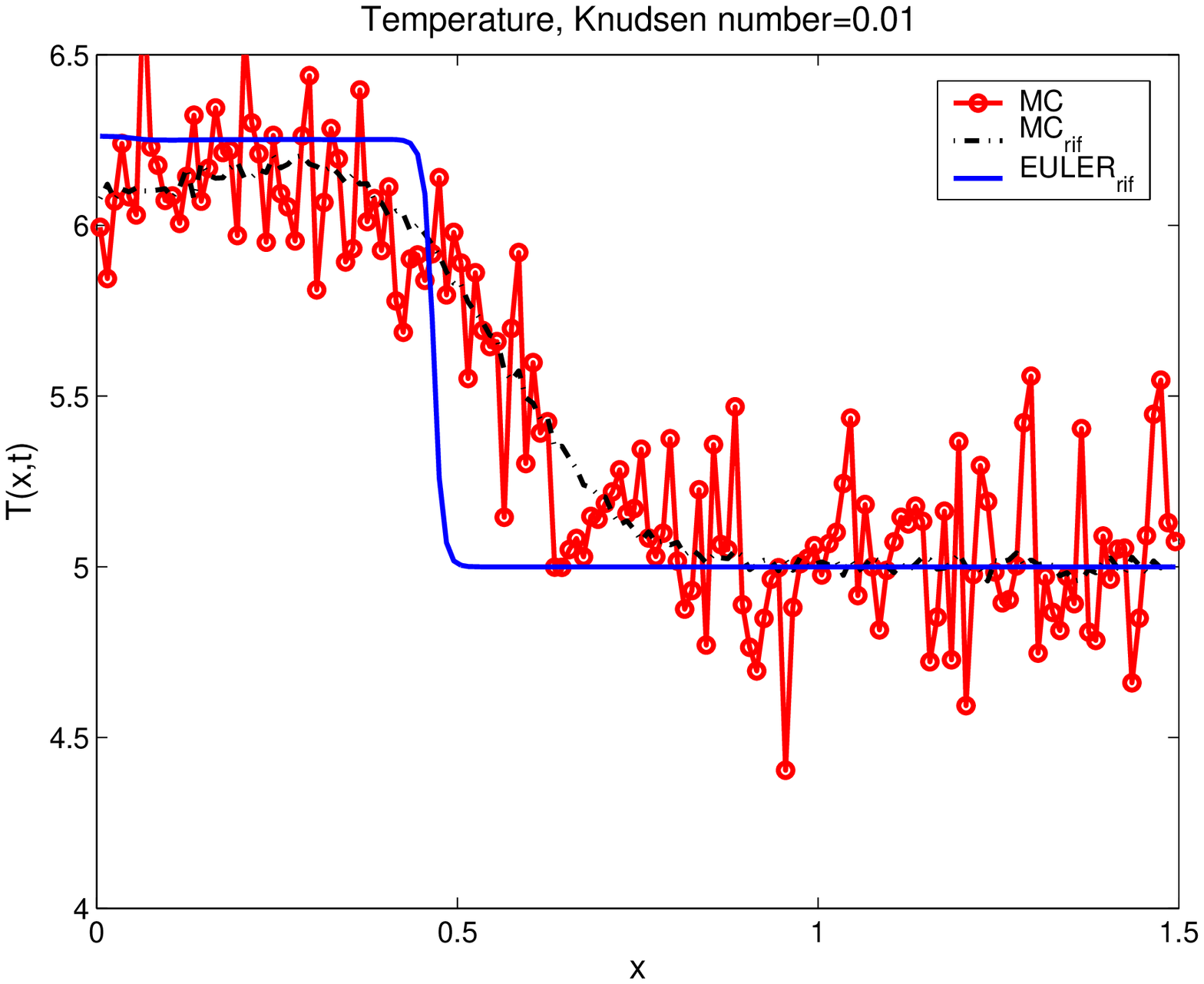}
\includegraphics[scale=0.39]{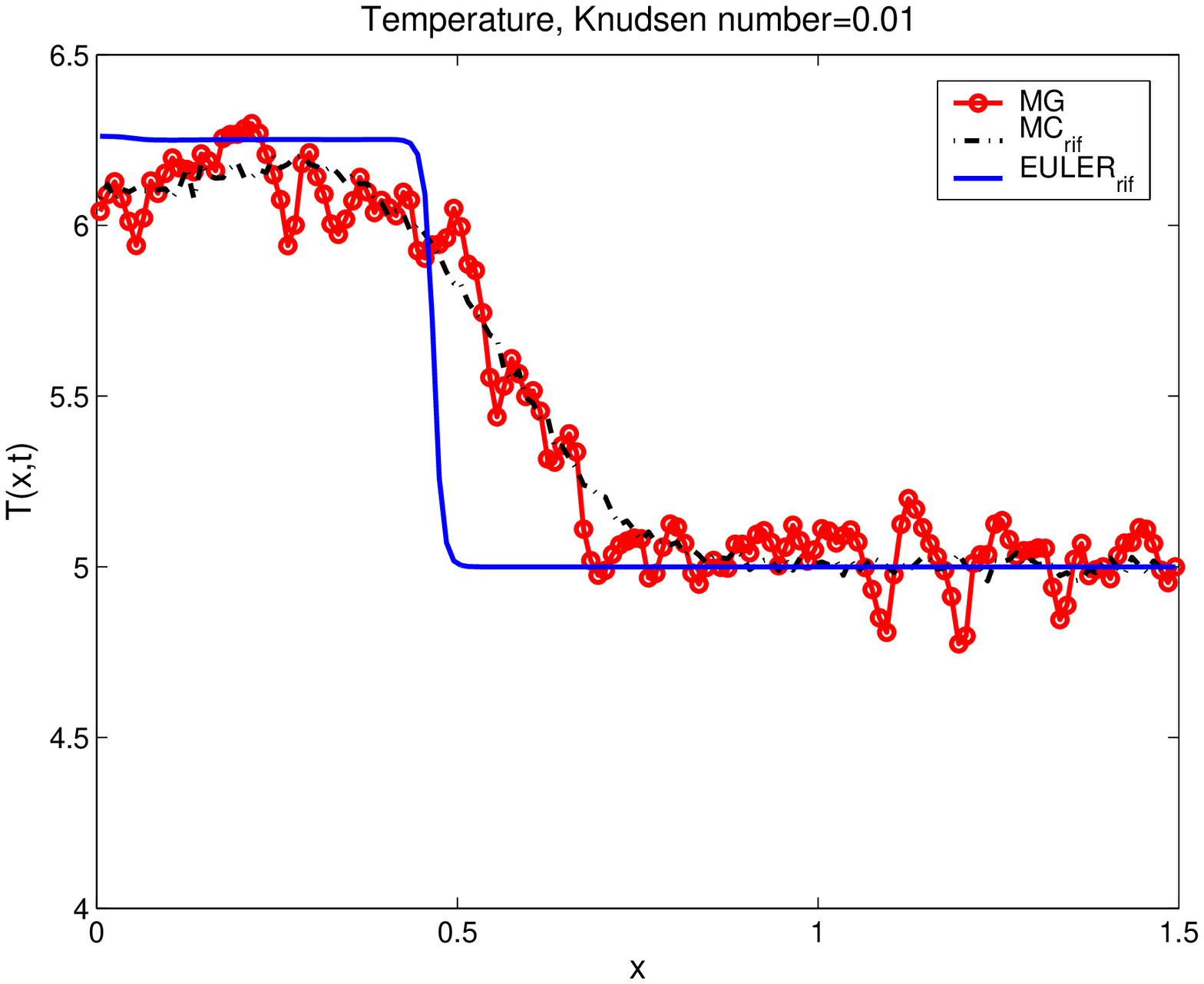}
\caption{Unsteady shock test: Solution at $t=0.18$ for the density
(top), velocity (middle) and temperature (bottom). MC method (left),
Moment Guided MG method (right). Knudsen number
$\varepsilon=10^{-2}$. Reference solution: dash dotted line. Euler
solution: continuous line. Monte Carlo or Moment Guided: circles
plus continuous line. 400 particles per cell.} \label{ST1}
\end{center}
\end{figure}

\begin{figure}
\begin{center}
\includegraphics[scale=0.39]{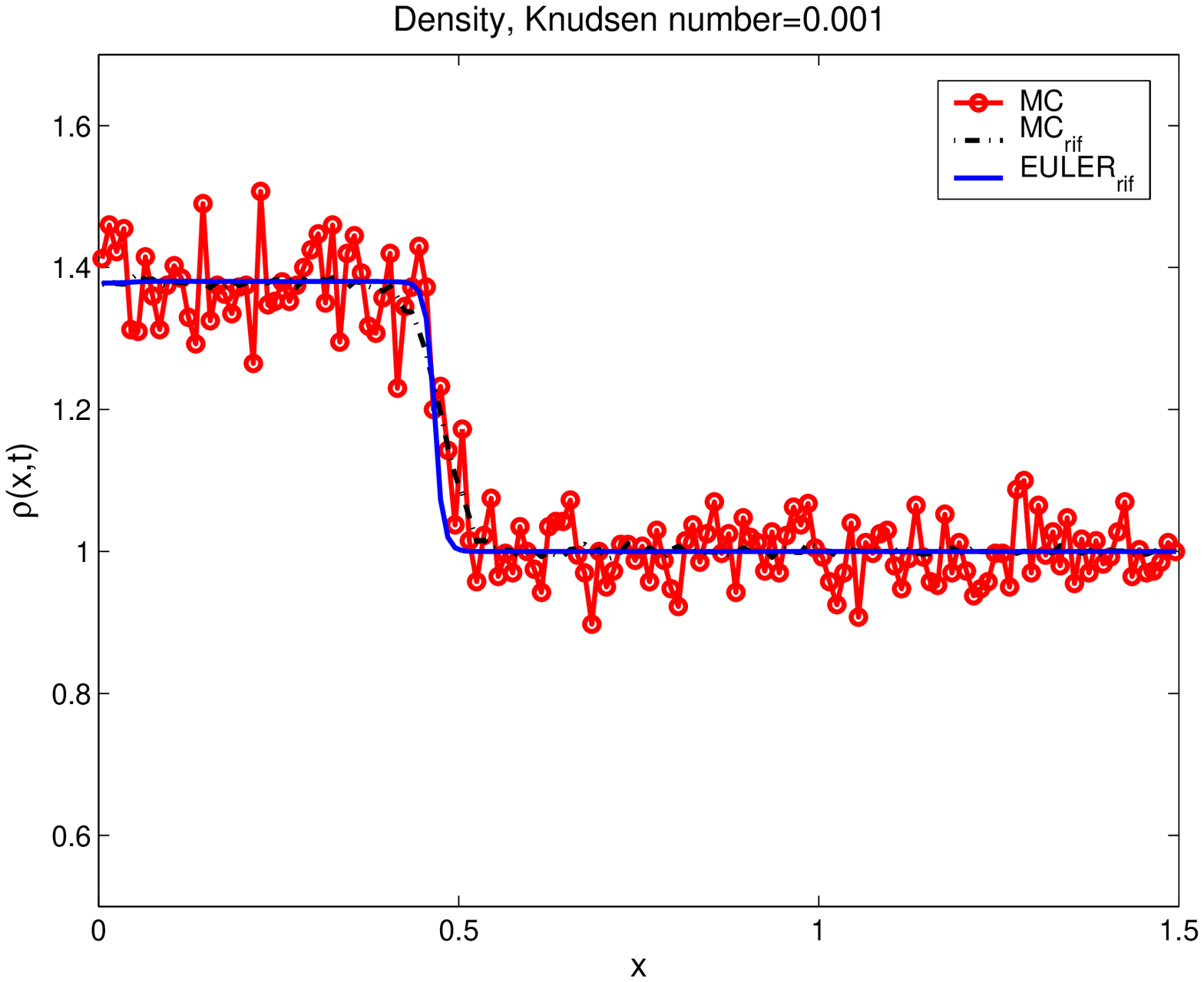}
\includegraphics[scale=0.39]{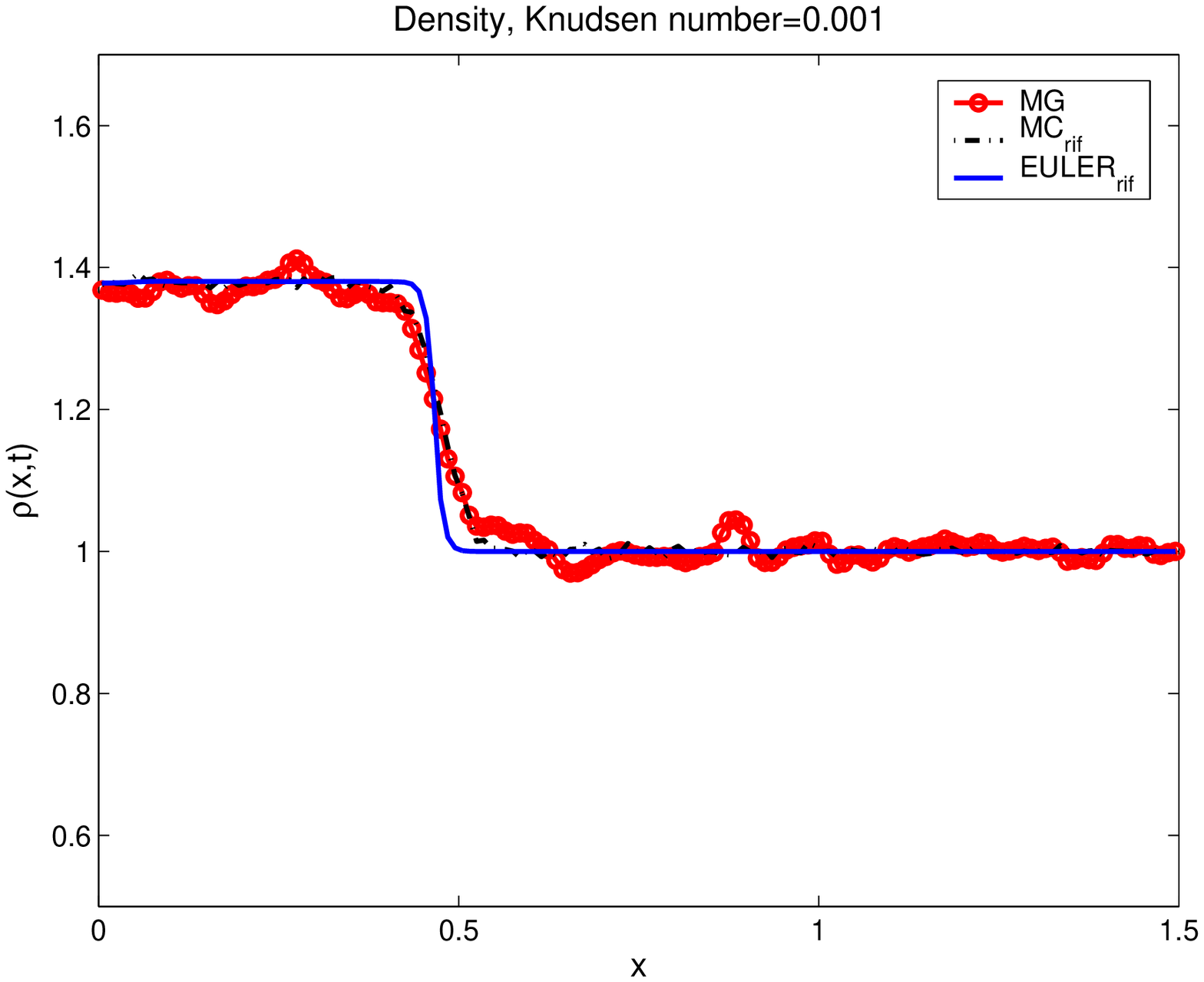}
\includegraphics[scale=0.39]{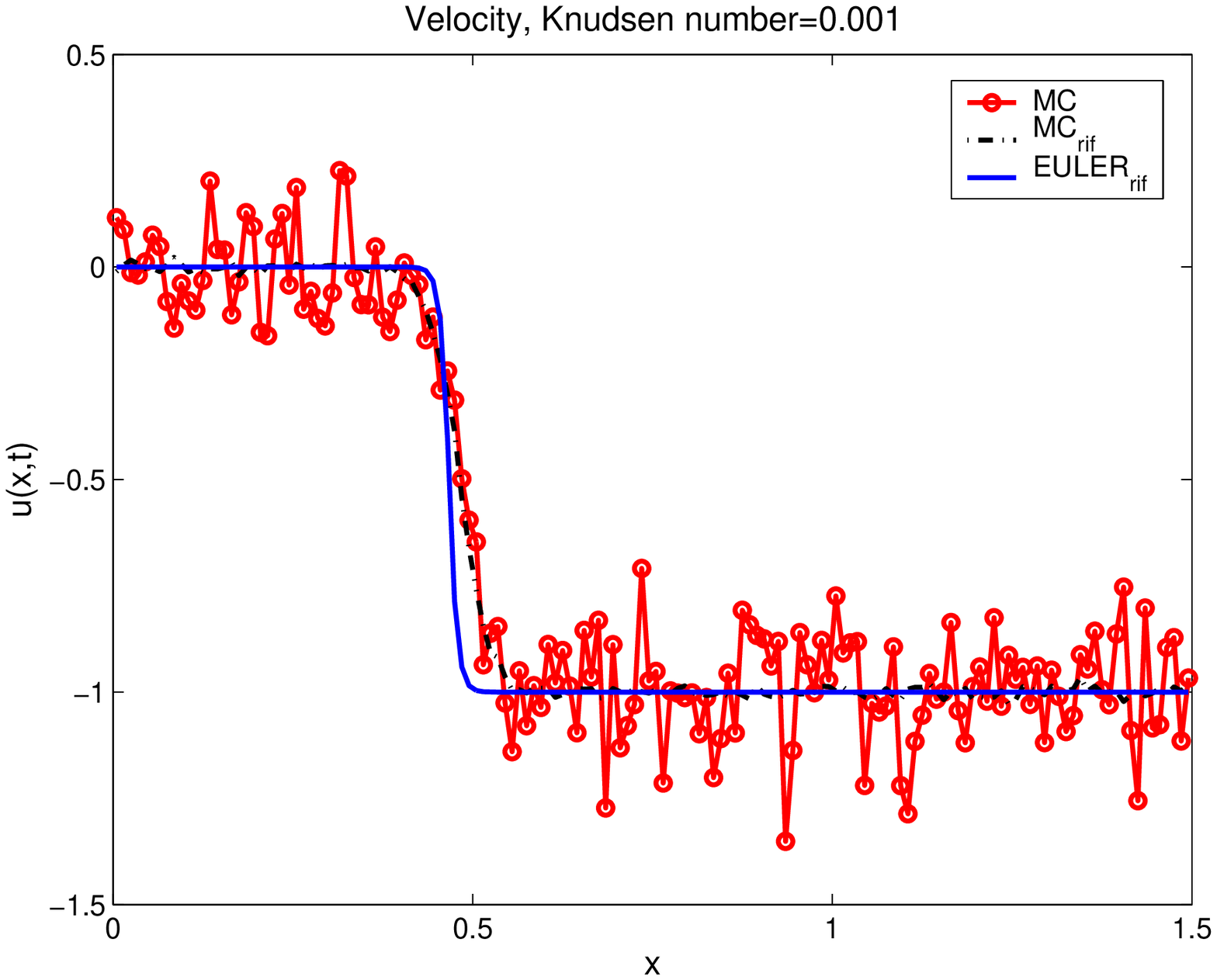}
\includegraphics[scale=0.39]{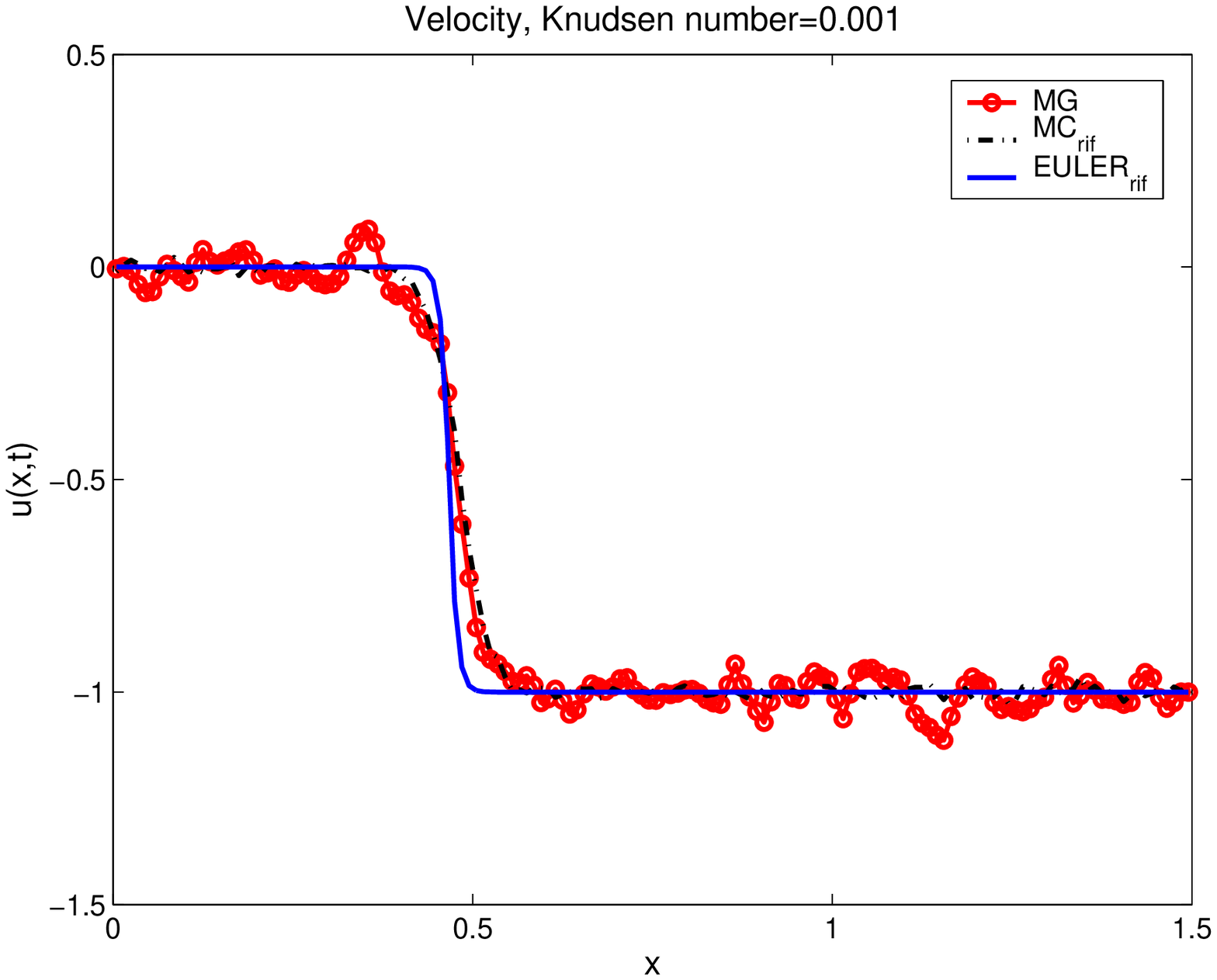}
\includegraphics[scale=0.39]{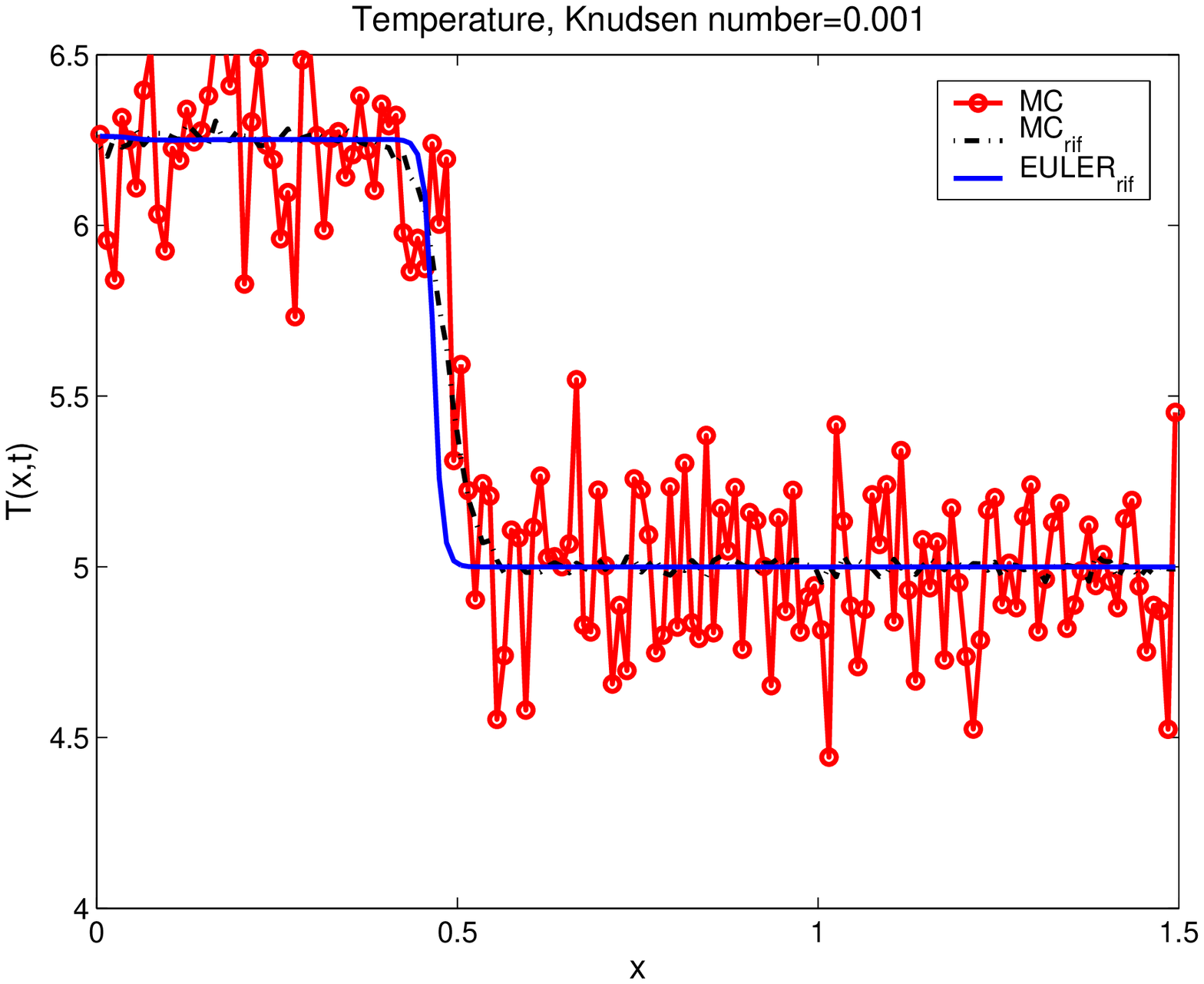}
\includegraphics[scale=0.39]{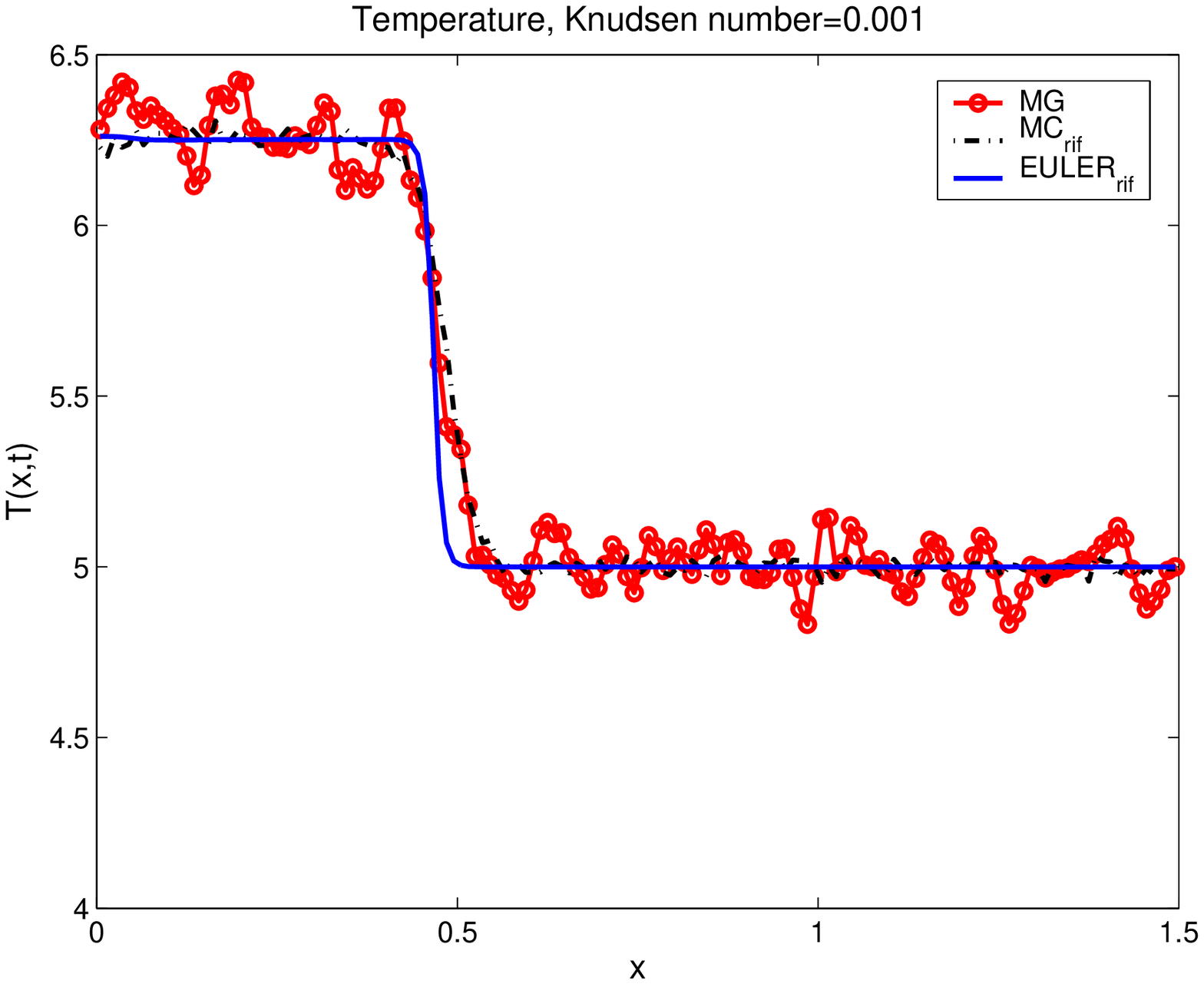}
\caption{Unsteady shock test: Solution at $t=0.18$ for the density
(top), velocity (middle) and temperature (bottom). MC method (left),
Moment Guided MG method (right). Knudsen number
$\varepsilon=10^{-3}$. Reference solution: dash dotted line. Euler
solution: continuous line. Monte Carlo or Moment Guided: circles
plus continuous line. 400 particles per cell.} \label{ST2}
\end{center}
\end{figure}

\begin{figure}
\begin{center}
\includegraphics[scale=0.39]{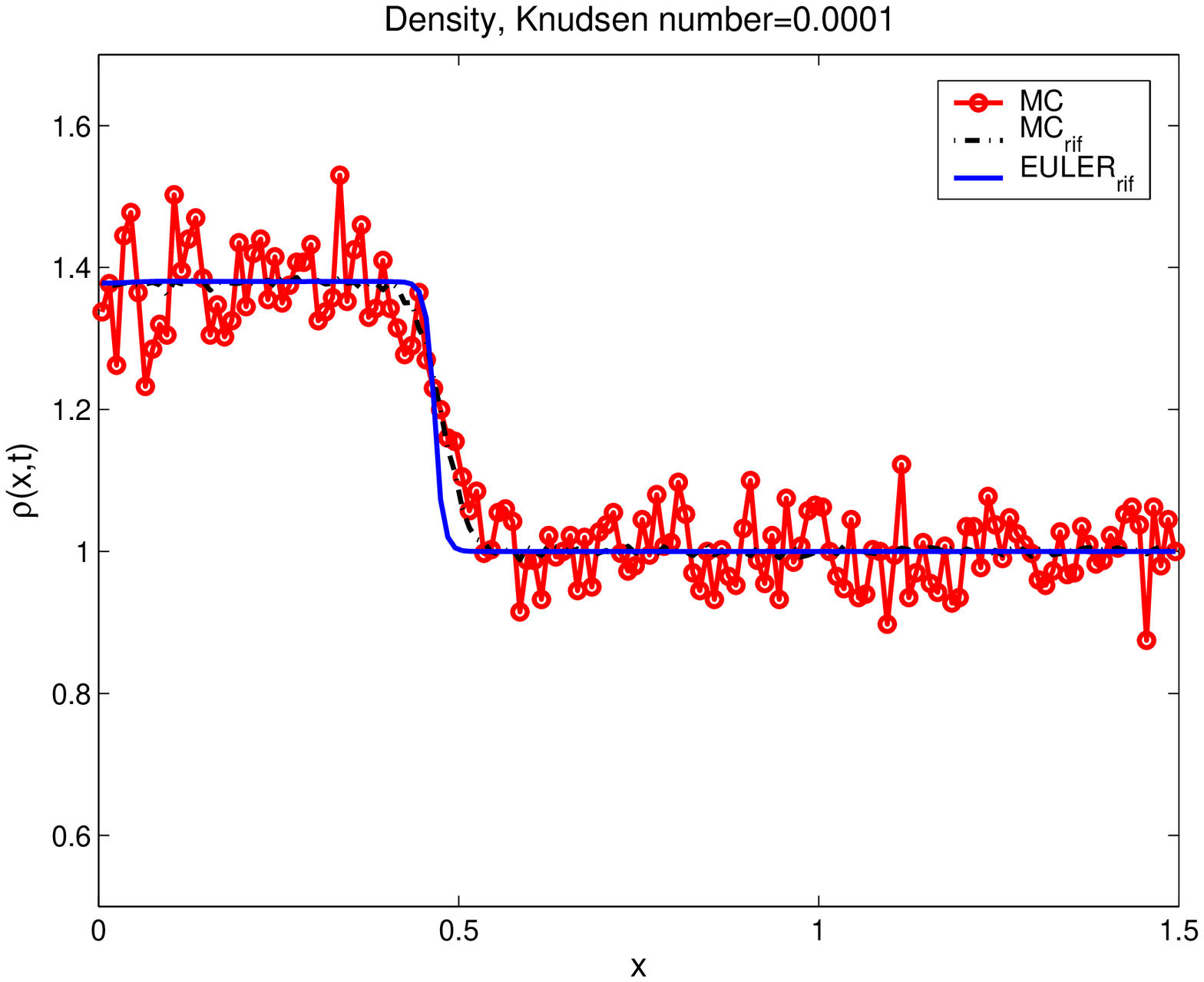}
\includegraphics[scale=0.39]{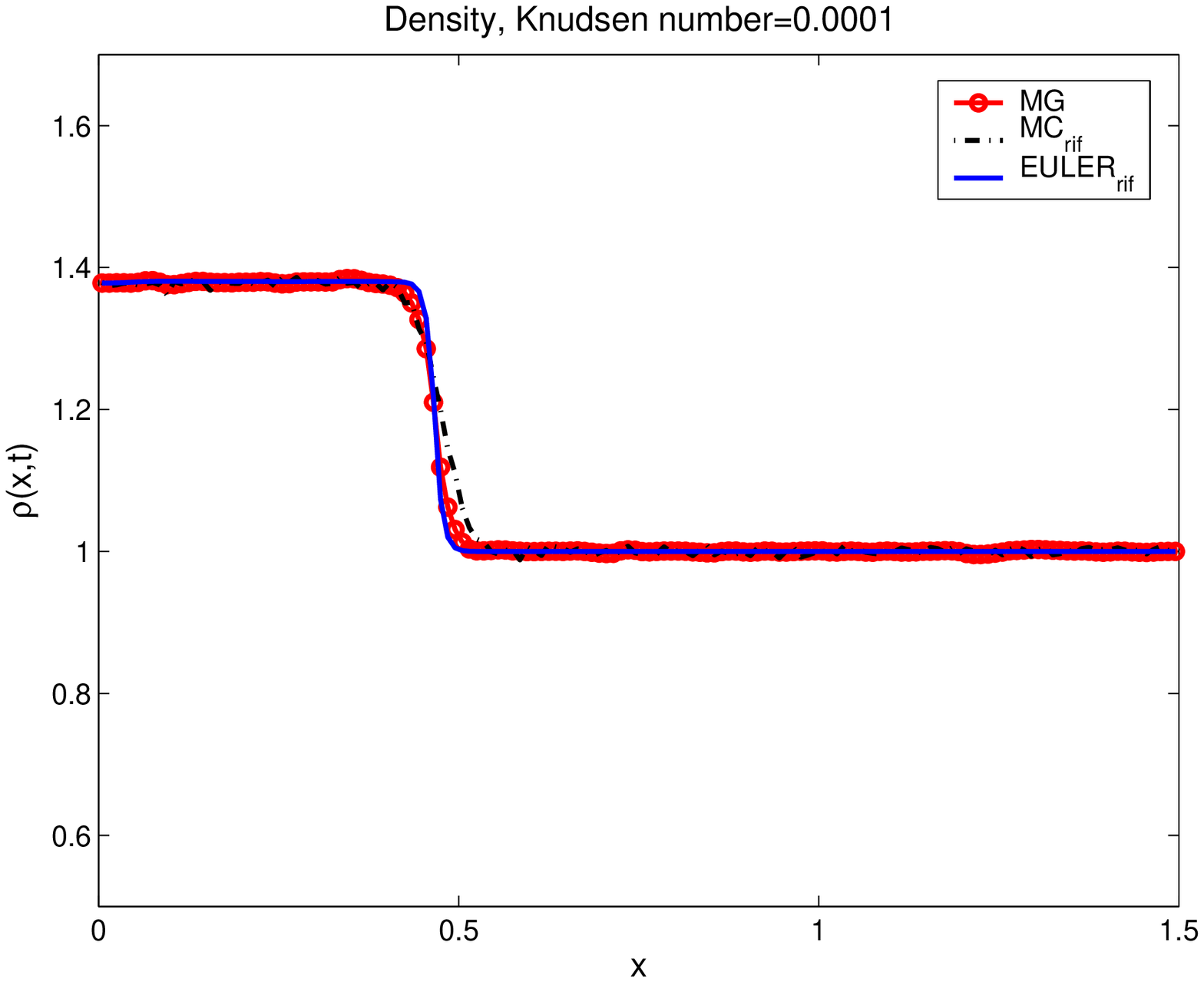}
\includegraphics[scale=0.39]{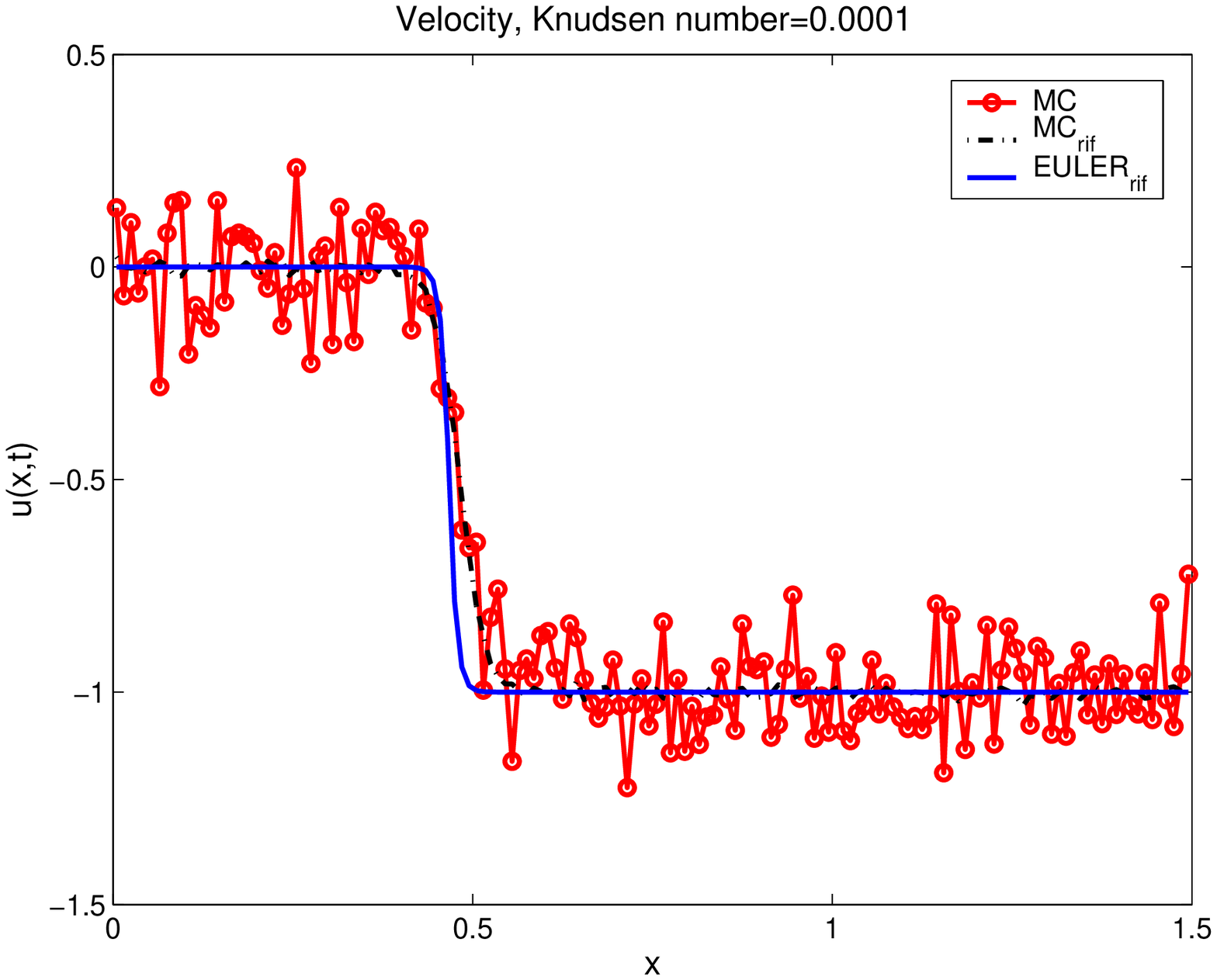}
\includegraphics[scale=0.39]{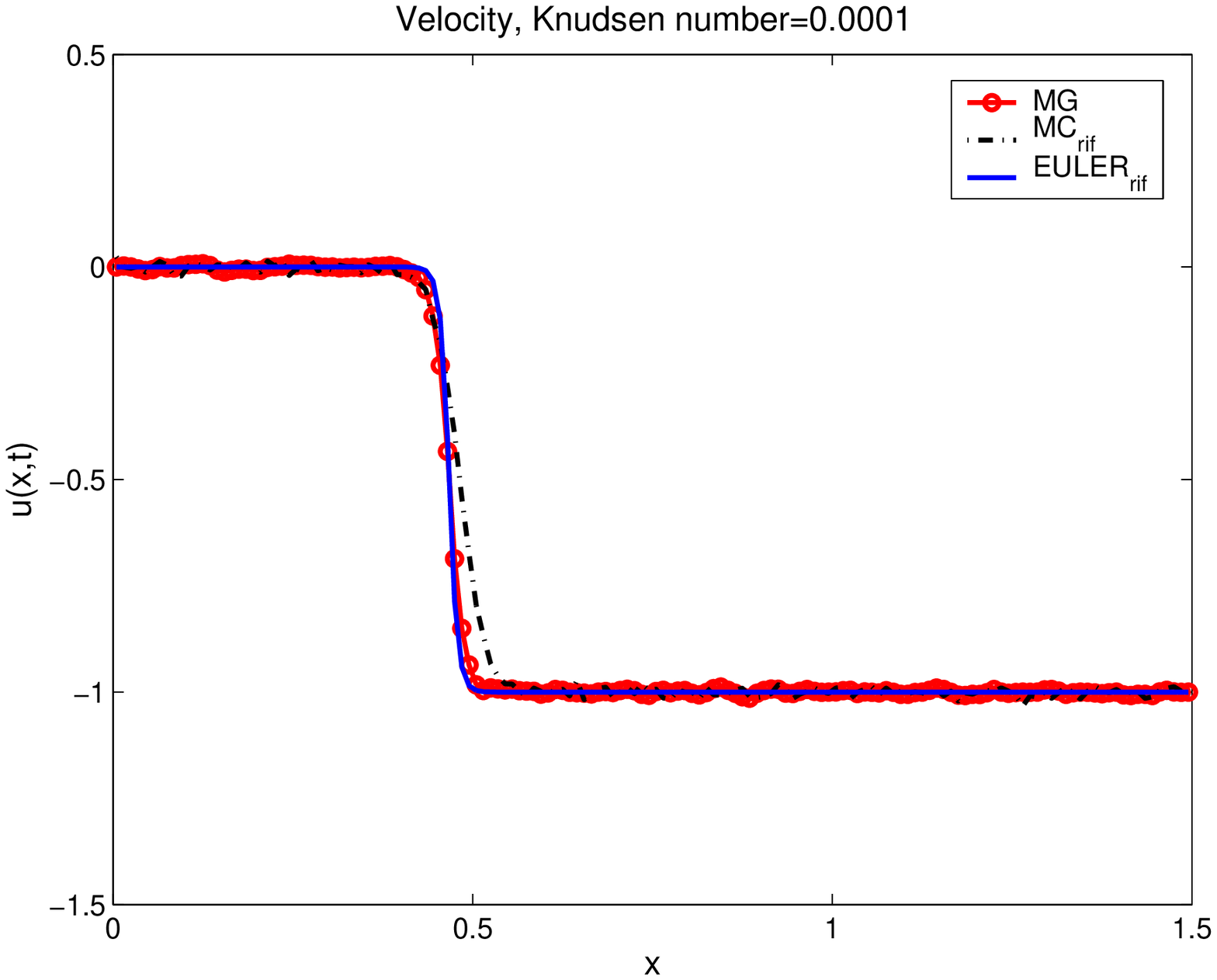}
\includegraphics[scale=0.39]{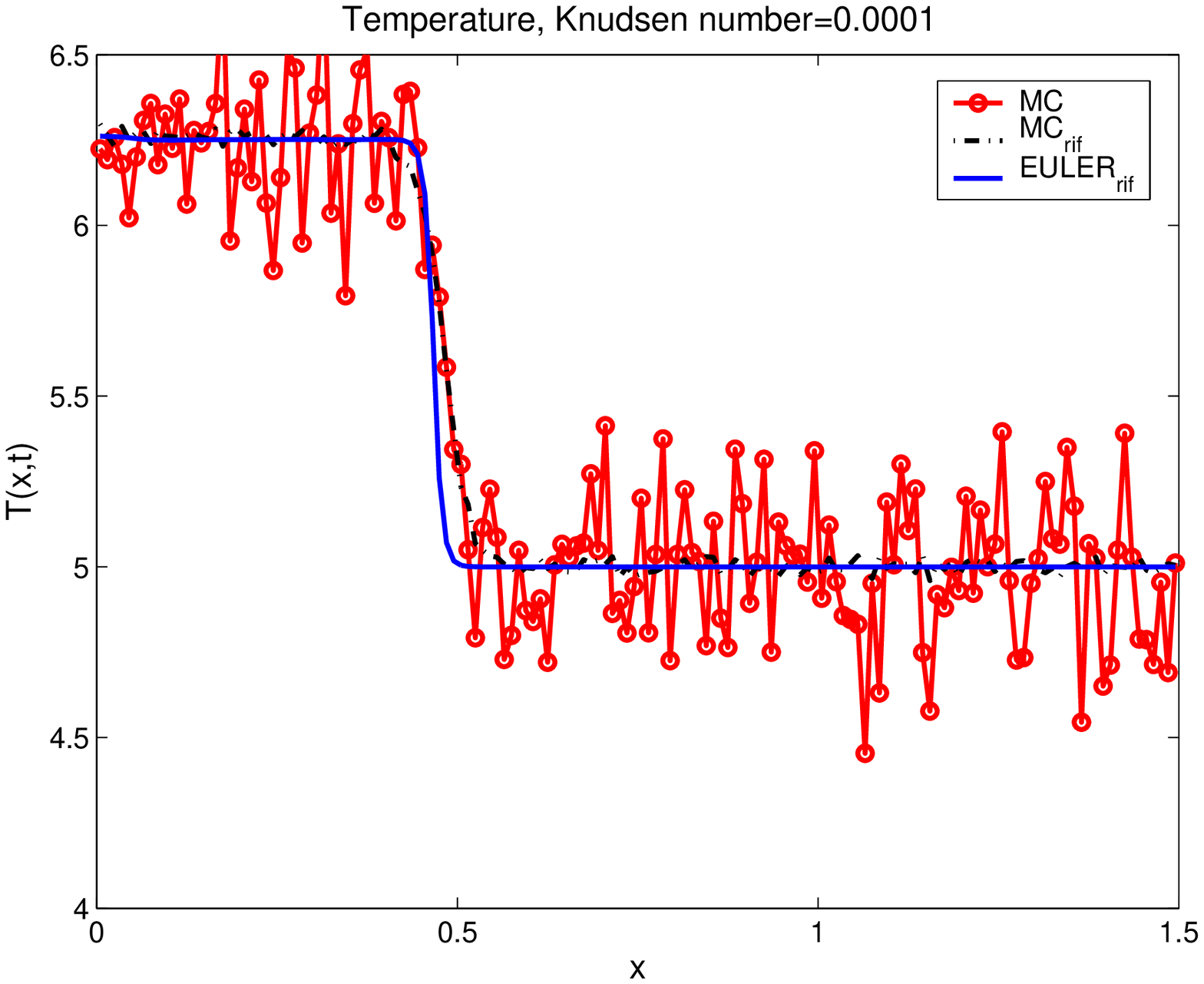}
\includegraphics[scale=0.39]{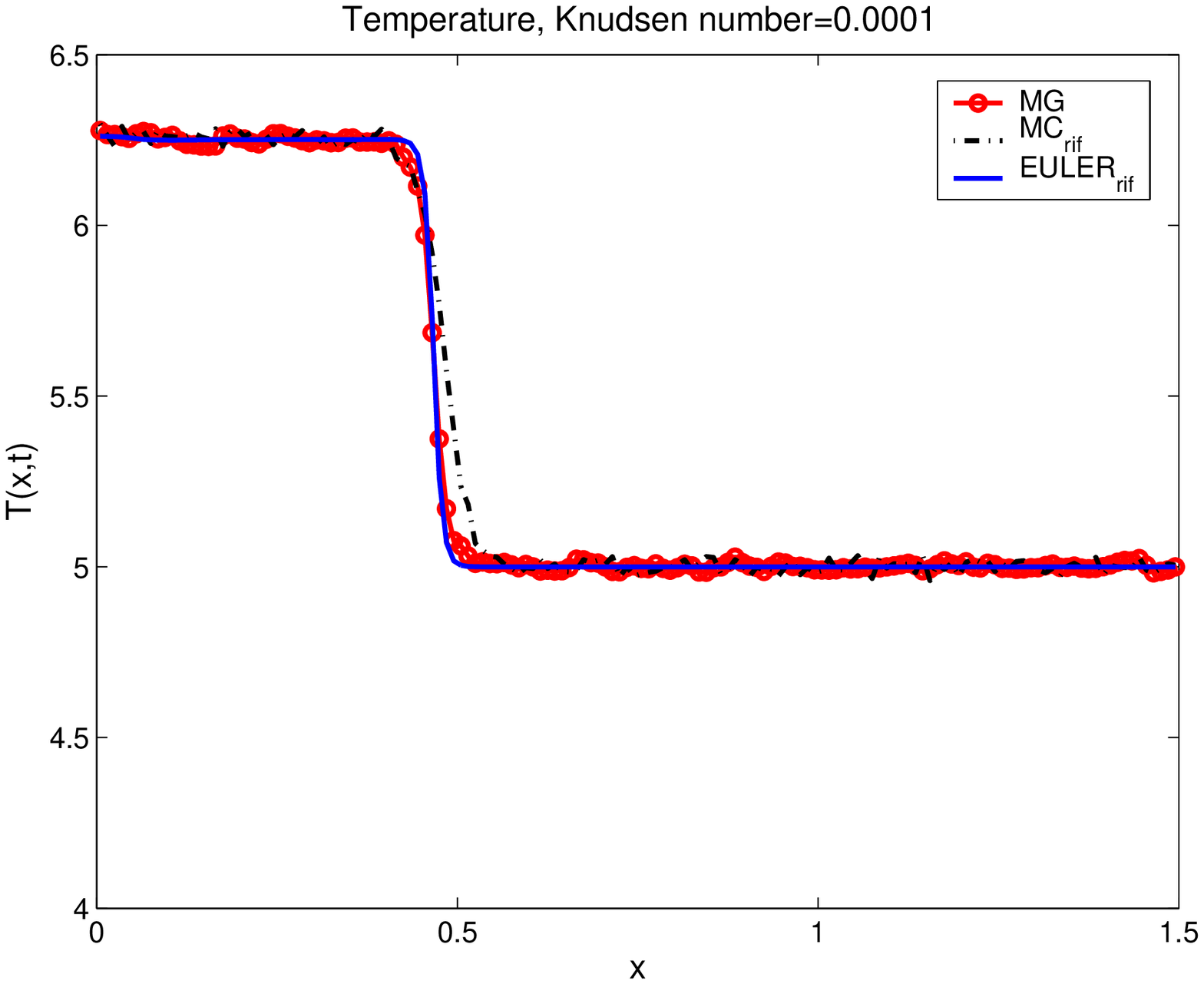}
\caption{Unsteady shock test: Solution at $t=0.18$ for the density
(top), velocity (middle) and temperature (bottom). MC method (left),
Moment Guided MG method (right). Knudsen number
$\varepsilon=10^{-4}$. Reference solution: dash dotted line. Euler
solution: continuous line. Monte Carlo or Moment Guided: circles
plus continuous line. 400 particles per cell.} \label{ST3}
\end{center}
\end{figure}

\section{Numerical results}

In this section we report some numerical results for the
method proposed on different test cases. First, we consider
two shock problems and then, we perform an
accuracy test using a smooth initial data with periodic boundary conditions. For the two shock tests, we compare the
moment guided (MG) solution with a Monte Carlo (MC)
solution which employs the first order exponential Runge-Kutta method together with a splitting technique between the transport and the collision parts. We report also in each figure the results for the limit compressible Euler equation using the same second order MUSCL scheme which has been used for the equilibrium part of the Moment Guided method. Finally, the reference solutions for the two shock test problems are obtained by using the same Monte Carlo method described above with the same number of cells, in which, however, the number of particles is such that the statistical noise is very small.

All the tests considered are relative to unsteady problems. We did this choice with the scope of highlighting
the fact that the method proposed in this paper is specifically designed for situations in
which the classical variance reduction techniques which employs time
averaging, which is the typical case of DSMC methods, cannot be used or turns out to be useless, since
time-averaging leads to the same computational effort of using more particles in one single simulation.

\subsection*{Unsteady shock test}

The first problem concerns the study of an unsteady shock.
The figures \ref{ST1} to \ref{ST3} consider the same initial data
for the density $\varrho=1$, the mean velocity $u=-1$ and the
temperature $T=1$ for different initial Knudsen number values,
ranging from $\varepsilon=10^{-2}$ to $\varepsilon=10^{-4}$. Specular reflection boundary conditions are used on the left side to produce the shock wave. The right boundary condition prescribes inlet flow. The
number of cells is $150$ while the time step is given by the minimum of the ratio of $\Delta x$ over the maximum
velocity owned by the particles and the ratio of $\Delta x$ over the
larger eigenvalues of the compressible Euler system. The final time is $0.18$. The Knudsen number value does not play any role in the choice of the time step, being the method asymptotic preserving and thus independent from the small scale constraint prescribed by $\varepsilon$. Each figure
depicts the density, the mean velocity and the temperature from top
to bottom, with the Monte Carlo solver on the left and the Moment
Guided method on the right. At the beginning of the simulation $400$ particles per cell are used for both methods Monte Carlo and Moment Guided, the solutions intentionally still contain some fluctuations. This is done to clearly show the difference in term of statistical error between a Monte Carlo method and our method. In addition, we report the solutions of
the compressible Euler equations and the reference solution computed by the same Monte Carlo method where $5 \ 10^{5}$ particles per cell are employed.

The three figures \ref{ST1}-\ref{ST3} show a large reduction of fluctuations for all cases analyzed. Observe that the reduction of the fluctuations depends on the Knudsen number value. In particular in the case $\varepsilon=10^{-4}$ the solution is very close to the limit solution. Observe also that the solution furnished by the Moment Guided method when $\varepsilon=10^{-4}$ is closer to the shock with respect to the reference Monte Carlo solution. This can be interpreted as an over relaxation problem of the Moment Guided method. However, an in deep simulation analysis we did not report, shows that, the difference between the MC and the MG methods when $\varepsilon=10^{-4}$ is due to the more diffusive behavior of the Monte Carlo method. This is caused by the large time steps $\Delta t$ allowed by the exponential method. These large time steps introduce a numerical diffusion in the schemes which is more important for the Monte Carlo method and less important for the Moment Guided one. 

\begin{figure}
\begin{center}
\includegraphics[scale=0.39]{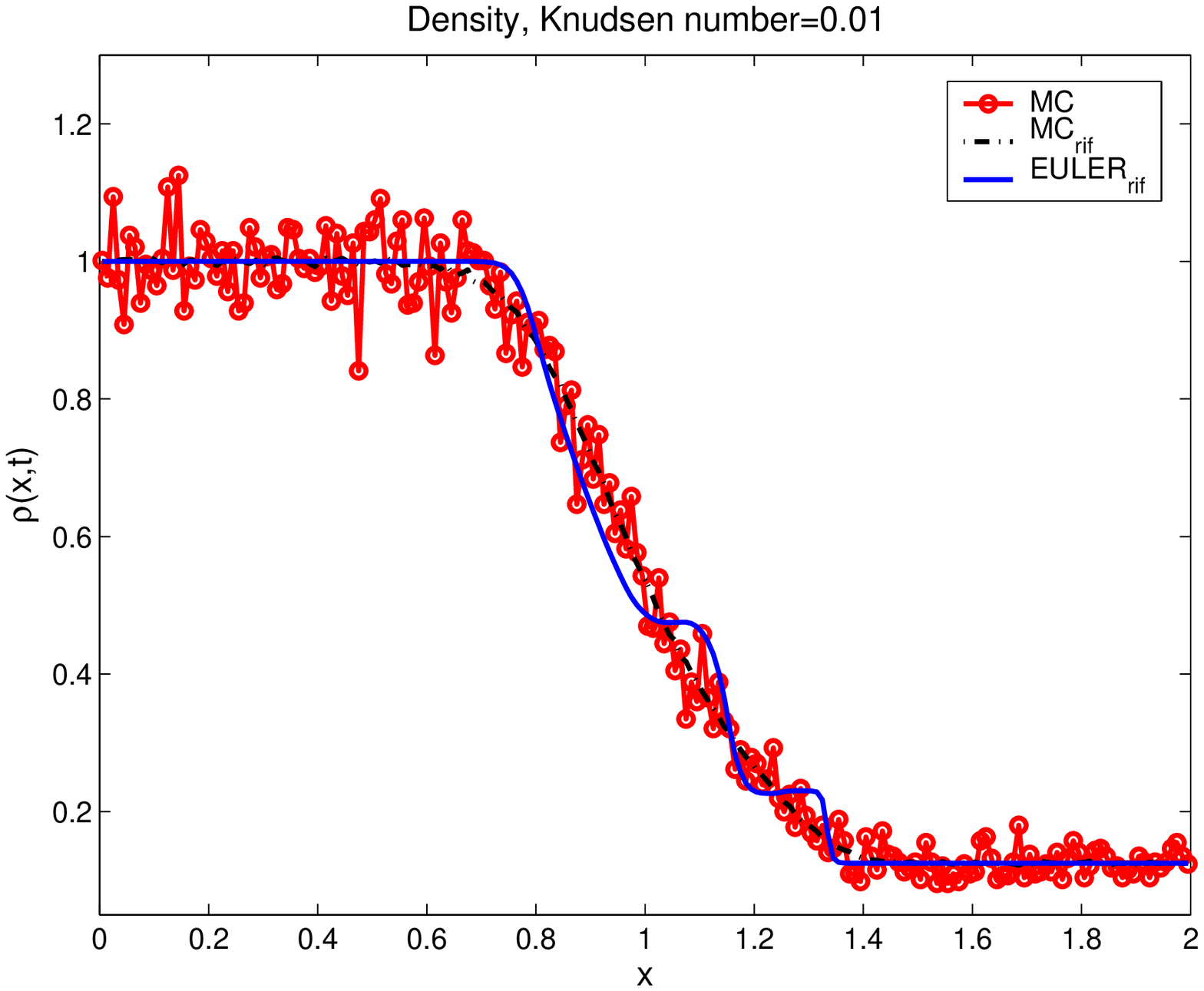}
\includegraphics[scale=0.39]{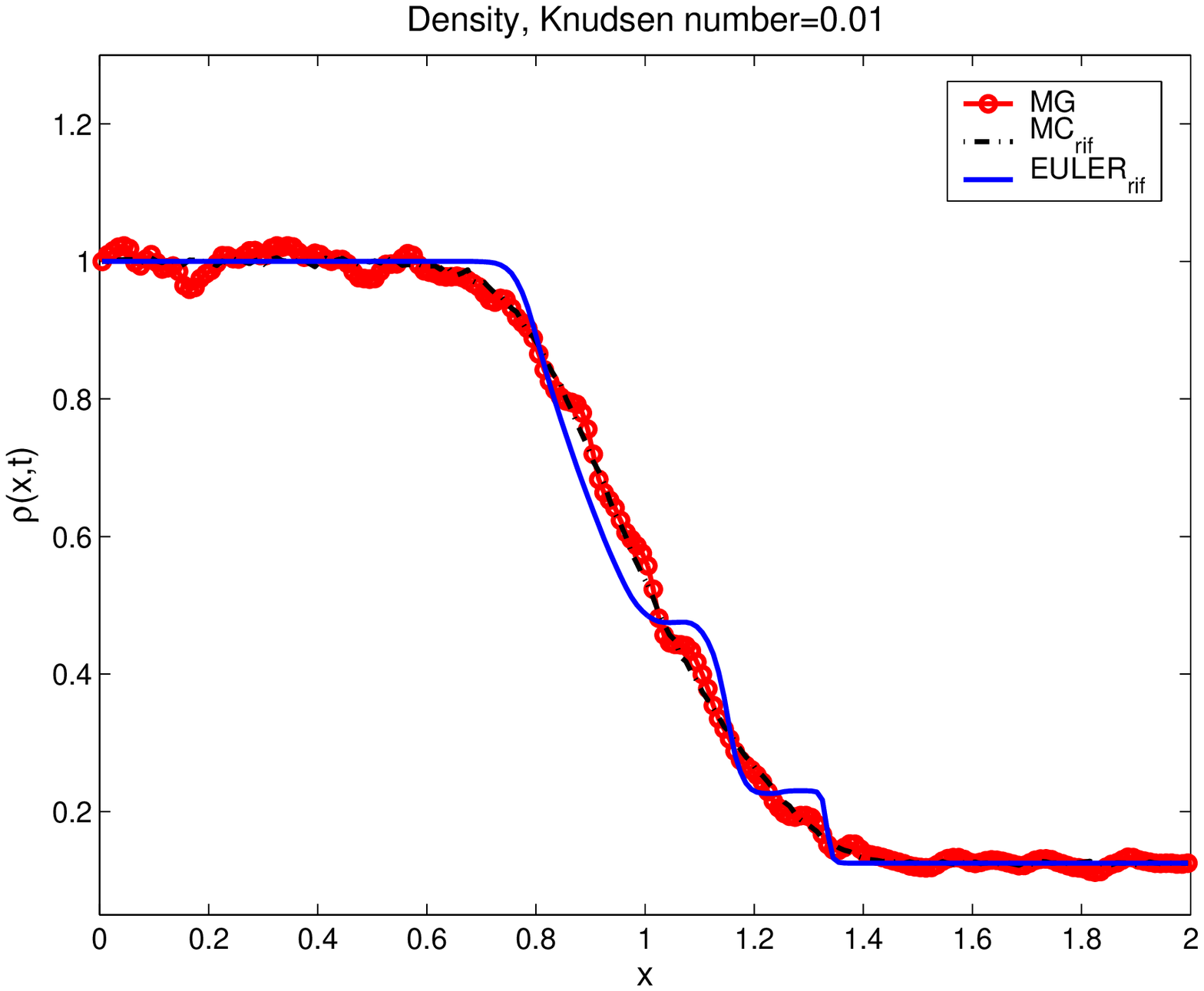}
\includegraphics[scale=0.39]{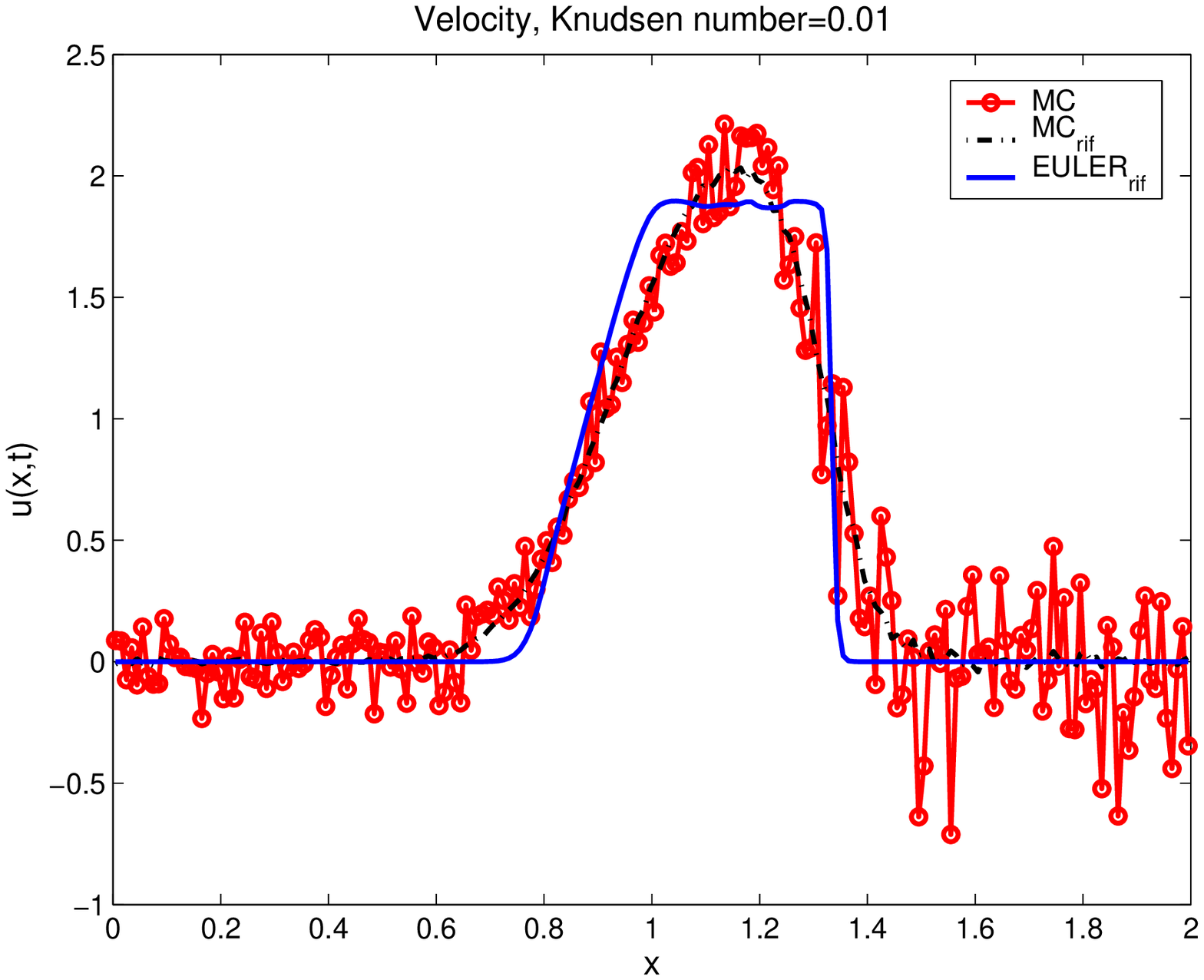}
\includegraphics[scale=0.39]{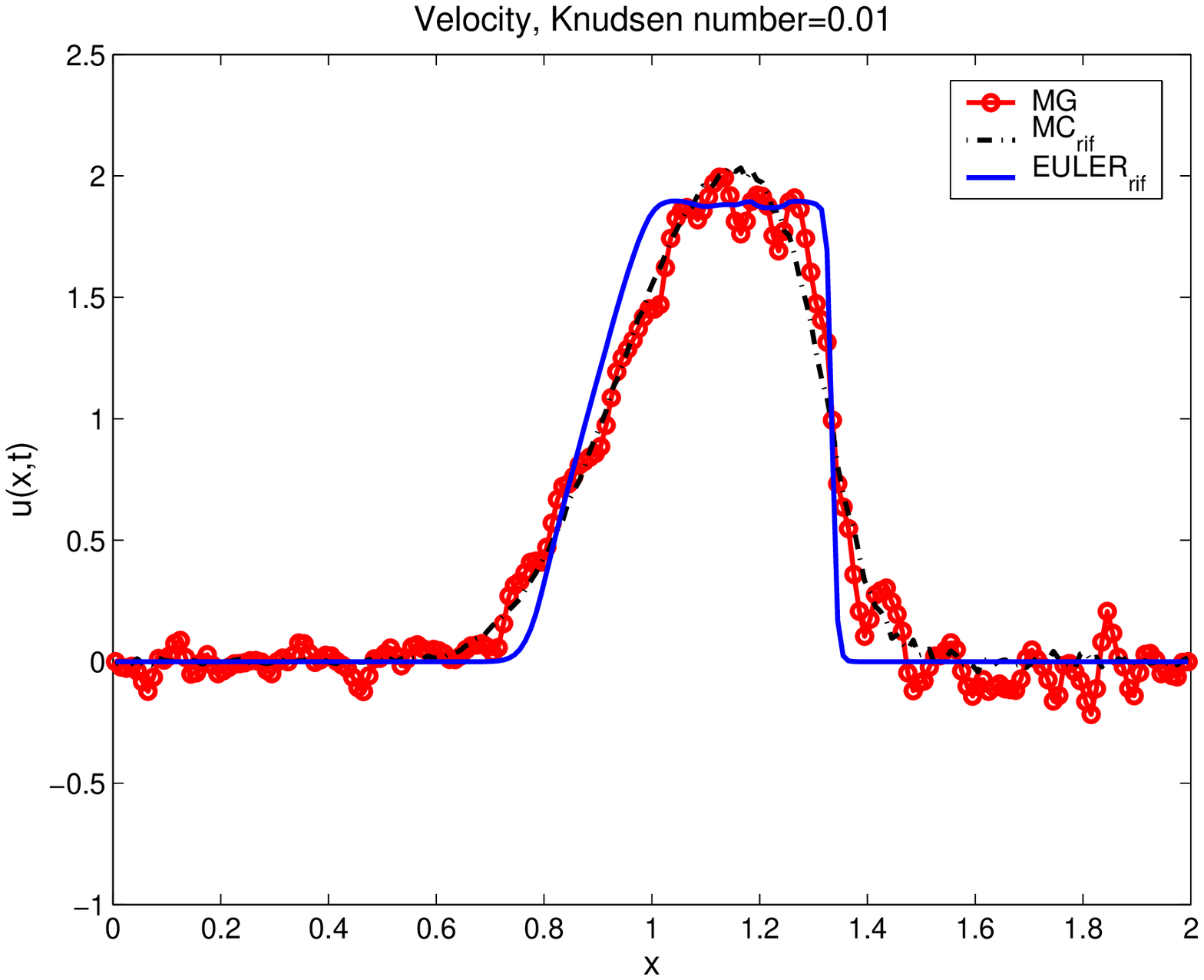}
\includegraphics[scale=0.39]{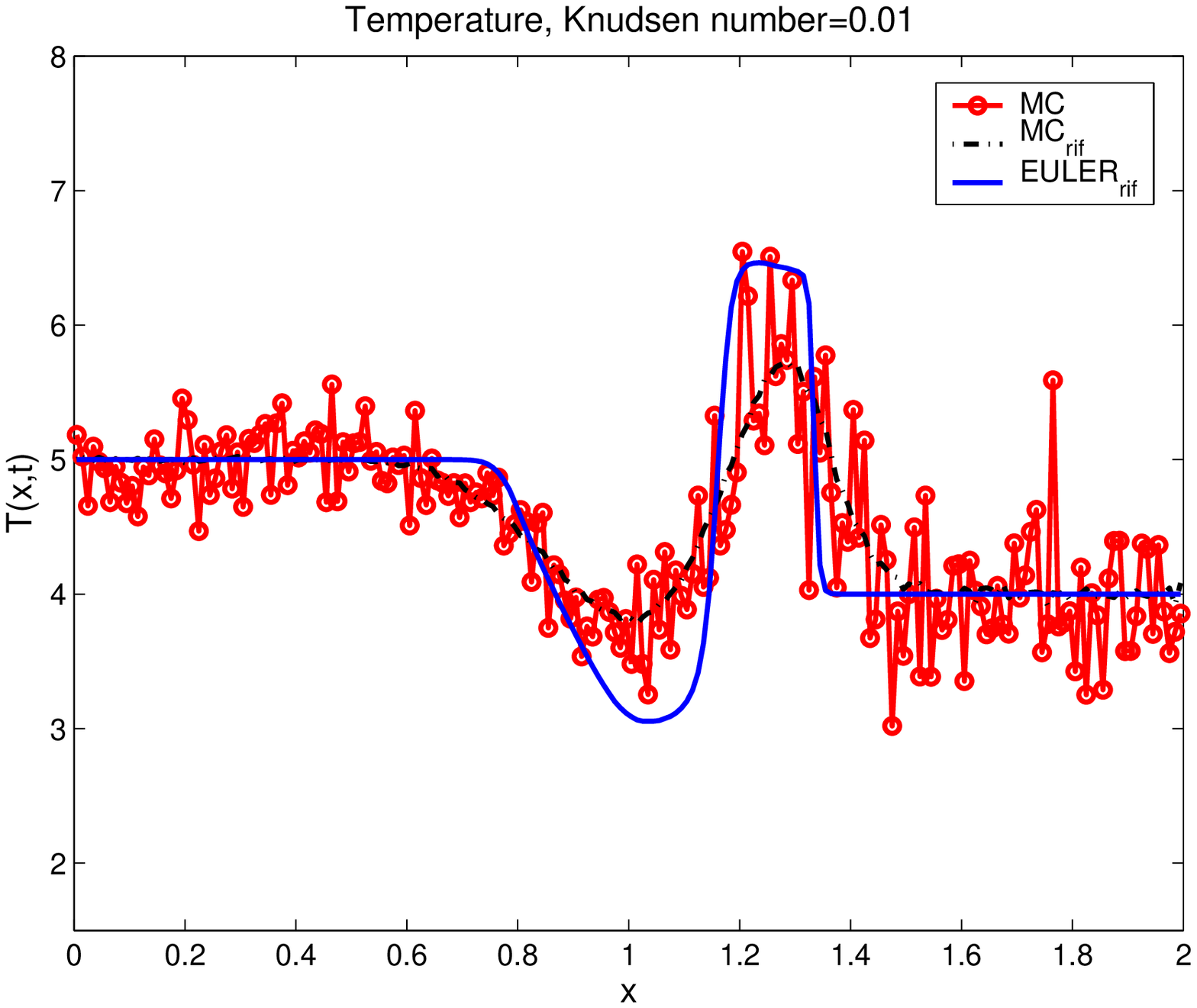}
\includegraphics[scale=0.39]{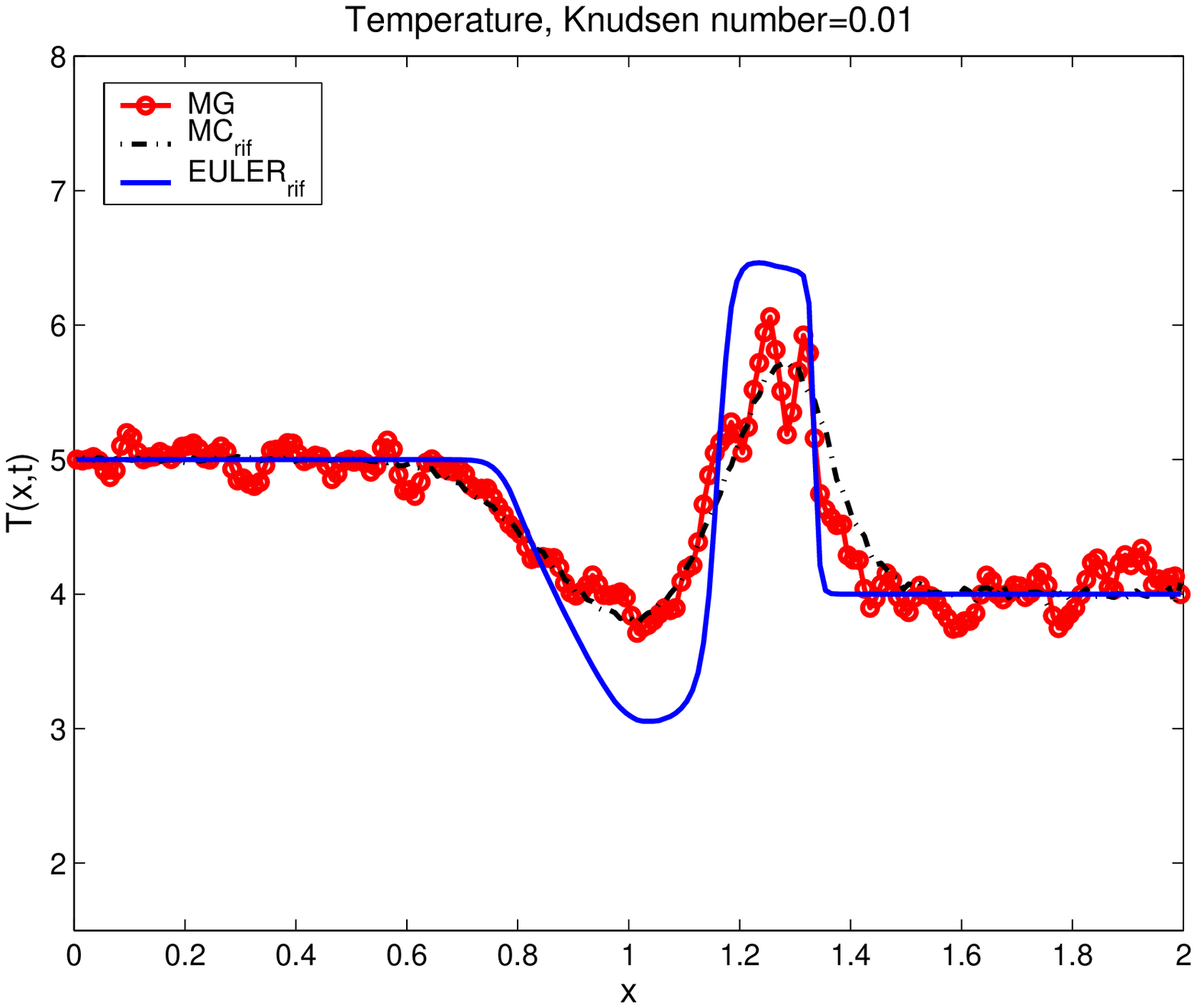}
\caption{Sod shock tube test: Solution at $t=0.08$ for the density
(top), velocity (middle) and temperature (bottom). MC method (left),
Moment Guided MG method (right). Knudsen number
$\varepsilon=10^{-2}$. Reference solution: dash dotted line. Euler
solution: continuous line. Monte Carlo or Moment Guided: circles
plus continuous line. 200 particles for cell.} \label{So1}
\end{center}
\end{figure}
\begin{figure}
\begin{center}
\includegraphics[scale=0.39]{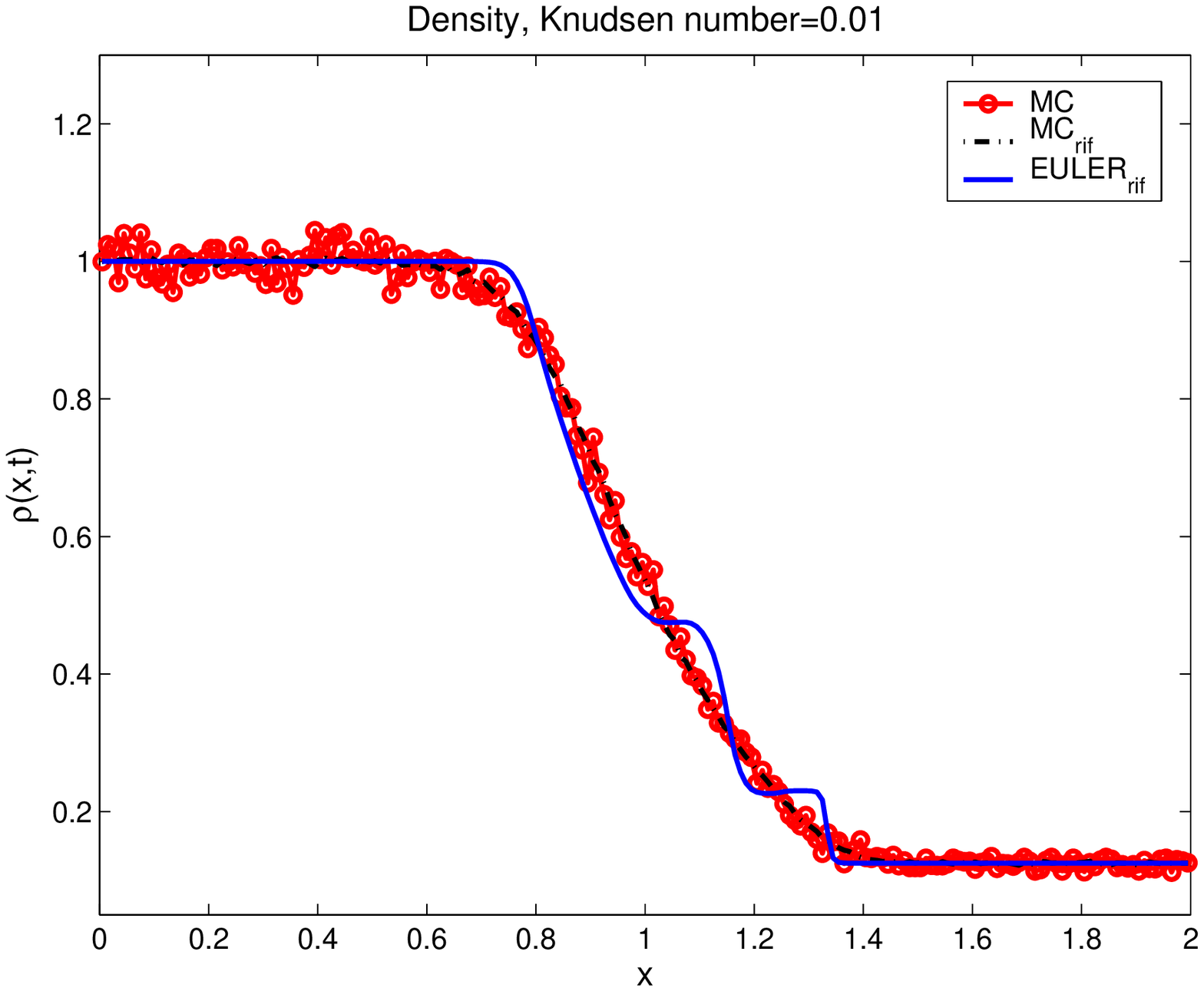}
\includegraphics[scale=0.39]{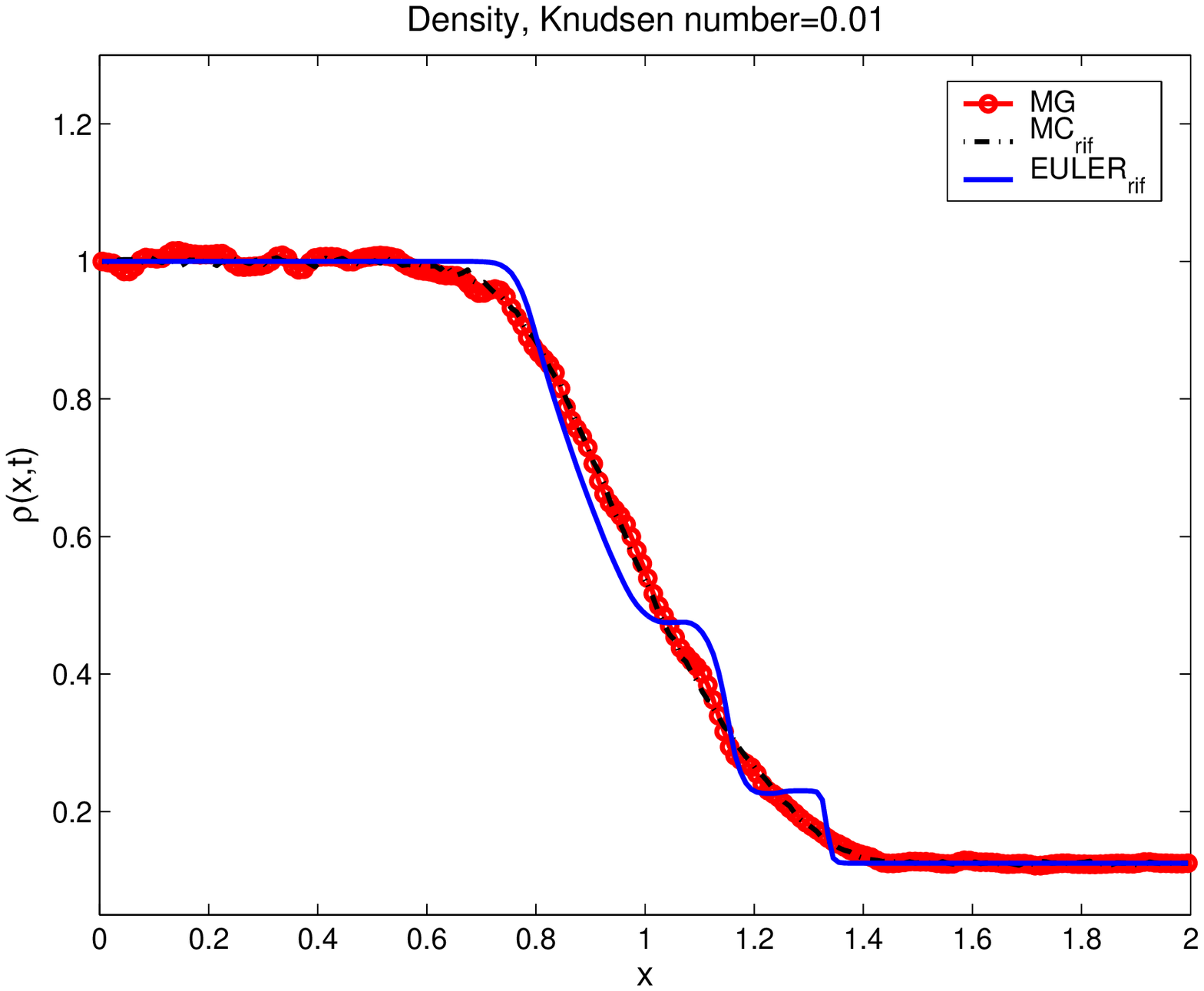}
\includegraphics[scale=0.39]{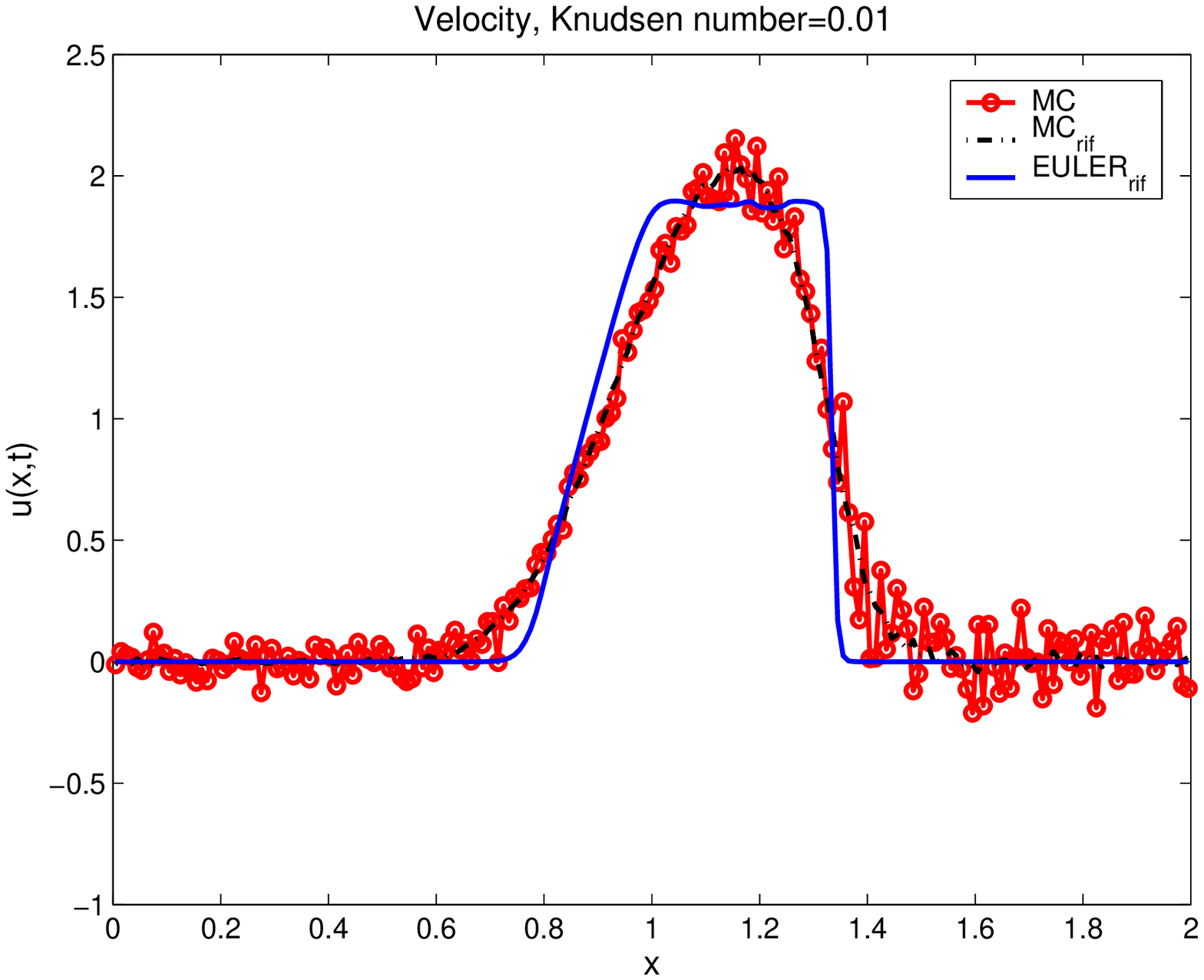}
\includegraphics[scale=0.39]{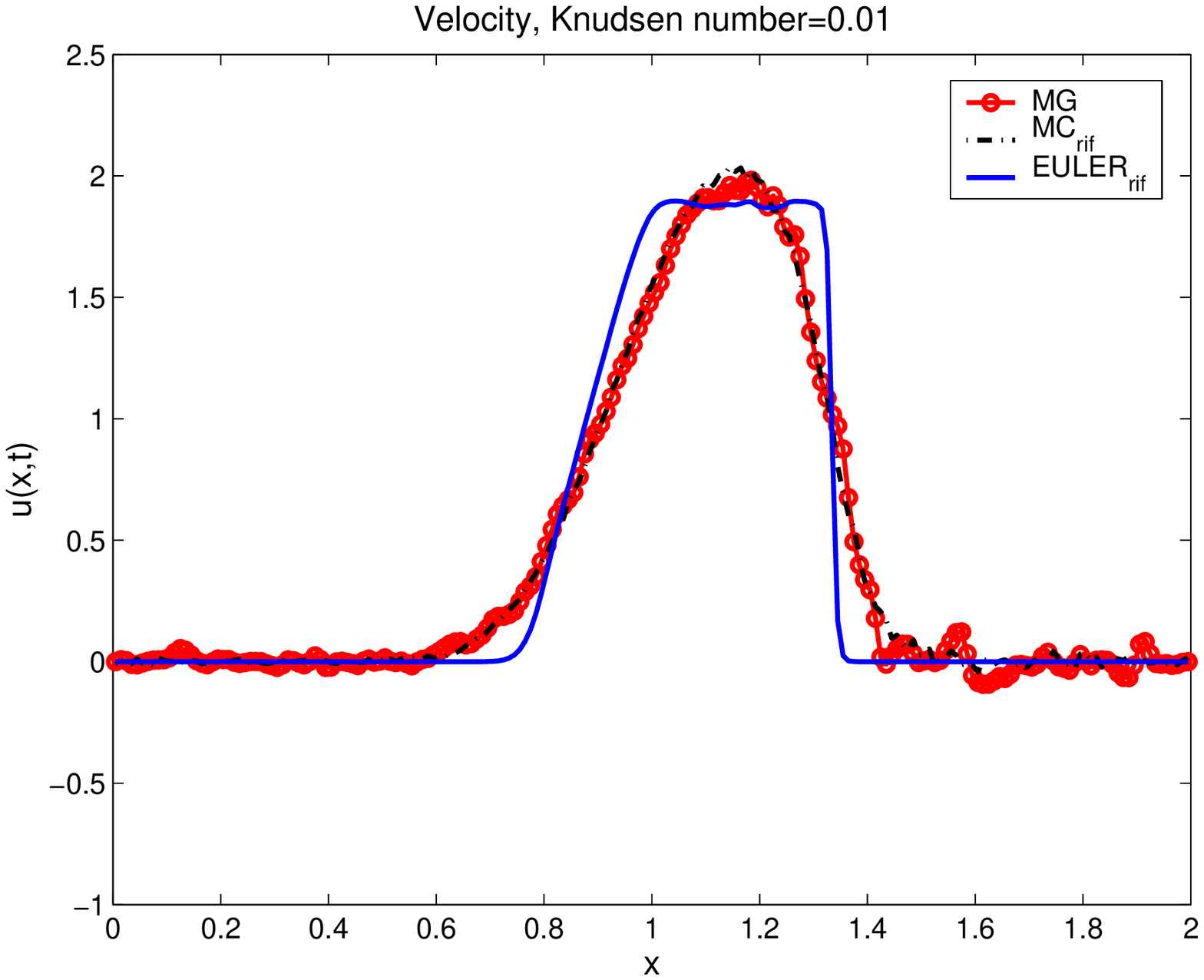}
\includegraphics[scale=0.39]{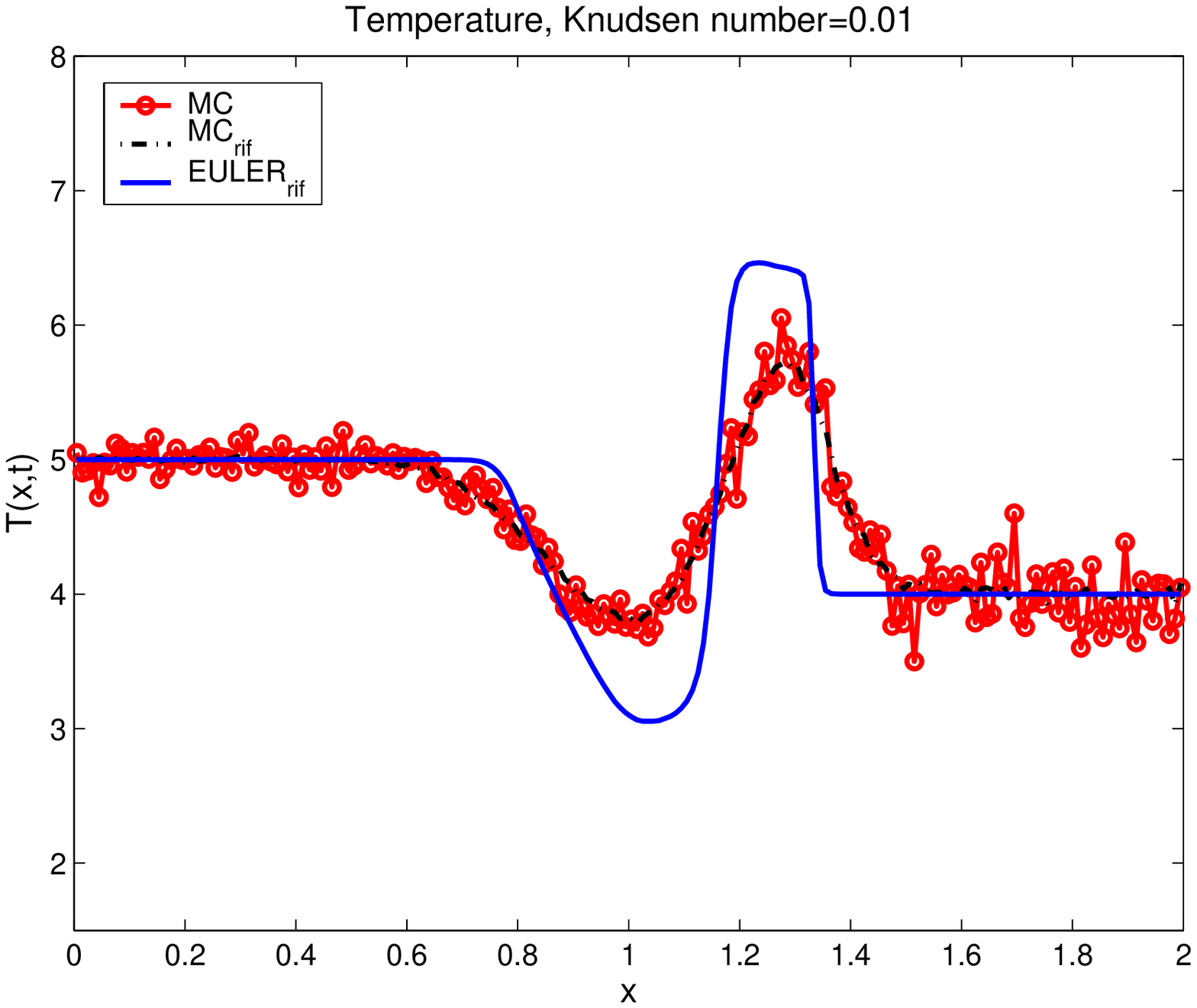}
\includegraphics[scale=0.39]{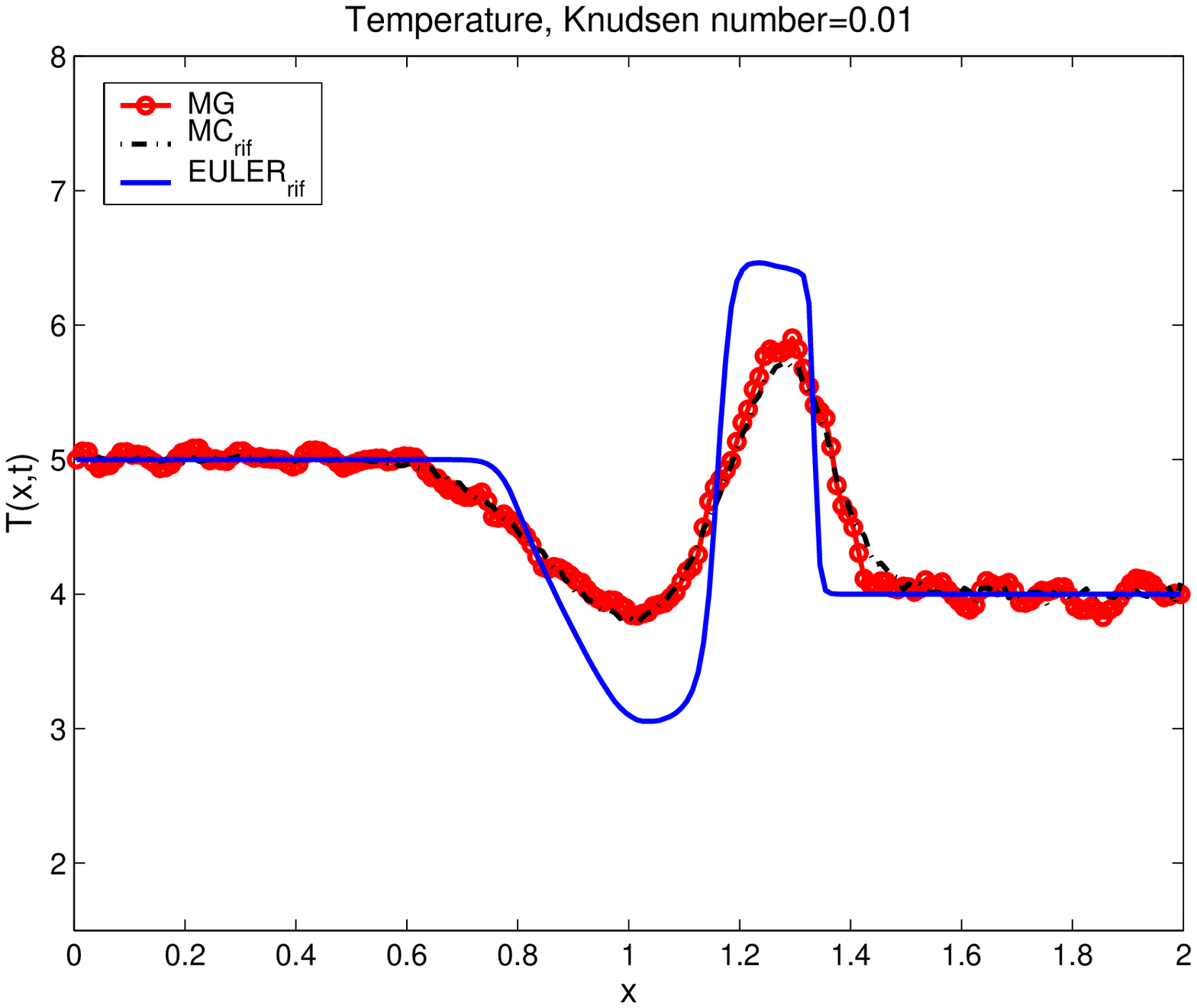}
\caption{Sod shock tube test: Solution at $t=0.05$ for the density
(top), velocity (middle) and temperature (bottom). MC method (left),
Moment Guided MG method (right). Knudsen number
$\varepsilon=10^{-2}$. Reference solution: dash dotted line. Euler
solution: continuous line. Monte Carlo or Moment Guided: circles
plus continuous line. 1000 particles for cell.} \label{So1b}
\end{center}
\end{figure}

\begin{figure}
\begin{center}
\includegraphics[scale=0.39]{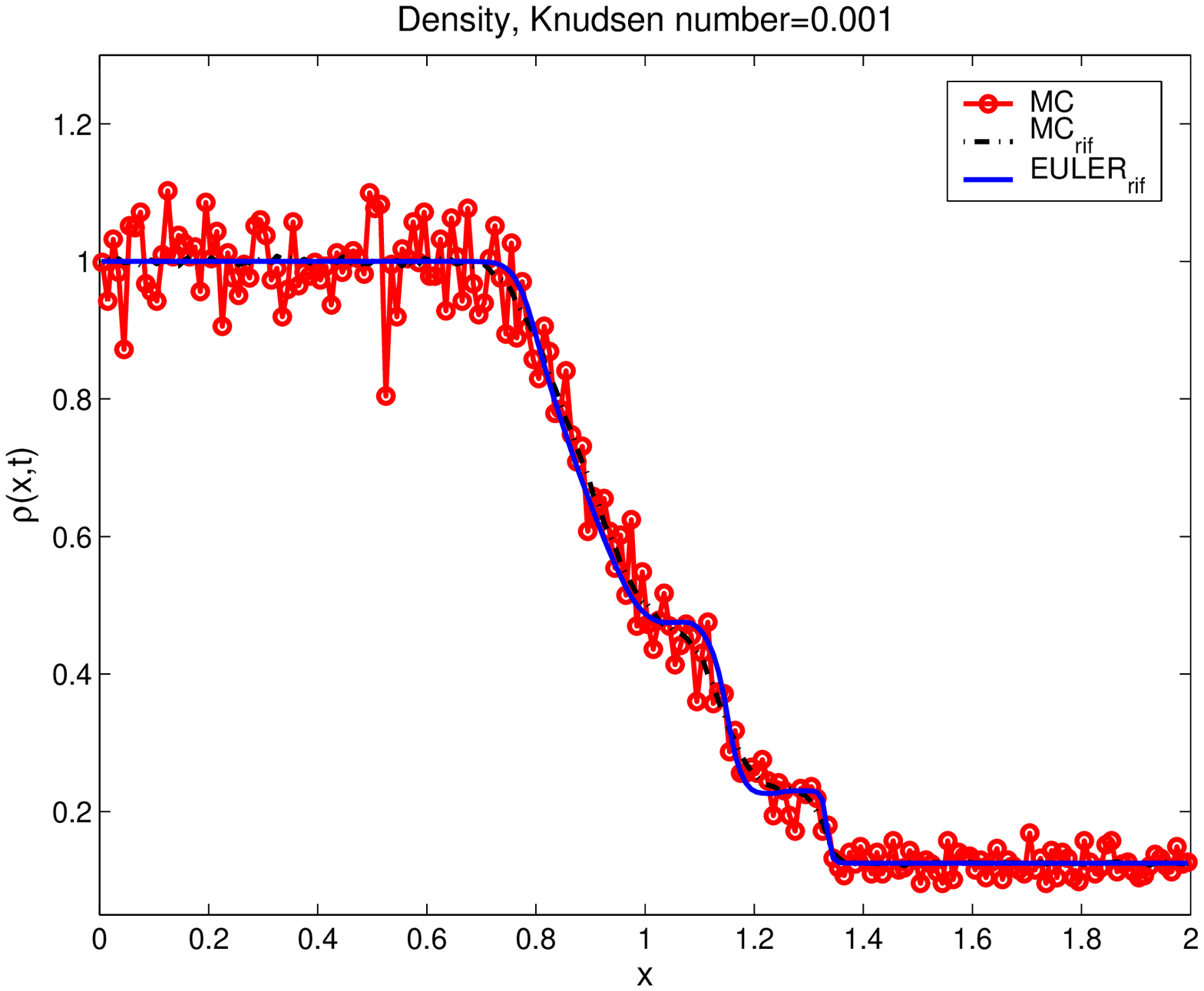}
\includegraphics[scale=0.39]{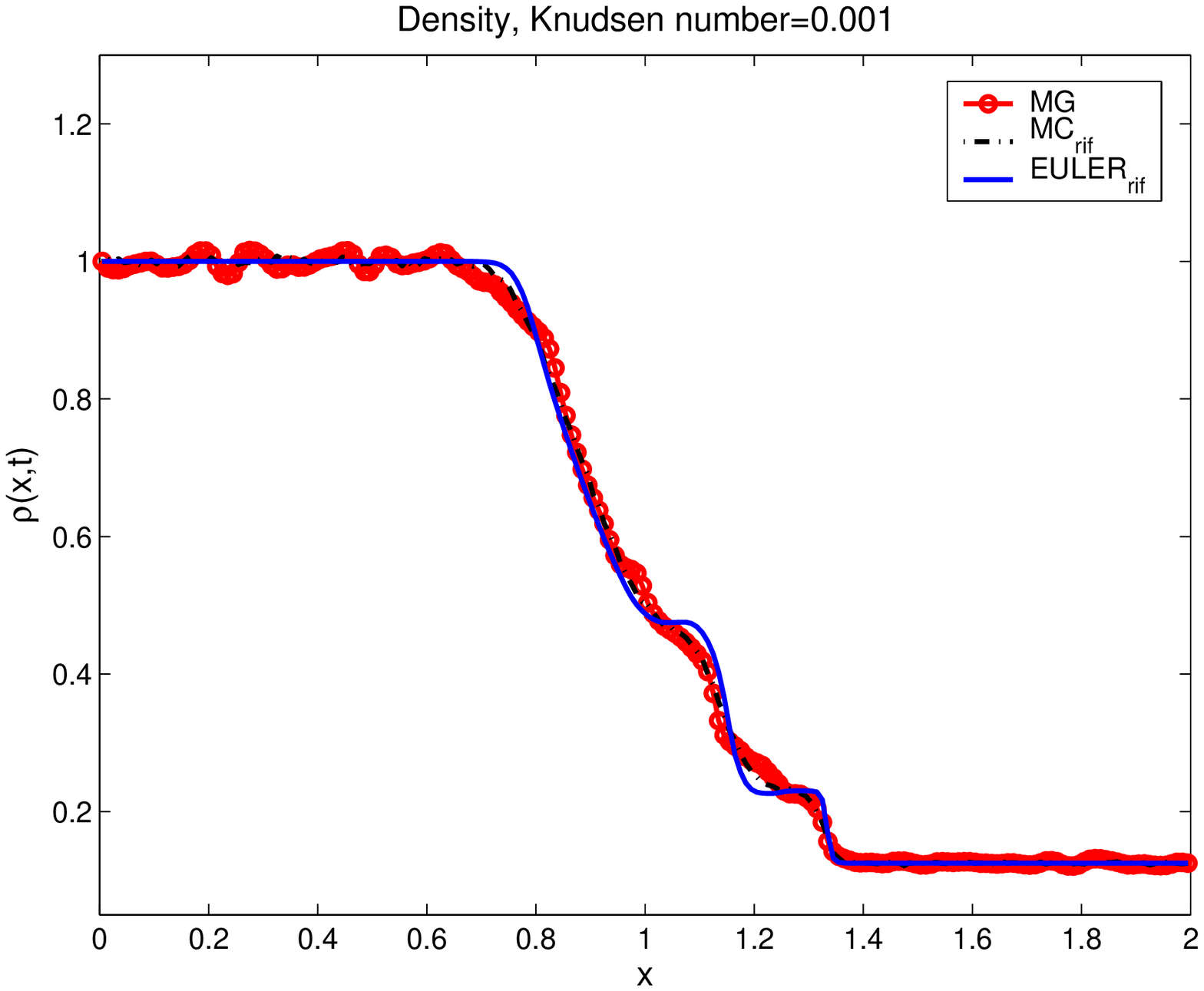}
\includegraphics[scale=0.39]{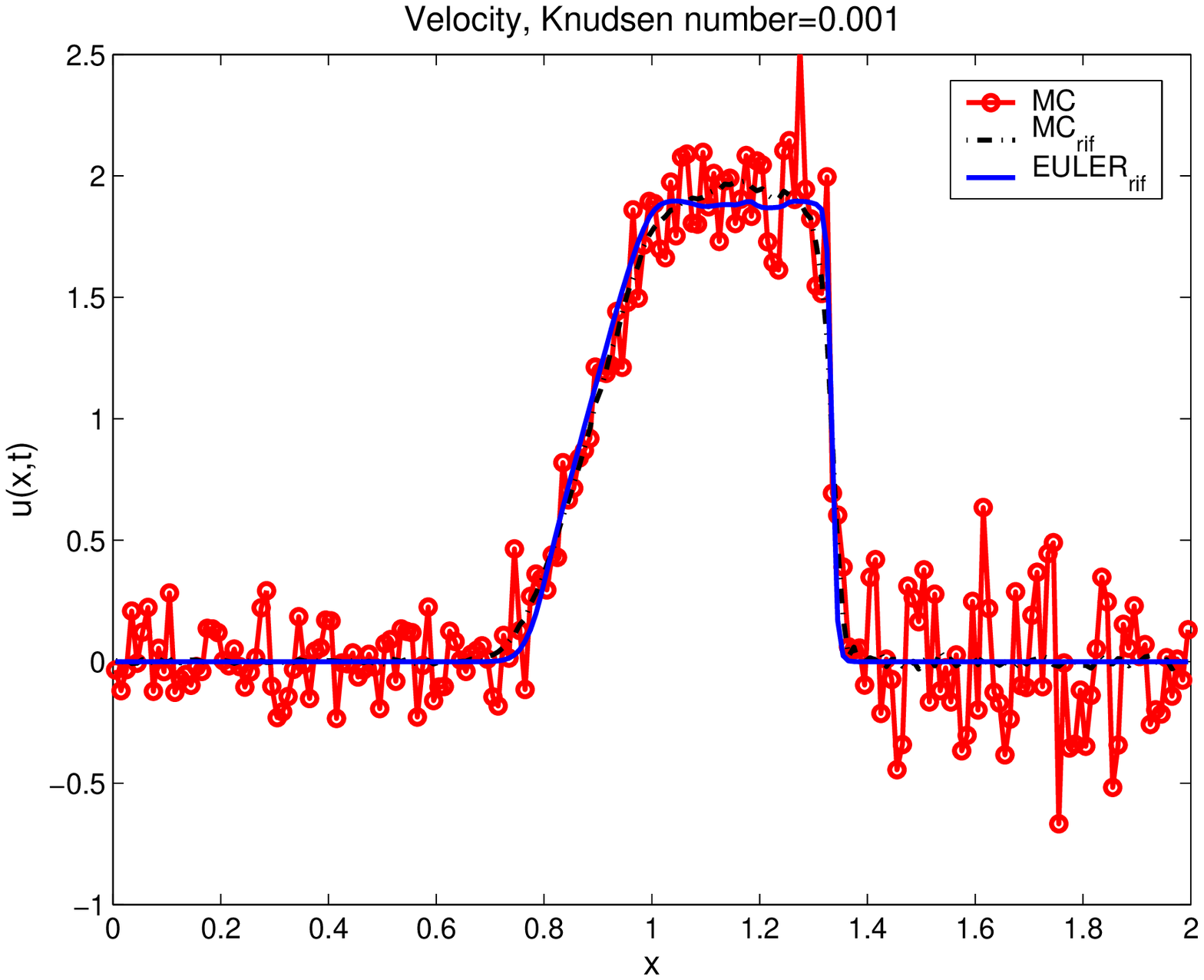}
\includegraphics[scale=0.39]{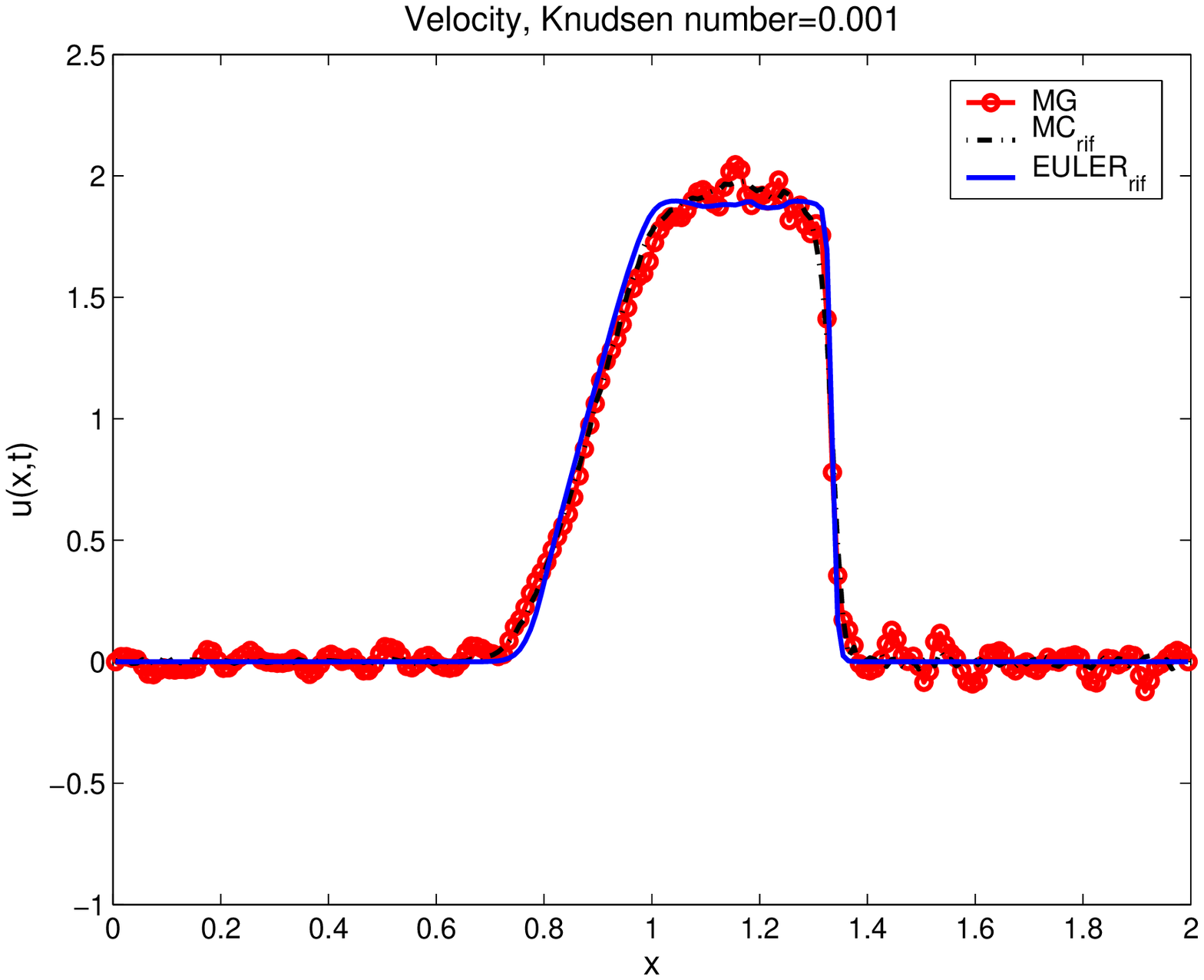}
\includegraphics[scale=0.39]{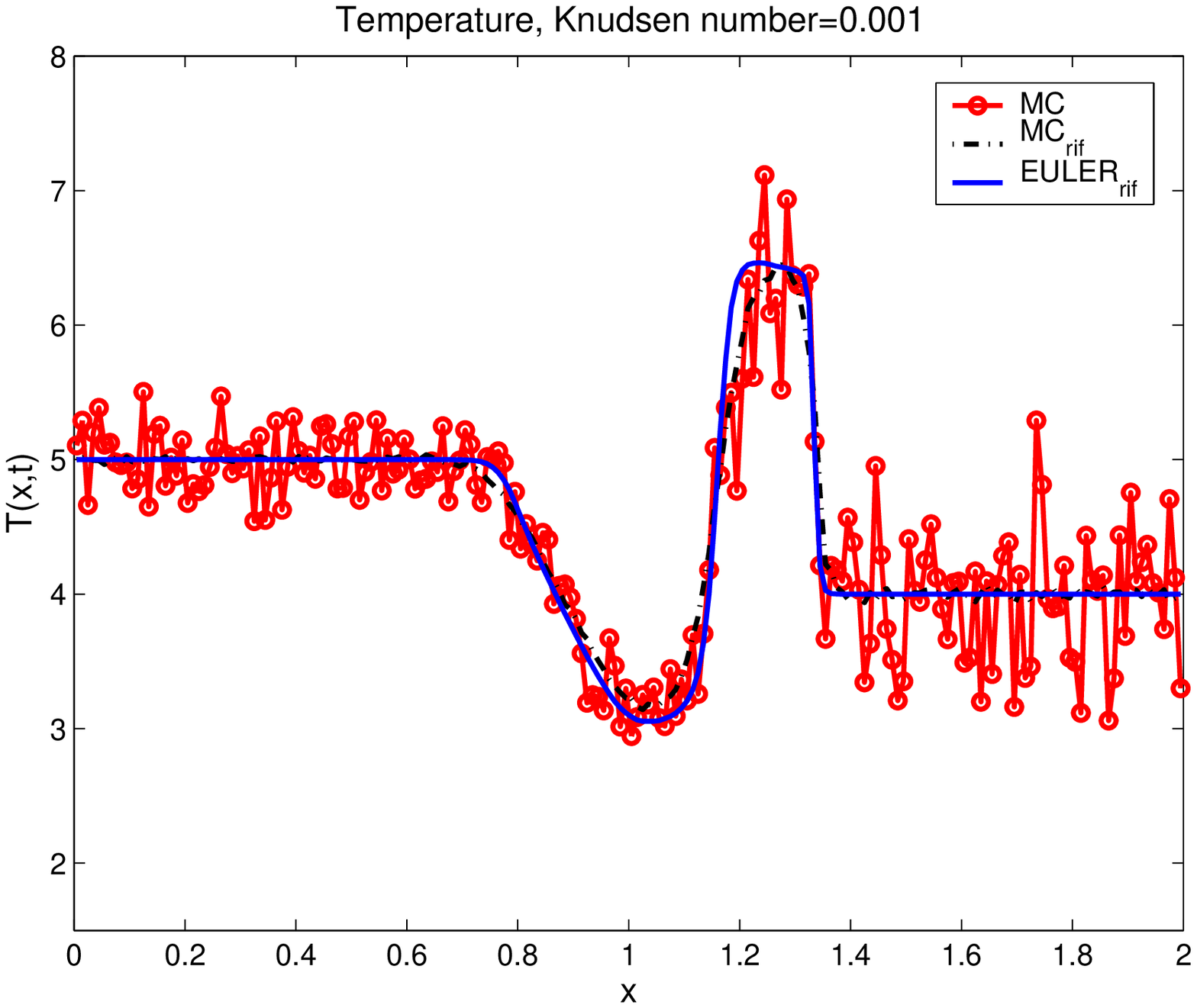}
\includegraphics[scale=0.39]{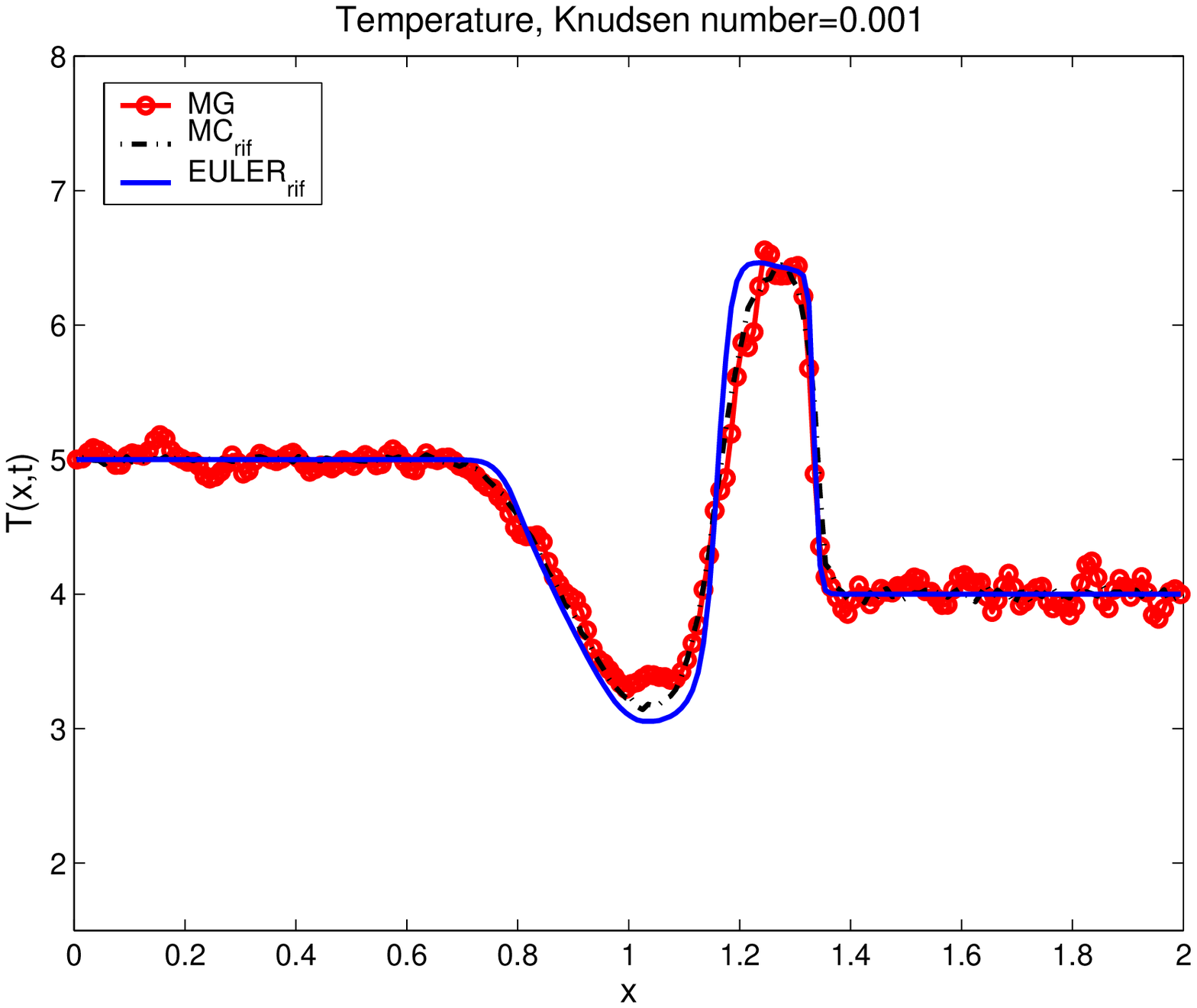}
\caption{Sod shock tube test: Solution at $t=0.08$ for the density
(top), velocity (middle) and temperature (bottom). MC method (left),
Moment Guided MG method (right). Knudsen number
$\varepsilon=10^{-3}$. Reference solution: dash dotted line. Euler
solution: continuous line. Monte Carlo or Moment Guided: circles
plus continuous line. 200 particles for cell.} \label{So2}
\end{center}
\end{figure}
\begin{figure}
\begin{center}
\includegraphics[scale=0.39]{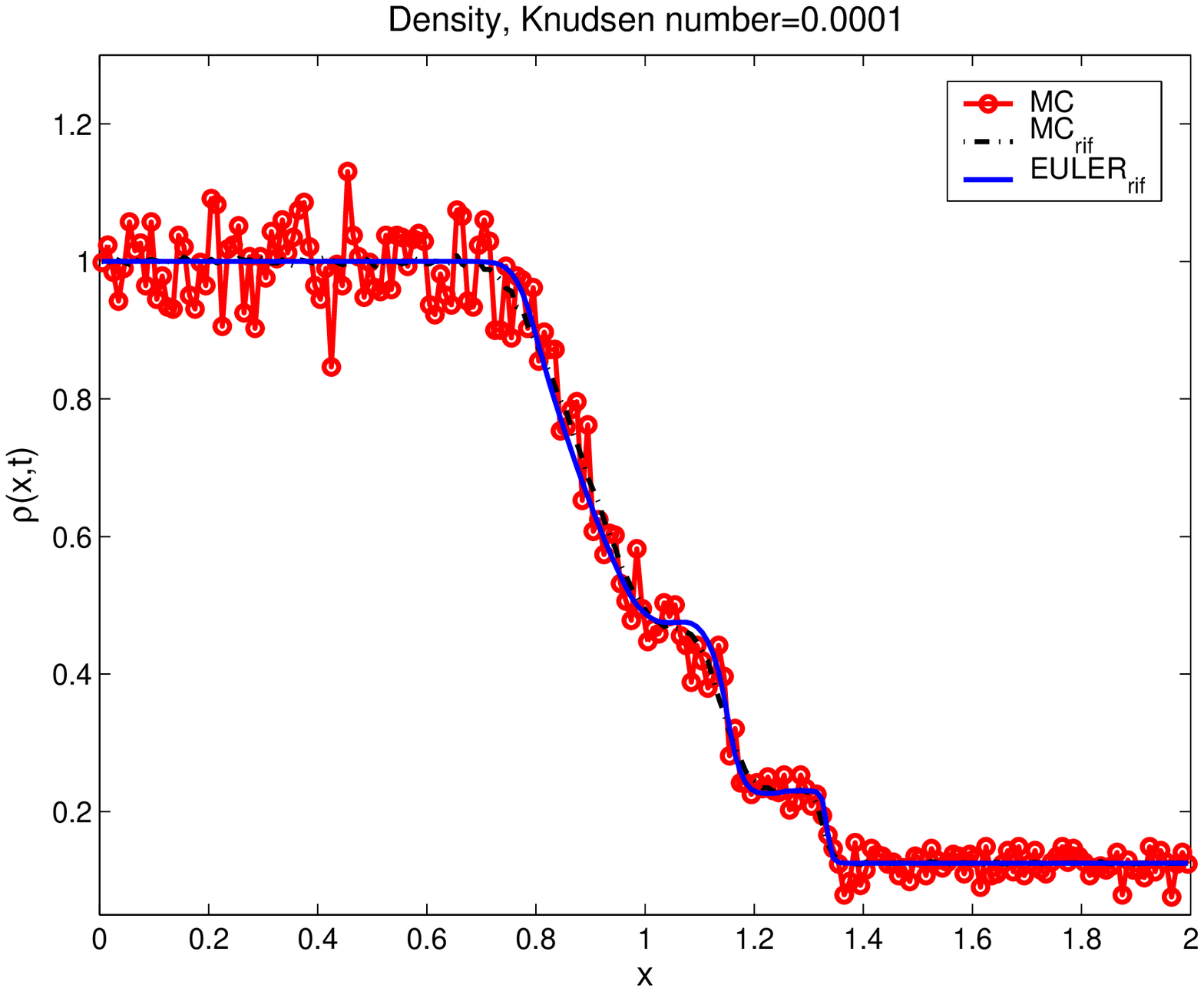}
\includegraphics[scale=0.39]{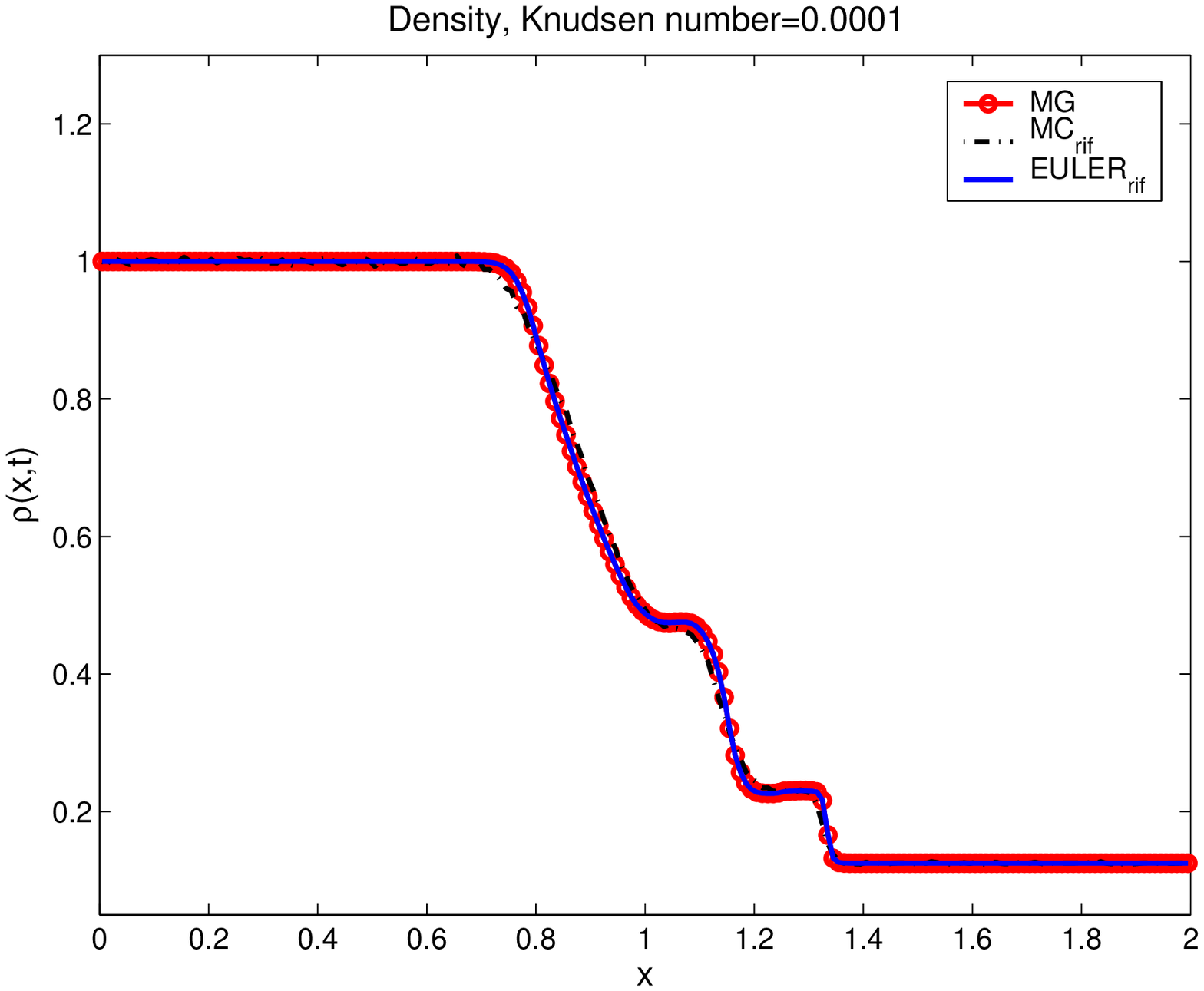}
\includegraphics[scale=0.39]{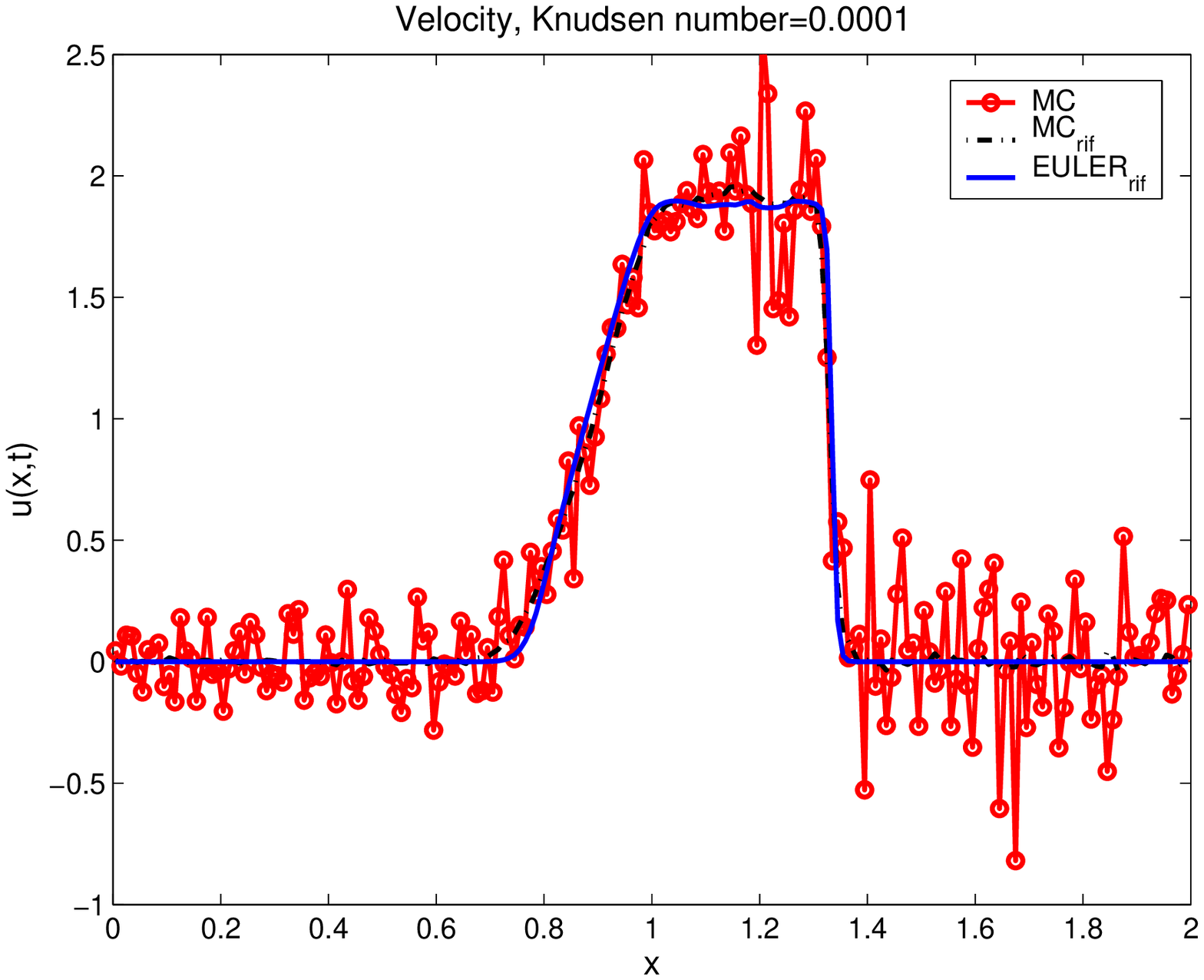}
\includegraphics[scale=0.39]{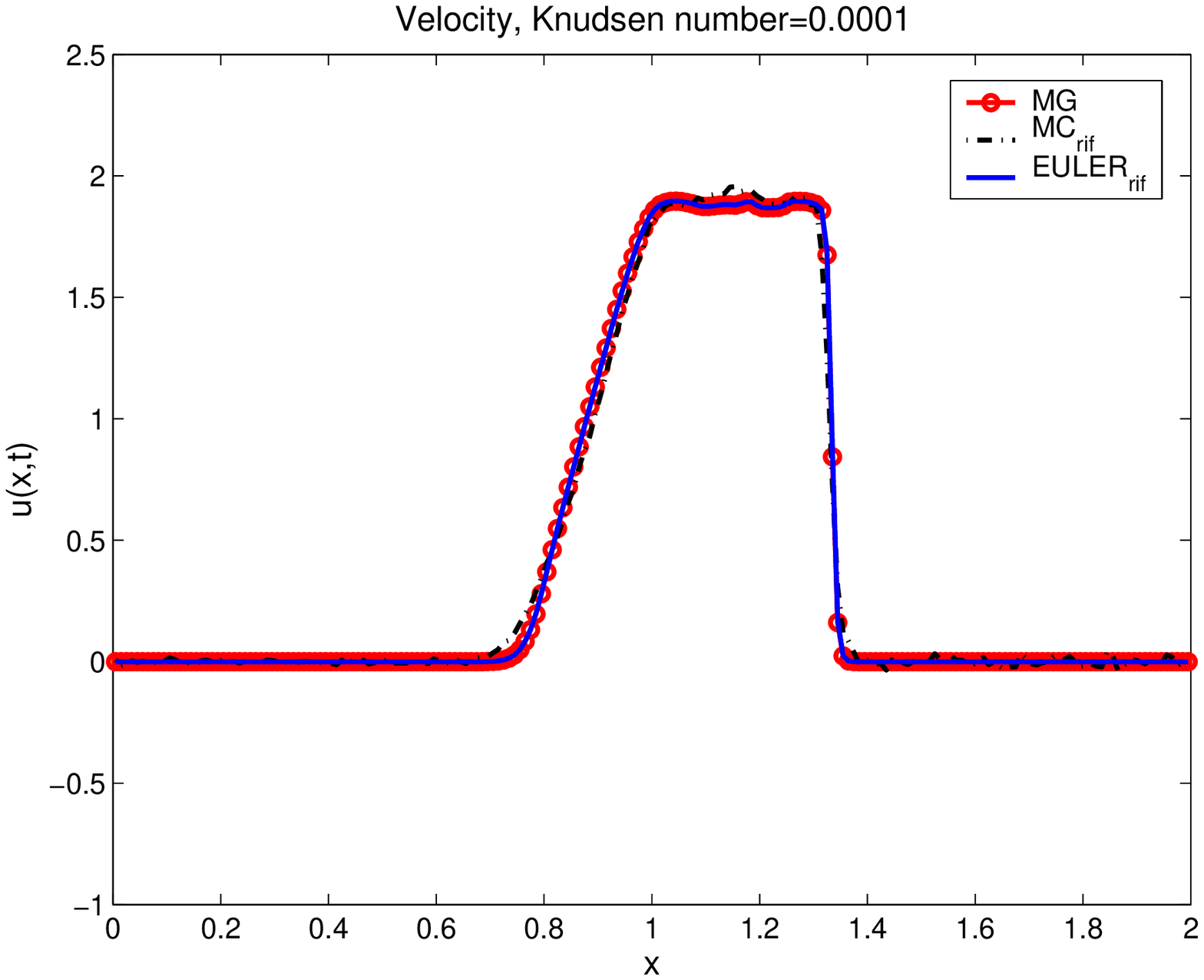}
\includegraphics[scale=0.39]{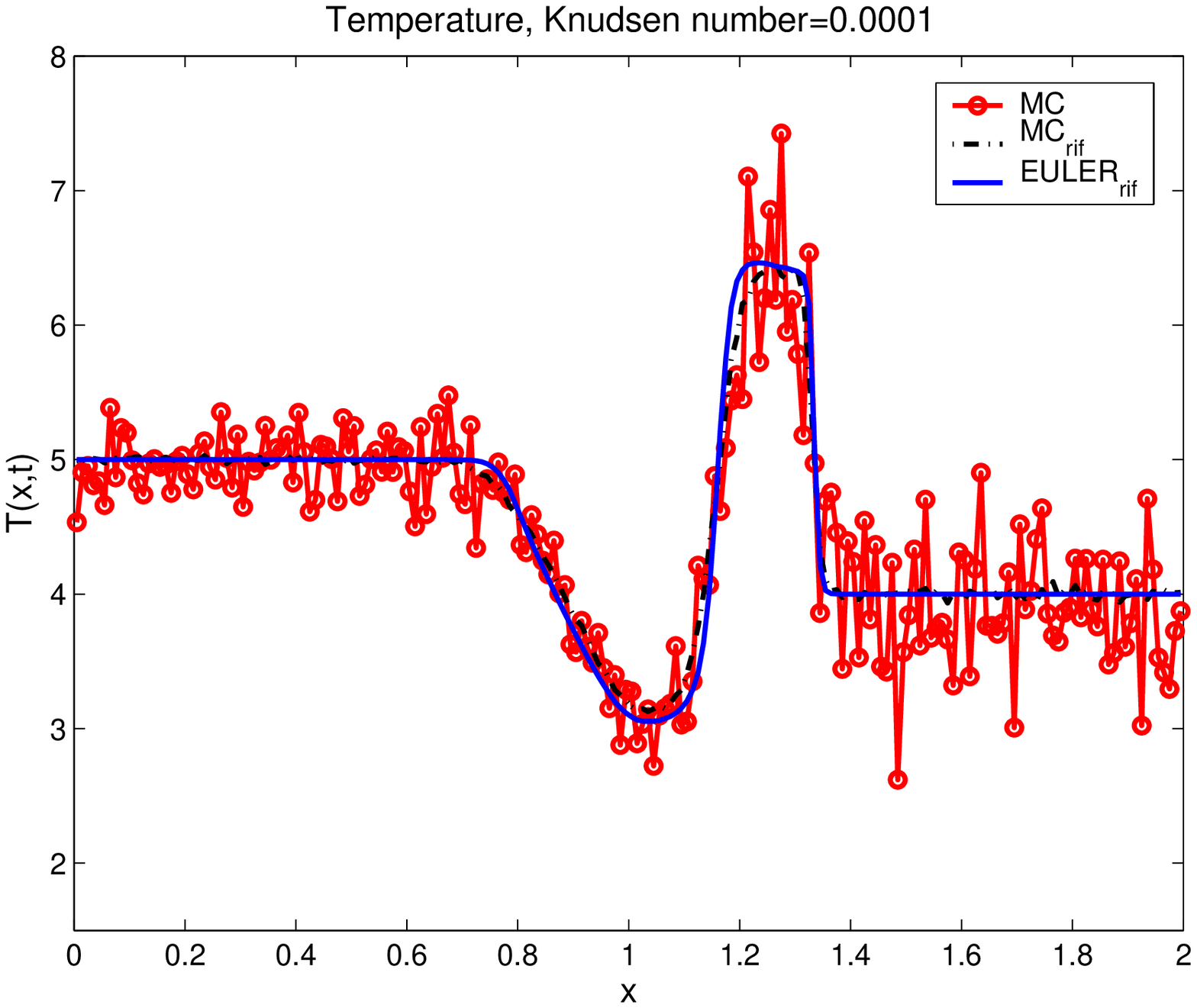}
\includegraphics[scale=0.39]{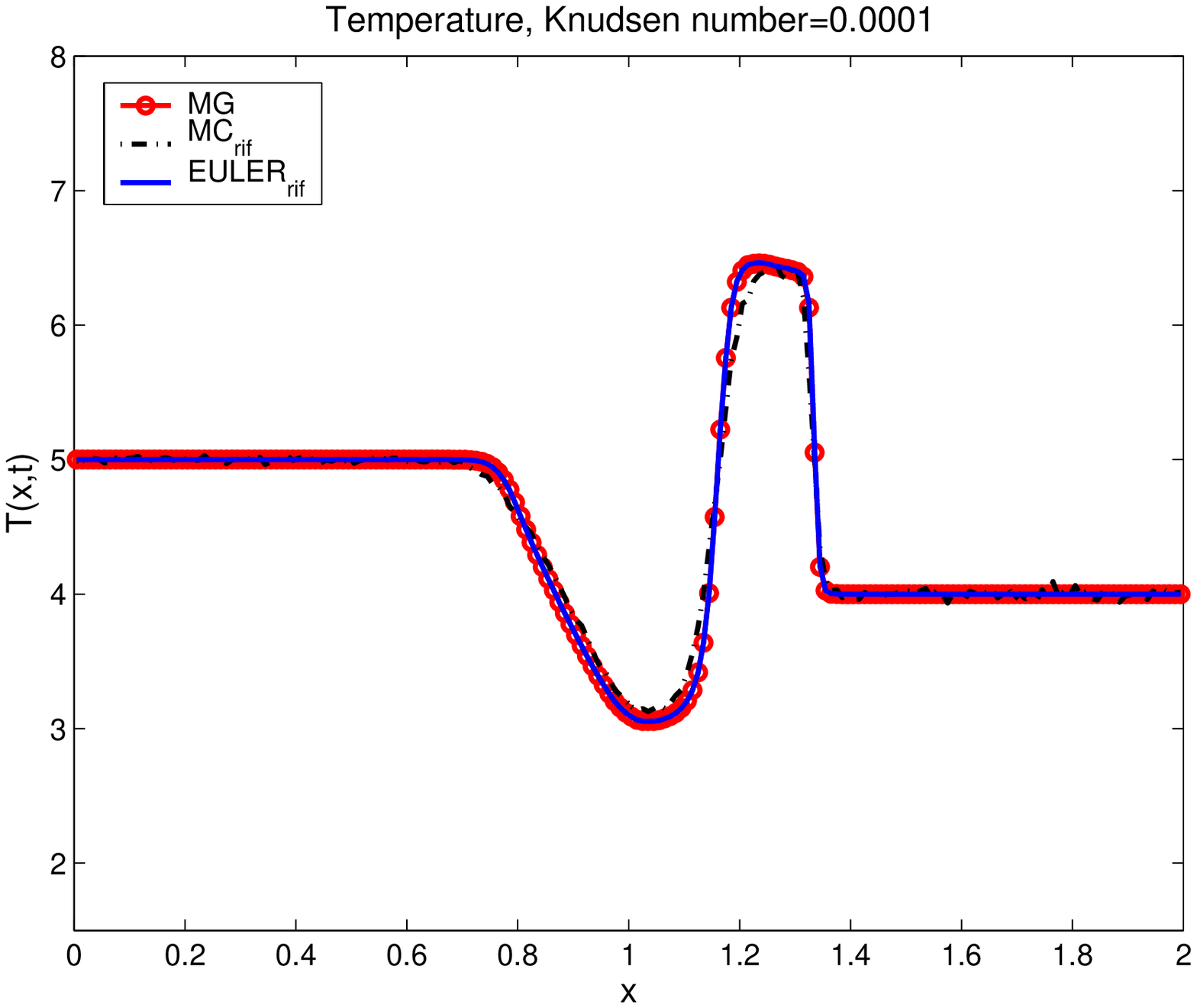}
\caption{Sod shock tube test: Solution at $t=0.08$ for the density
(top), velocity (middle) and temperature (bottom). MC method (left),
Moment Guided MG method (right). Knudsen number
$\varepsilon=10^{-4}$. Reference solution: dash dotted line. Euler
solution: continuous line. Monte Carlo or Moment Guided: circles
plus continuous line. 200 particles for cell.} \label{So3}
\end{center}
\end{figure}

\subsection*{Sod shock tube}

The second problem analyzed is the Sod tube test. Again this
choice relies on the fact that the method is specifically aimed to the study
of unsteady problems. In this case, we want to study the capability of the method in describing different waves with respect to the simple shock wave. The figures \ref{So1} to \ref{So3} consider
the same initial data for the density $\varrho=1$ upstream and
$\varrho=0.125$ downstream of the initial shock, the mean velocity
$u=0$ in all the domain and the temperature: $T=5$ upstream and
$T=4$ downstream of the shock. Different initial Knudsen number
values are considered in the test, which range from
$\varepsilon=10^{-2}$ to $\varepsilon=10^{-4}$ while the number of
cell is $200$. In figures \ref{So1}, \ref{So2} and \ref{So3} $200$ particles per cell are also used on average, more
precisely $200^{^2}$ particles are present at the beginning of the
simulation and then distributed accordingly to the density in each
cell. The final time is $T_{fin}=0.08$, while the time step as previously is given by the minimum between the the ratio of $\Delta x$ over the maximum velocity of
the particles and the ratio of $\Delta x$ over the larger
eigenvalue of the compressible Euler equations. Again the Knudsen number value does not constrain the time step, thanks to the exponential method employed in the resolution of the collision integral.

As for the previous test each figure depicts the density, the mean velocity and
the temperature from top to bottom, with the Monte Carlo solver
(left) and the Moment Guided method (right). The reference solution is
obtained through the Monte Carlo method in which the number of particles is taken equal to $2 \ 10^{7}$ with a solution averaged $5$ times, this gives
a final number of $10^{8}$ particles employed for constructing the reference solution.
The figures show good results for all
ranges of Knudsen numbers in terms of reduction of fluctuations. In particular, we clearly see a strong reduction of fluctuations for $\varepsilon=10^{-3}$ (figure \ref{So2}), while for $\varepsilon=10^{-4}$ (figure \ref{So3}) the fluctuations are almost completely disappeared. In the case of $\varepsilon=10^{-2}$ (figure \ref{So1}), we see that even if the statistical error is smaller for the Moment Guided method than for the Monte Carlo method, both solutions are still far from the converged solution indicated by the dash-dotted line. Thus, we reported the same simulation in figure \ref{So1b} in which $1000$ particle per cell on average are employed for both methods.
This figure permits to show the faster convergence of the Moment Guided method towards the reference
solution also in the case of larger Knudsen. Finally,
thanks to the high order discretization of the moments equations, observe that when the gas is close to the
fluid limit (figure \ref{So3}), the shock is very well represented.

\subsection*{Accuracy test}
\begin{figure}
\begin{center}
\includegraphics[scale=0.39]{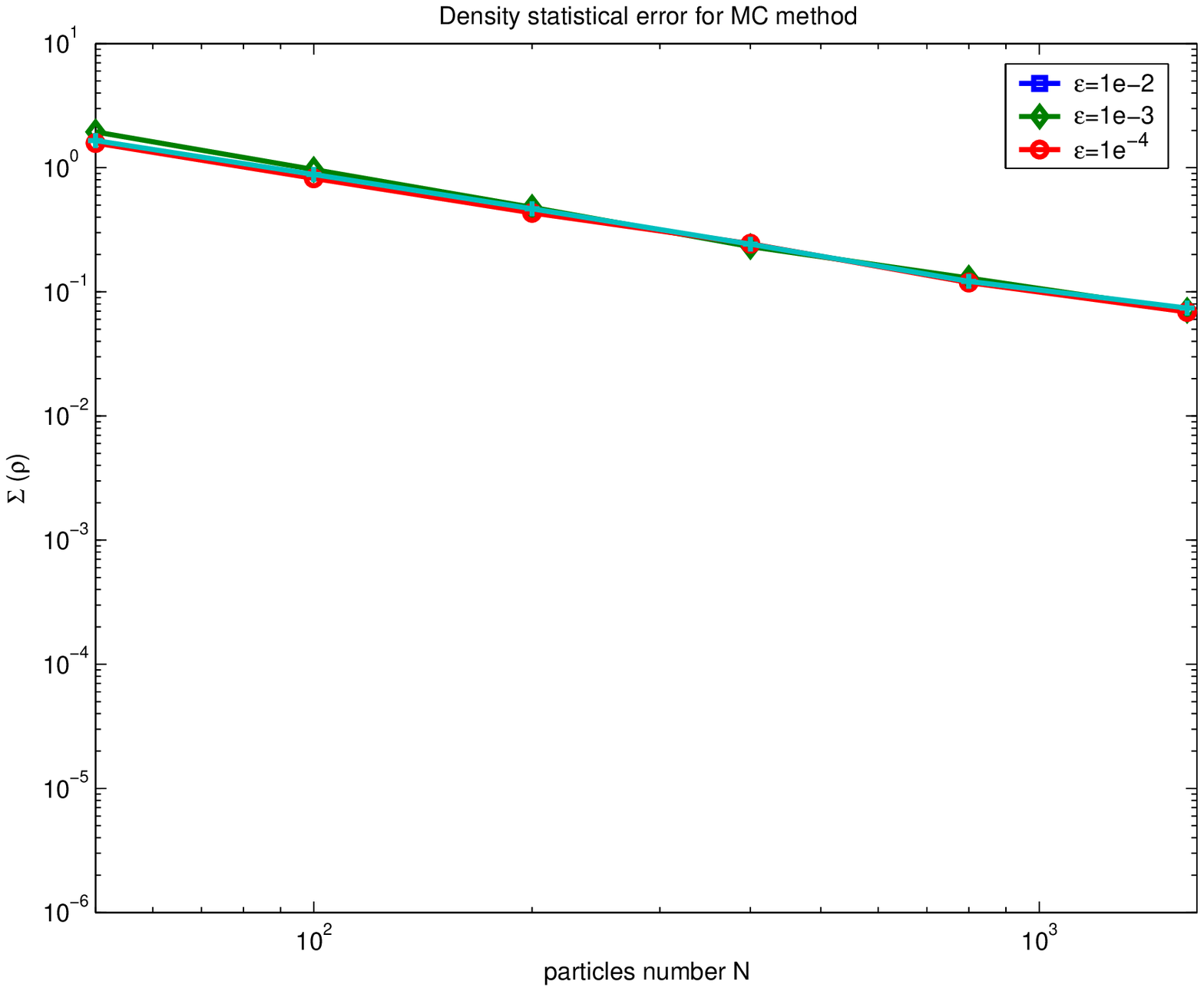}
\includegraphics[scale=0.39]{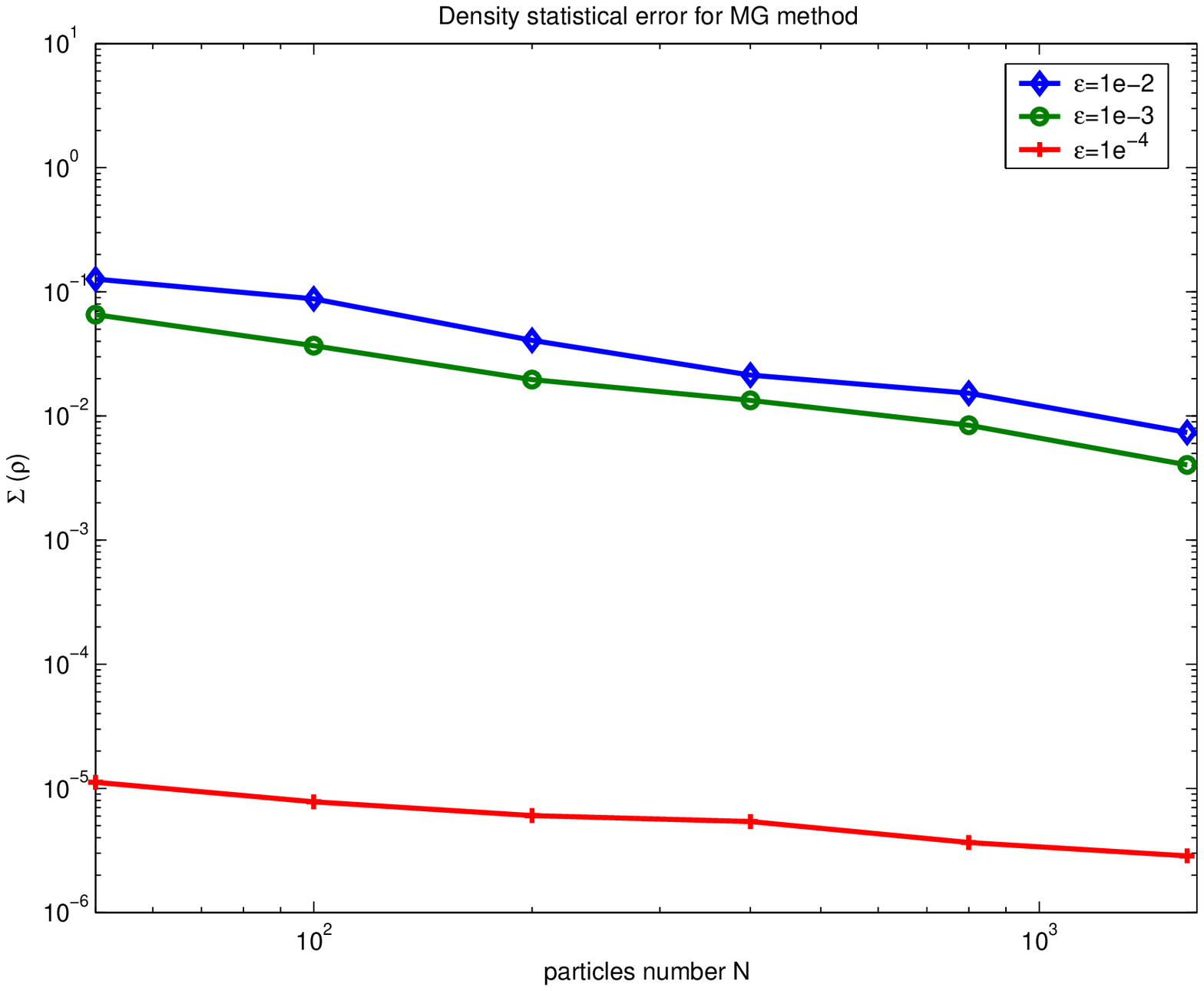}
\includegraphics[scale=0.39]{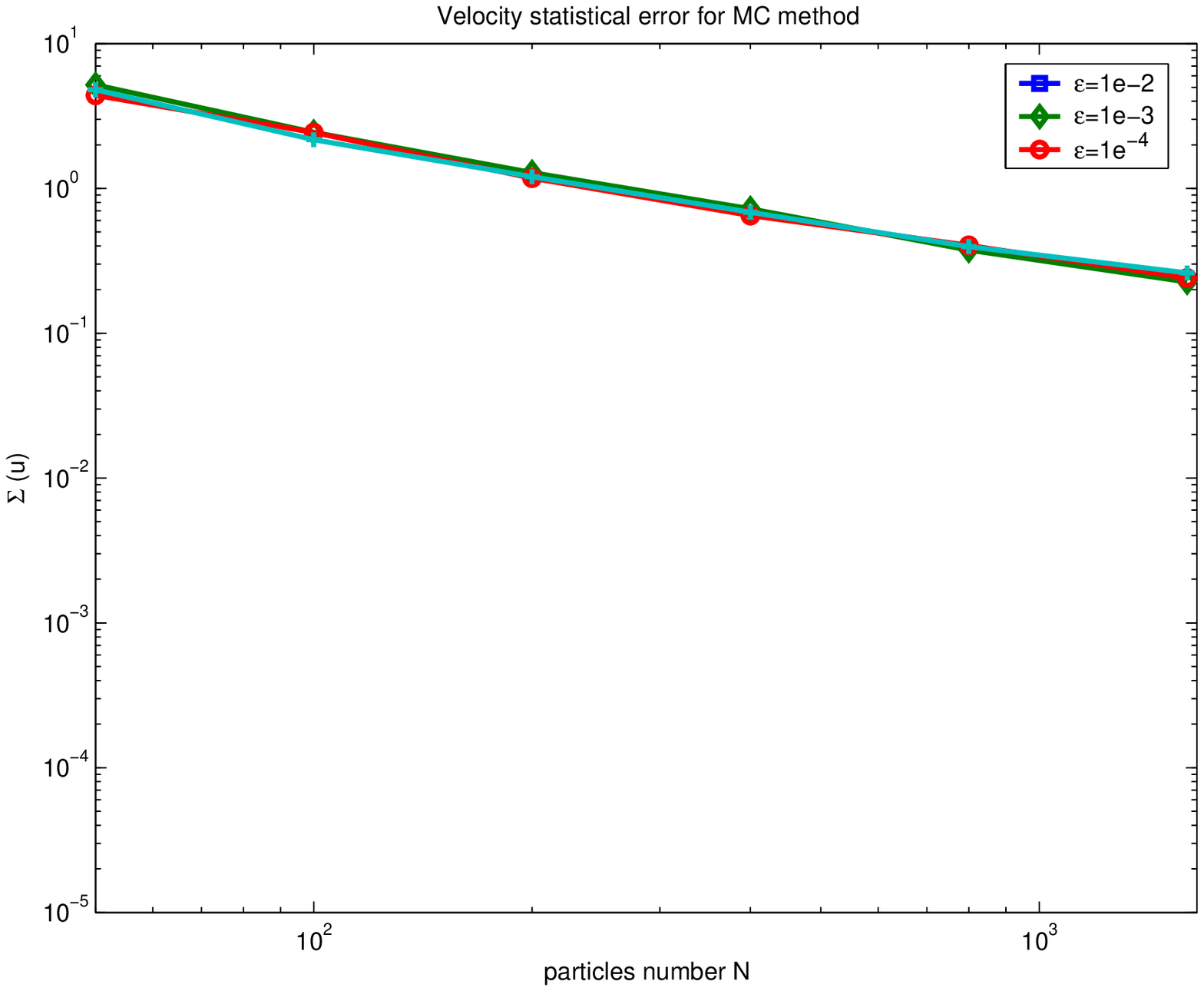}
\includegraphics[scale=0.39]{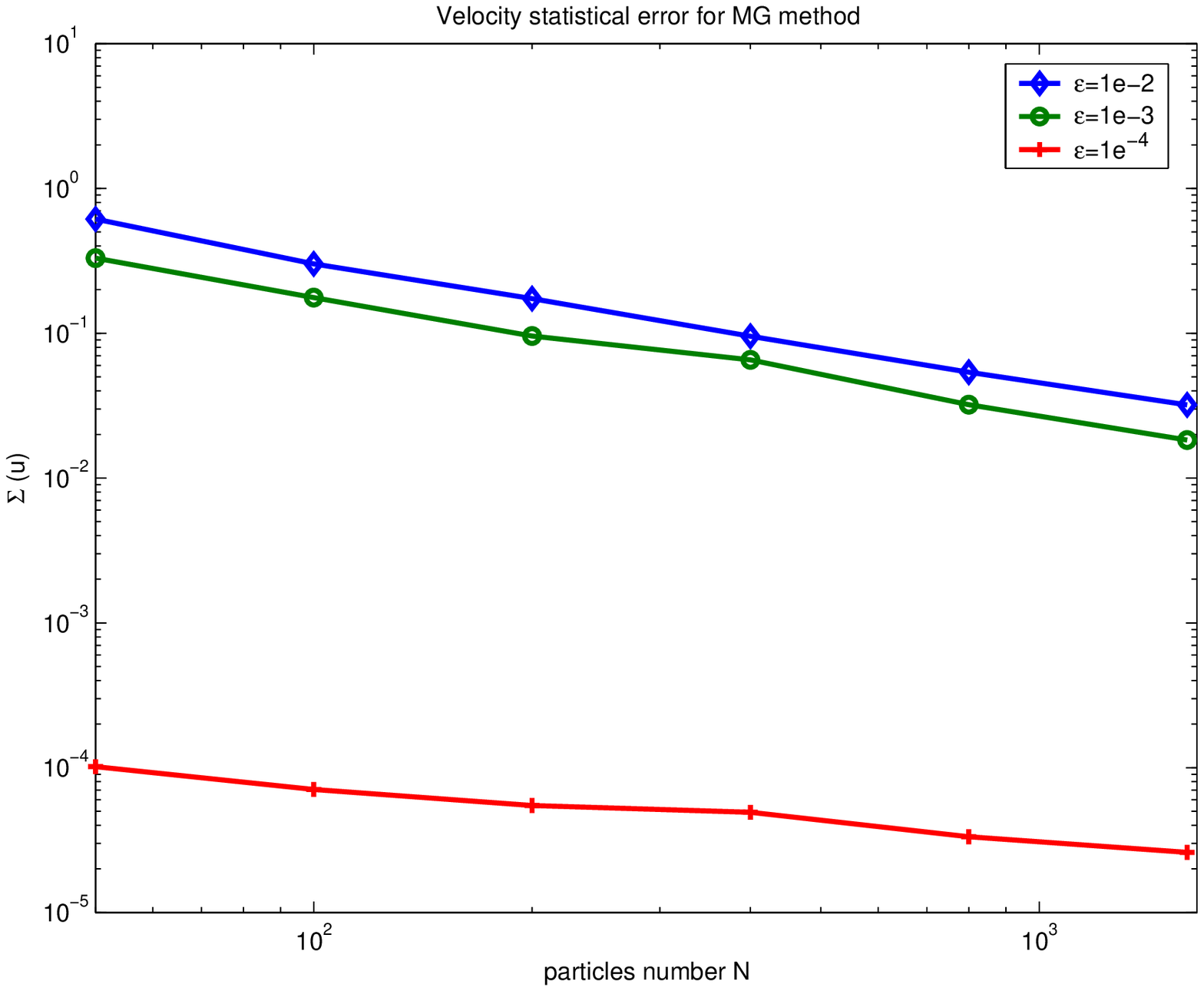}
\includegraphics[scale=0.39]{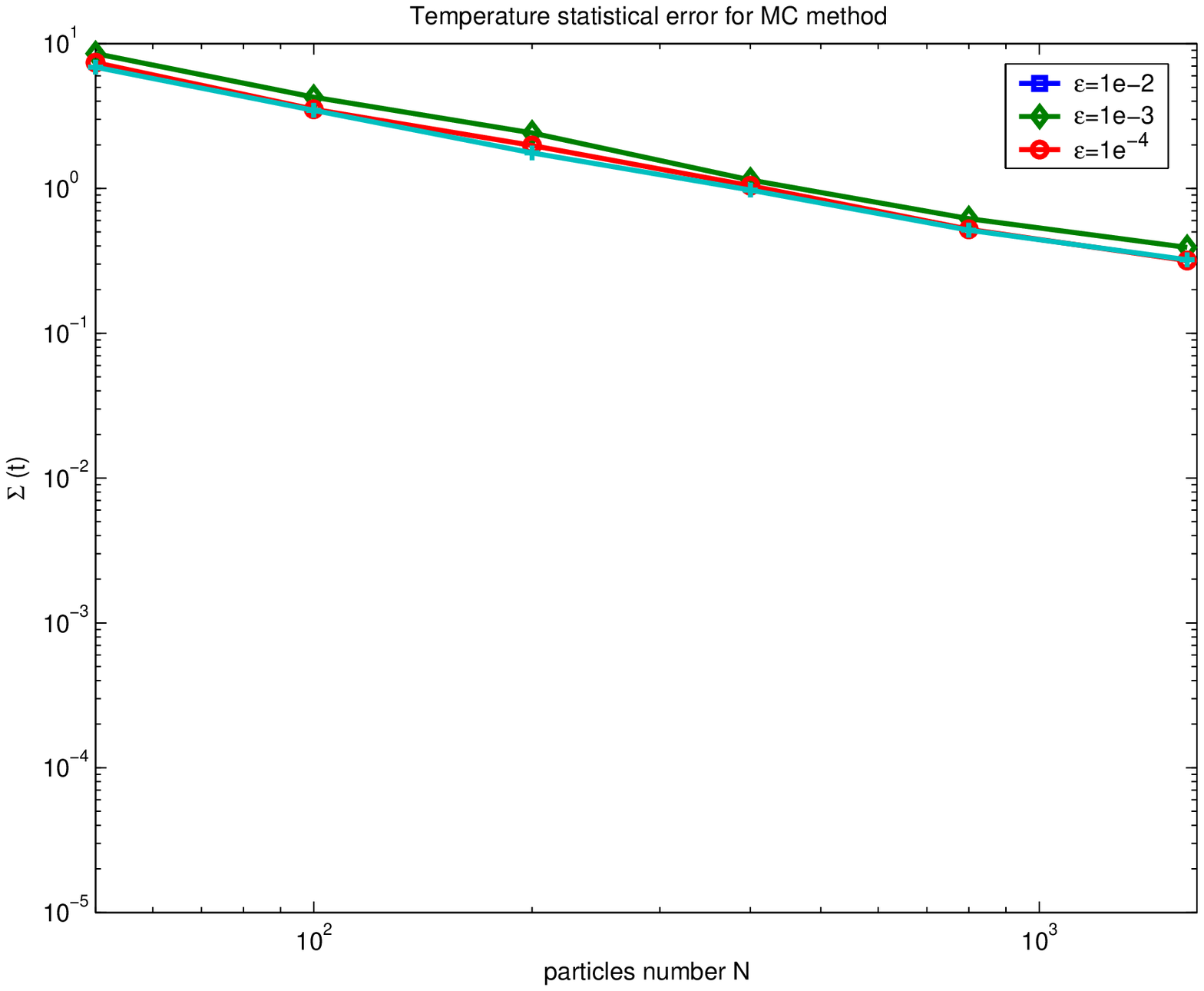}
\includegraphics[scale=0.39]{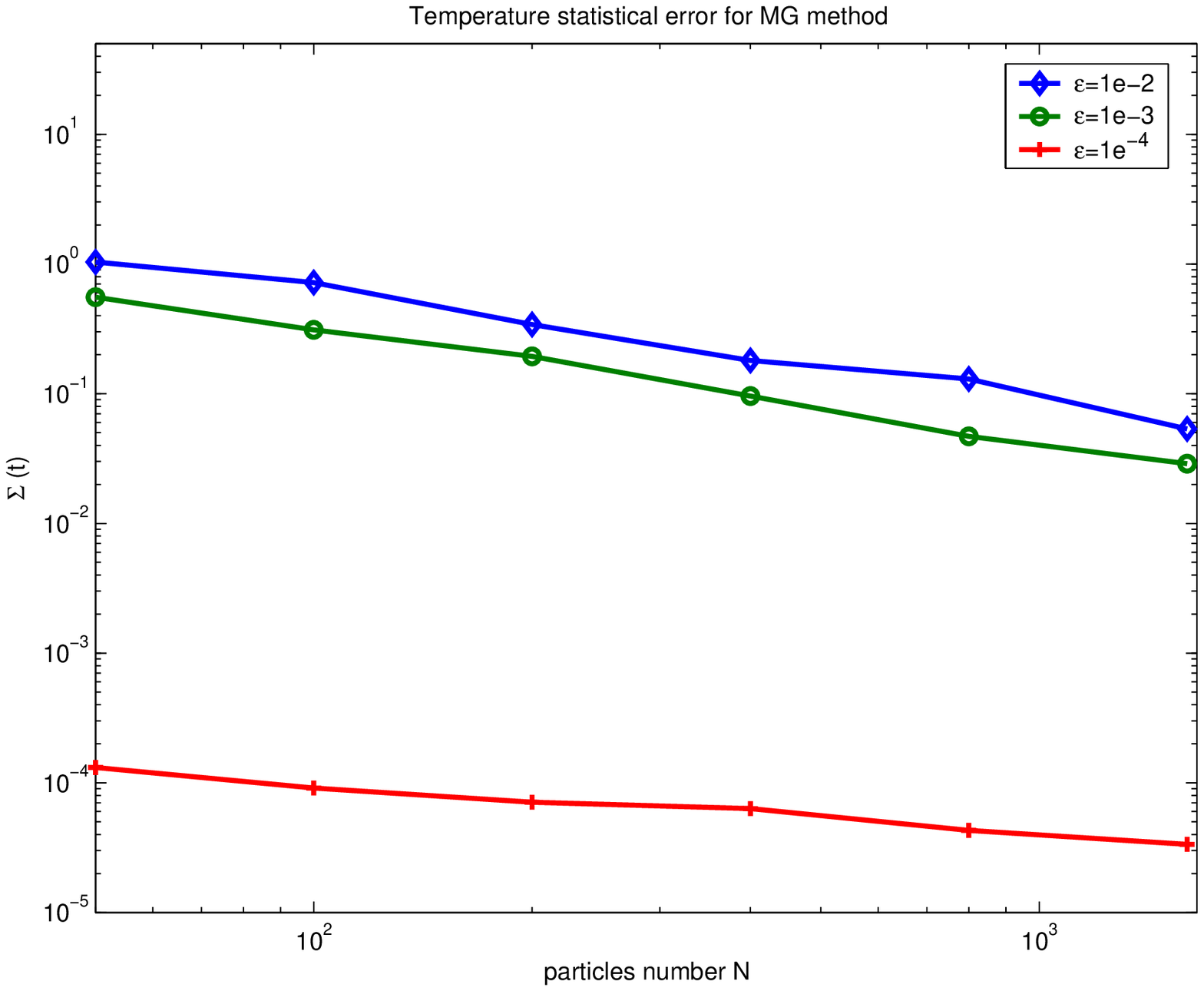}
\caption{Statistical error test: Solution at $t=0.05$ for density
(top), velocity (middle) and temperature (bottom). MC method (left),
Moment Guided MC method (right). Knudsen number vary from
$\varepsilon=10^{-2}$ to $\varepsilon=10^{-4}$. Squares indicate
errors for $\varepsilon=10^{-2}$, diamonds for
$\varepsilon=10^{-3}$, circles for $\varepsilon=10^{-4}$.} \label{SE}
\end{center}
\end{figure}

In this final part, we report the results of a
an error analysis with respect to the number of particles.
In this case the reference solution is obtained by averaging $M$ independent
realizations for the for the Monte Carlo method (subscript $MC$) \be
\overline{U}_{MC}=\frac{1}{M}\sum_{k=1}^{M}U_{k,MC}\ee and for the Moment Guided
method (subscript $MG$) \be
\overline{U}_{MG}=\frac{1}{M}\sum_{k=1}^{M}U_{k,MG}.\ee  Two
different reference solutions have been used since the two schemes present
different discretization errors and thus we expect them to converge, when the
number of particles goes to infinity, to different discretized (in time and space)
solutions. The difference between the two limit
solutions (Monte Carlo and Moment Guided Monte Carlo) is mainly due
to the different numerical diffusion introduced by the two
methods. This has been made clear for instance in the unsteady shock test case when $\varepsilon=10^{-4}$. In this case, as already explained, the Monte Carlo solution is more diffusive with respect to the Moment Guided one and thus, in order to measure only the statistical error, we cannot consider the same reference solution for the two schemes. The two
reference solutions have been obtained by fixing the time step and mesh
size and letting the number of particles grows to eliminate the statistical error. The
number of samples used to compute the reference solution is on average
$10^{5}$ per cell, while the number of realizations is M=5. The number of spatial cells is $100$, the time step
is fixed to $10^{-3}$ for all the values of the Knudsen number. Both reference solutions still contain space and time
discretization errors but the amount of such errors does not
change when the number of particles varies. Therefore, the error we measure by fixing the space and time discretization is only due to the stochastic nature of the method.

In practice, in our test we measured the quantity
\be
\Sigma^2(N)=\frac{1}{M}\sum_{k=1}^{M}\sum_{j=1}^{N}(U_{k,j}-\overline{U_{j}})^2\ee
where $\overline{U_j}$ represents the reference solutions,
$N=100$ and $M=10$. The
following initial data have been used \be \nonumber \varrho(x,0)=1+a_\varrho
\sin\frac{2\pi x}{L} \ee \be u(x,0)=1.5+a_u \sin\frac{2\pi x}{L} \ee
\be \nonumber E(x,0)=2.5+a_W
\sin\frac{2\pi x}{L} \ee with \be \nonumber a_\varrho=0.3 \ \
a_u=0.1 \ \ a_W=1. \ee The boundary conditions are periodic while the time interval in which the simulation has been performed is $t \in [0,5\times 10^{-2}]$. The
results of this test in log-log scale are shown in Figure \ref{SE}.
On the left, we reported the error for the Monte
Carlo while on the right for the Moment Guided method. For the Monte Carlo
method, as expected the stochastic error does not change with
respect to the Knudsen number. At variance, for the Moment Guided
method, the error decreases as the Knudsen number diminishes. In particular observe that for the test studied the stochastic error
is between $8$ and $12$ times smaller of the error of Monte Carlo method for $\varepsilon=10^{-2}$ and $\varepsilon=10^{-3}$ while is $10^{-4}$ smaller when $\varepsilon=10^{-4}$.

\section{Conclusions} We have extended a new class of hybrid methods, first developed in \cite{dimarco1},
which aim at reducing the variance in Monte Carlo schemes to the general case of the Boltzmann equation.
The key idea consists in closing the set of macroscopic moments equations through a coupling with the kinetic equation solved by means of Monte Carlo methods. In order to correctly close the macroscopic system, we need to constrain particle positions and velocities in such a way that the moments given by the solution of the Boltzmann equation exactly match the moments given by the solution of the macroscopic set of moment equations. The new Moment Guided method is constructed in such a way that it avoids the problem of stiff regimes or equivalently of very small Knudsen numbers. In particular, the method benefits of the asymptotic preserving property which permits to circumvent the small scale resolution capturing the limit solution of the problem.
A third important achievement consists in progressively reduce the contribution of the kinetic solution when the Knudsn number tends to zero and to automatically obtain in the fluid limit an high order method for the compressible Euler equations without any stochastic contribution.

The numerical results performed show large reductions of fluctuations in all
regimes, from dense to rarefied, compared to Monte Carlo methods for the Boltzman equation. The reduction becomes stronger as we
approach equilibrium. In addition when $\varepsilon=0$ the stochastic error becomes automatically zero. The numerical convergence tests show that we effectively get better
performances for the Moment Guided method with respect to
pure Monte Carlo schemes for all the cases studied. In particular, for unsteady problems, in which time averaging techniques lose their efficiency the method developed in this paper seems very promising. The computational cost of the method is
comparable to the cost of a traditional Monte Carlo solver with the
addition of the cost of a macroscopic solver for the compressible
Euler equations. This latter is usually computationally much less expensive
than any type of solver applied to kinetic equations.

Currently, we are working on extensions of the present method using higher order closures of
the moment hierarchy in order to solve a larger set of hydrodynamics
equations and thus to additionally reduce the stochastic errors. We are also working to the extension of the method to plasma physics problems such as the Vlasov and the collisional Vlasov models.

\vskip 1cm \textbf{Acknowledgement}. This work was supported by the ANR Blanc project BOOST. The author would like to thank Prof.
L. Pareschi and Prof. P. Degond for the stimulating discussions.



\begin{thebibliography} 
\small



%

\bibitem{Babovsky}
{\sc H.~Babovsky}, {\it On a simulation scheme for the Boltzmann
equation}, Math. Methods Appl. Sci., 8 (1986), pp.~223--233.


\bibitem{bird}
{\sc{G.A.Bird}}, {{Molecular gas dynamics and direct simulation
of gas flows}}, Clarendon Press, Oxford (1994).

\bibitem{birsdall}
{\sc{C.K. Birsdall, A.B. Langdon}}, {\it{Plasma Physics Via Computer
Simulation}}, Institute of Physics (IOP), Series in Plasma Physics
(2004).



\bibitem{BLT}
{\sc J.~F.~Bourgat, P.~LeTallec, M.D.~Tidriri}, {\it Coupling
Boltzmann and Navier-Stokes equations by friction}. J. Comput. Phys.
127, vol. 2 (1996), pag.~227--245.


%

\bibitem{Boyd2}
{\sc J. Burt, I. Boyd}, {\em A hybrid particle approach for
continuum and rarefied flow simulation}, J. Comput. Phys., Vol. 228,
(2009), pp. 460-475.




\bibitem{Cf}
{\sc R.~E.~Caflisch}, {\em Monte Carlo and Quasi-Monte Carlo
Methods}, Acta Numerica (1998)  pp.~1--49.

\bibitem{CPima}
{\sc R.~E.~Caflisch, L.~Pareschi}, {\em Towards an hybrid method for
rarefied gas dynamics}, IMA Vol. App. Math., vol. 135 (2004),
pp.~57--73.


%




\bibitem{cercignani}
{\sc C.~Cercignani}, {The Boltzmann Equation and Its
Applications}, Springer-Verlag, New York, (1988).

%




%
\bibitem{Lemou}
{\sc A. Crestetto, N.Crouseilles, M.Lemou}, {\em
Kinetic/fluid micro-macro numerical schemes for
Vlasov-Poisson-BGK equation using particles}, To appear on KRM 2012.

\bibitem{degond}
{\sc N. Crouseilles, P. Degond, M. Lemou}, {\em A hybrid kinetic-fluid
model for solving the gas-dynamics Boltzmann BGK equation}, J.
Comput. Phys., vol. 199 (2004), pp. 776-808.

\bibitem{dimarco1}
{\sc P. Degond, G.~Dimarco, L. Pareschi}, {\it The Moment Guided
Monte Carlo Method}, Int. J.
Num. Meth. Fluids, vol. 67 (2011), pp. 189-213.




\bibitem{degond3}
{\sc P. Degond, J.-G. Liu, L. Mieussens}, {\em Macroscopic fluid
models with localized kinetic upscaling effects}, SIAM MMS, vol. 5
(2006), pp. 940-979.

\bibitem{degond4}
{\sc P. Degond, G. Dimarco}, {\it Fluid simulations with localized Boltzmann upscaling by direct Monte Carlo.}
J. Comp. Phys., vol. 231 (2012), pp. 2414-2437.


\bibitem{dimarco3}
{\sc P. Degond, G. Dimarco, L. Mieussens.}, {\it A Multiscale
Kinetic-Fluid Solver With Dynamic Localization Of Kinetic Effects.}
J. Comp. Phys., Vol. 229, pp. 4907-4933, (2010).


\bibitem{dimarco4}
{\sc G.~Dimarco and L.~Pareschi}, {\it Hybrid multiscale methods I.
Hyperbolic Relaxation Problems}, Comm. Math. Sci., 1, (2006),
pp.~155-177.
%

\bibitem{dimarco5}
{\sc G.~Dimarco, L.~Pareschi}, {\it A Fluid Solver Independent
Hybrid method for Multiscale Kinetic Equations}, SIAM J. Sci.
Comput. Vol. 32 issue 2, pp. 603-634 (2010).


\bibitem{dimarco7}
{\sc G.~Dimarco and L.~Pareschi}, {\it Exponential Runge-Kutta methods for stiff kinetic equations}, SIAM
J. Num.  Anal., vol. 49 , pp. 2057-2077 (2011).

\bibitem{Filbet} {\sc F.~Filbet, S.~Jin}, {\it A class of asymptotic preserving schemes for kinetic equations and related
problems with stiff sources}, J. Comp. Phys. 229 (2010), 7625-7648.

\bibitem{Filbet3} {\sc F.~Filbet, S.~Jin}, {\it An Asymptotic Preserving Scheme for the ES-BGK Model of the Boltzmann Equation}.
J. Sci. Comput., Vol. 46 (2011), 204-224.

%
\bibitem{WE1}
{\sc W.~E, B.~Engquist}, {\it The heterogeneous multiscale methods},
Comm. Math. Sci., vol. 1 (2003), pp.~87-133.

%









\bibitem{HashHassan}
{\sc D.B.~Hash and H.A.~Hassan},
\newblock {\em Assessment of schemes for coupling Monte Carlo and Navier-Stokes solution
methods}, J. Thermophys. Heat Transf., 10, (1996), pp.~242--249.


\bibitem{Hadji}
{\sc T. Homolle, N. Hadjiconstantinou}, {\em A low-variance
deviational simulation Monte Carlo for the Boltzmann equation.}
J. Comp. Phys. 226 (2007), pp 2341--2358.

\bibitem{Hadji1}
{\sc T. Homolle, N. Hadjiconstantinou}, {\em Low-variance
deviational simulation Monte Carlo.} Phys. Fluids 19 (2007)
041701.




%
%
%

%
%

\bibitem{Jin2}
{\sc S.~Jin}, {\it Efficient Asymptotic-Preserving (AP) schemes for some multiscale kinetic equations},
SIAM J. Sci. Comput., 21 (1999), 441--454.

\bibitem{JX}
{\sc S.~Jin and Z.~P.~Xin}, {\it Relaxation schemes for systems of
conservation laws in arbitrary space dimensions}, Comm. Pure Appl.
Math., vol. 48 (1995), pp.~235--276.




%
%



%

\bibitem{Letallec}
{\sc P.~LeTallec and F.~Mallinger},
\newblock {\em Coupling Boltzmann and Navier-Stokes by half
fluxes}
\newblock {J. Comput. Phys.}, vol .136  (1997), pp.~51--67.

\bibitem{leveque:numerical-methods} {\sc R.~J.~LeVeque},
   {Numerical Methods for Conservation Laws},
   Lectures in Mathematics, Birkhauser Verlag, Basel (1992).


\bibitem{sliu}
{\sc S.~Liu}, {Monte Carlo strategies in scientific computing},
Springer, (2004).





\bibitem{Nanbu80}
{\sc K.~Nanbu}, {\it Direct simulation scheme derived from the
Boltzmann equation}, J. Phys. Soc. Japan, vol. 49 (1980),
pp.~2042--2049.











%



%
%

%


\bibitem{Pullin78}
{\sc D.~I.~Pullin}, {\it Direct simulation methods for compressible
inviscid ideal gas flow}, J. Comput. Phys., 34 (1980), pp.~231--244.



%





%
\bibitem{tiwari_JCP}
  {\sc S. Tiwari}, {\it Coupling of the Boltzmann and Euler equations with automatic domain decomposition},
  J. Comput. Phys., vol. 144, 1998, 710--726.

\bibitem{tiwari_JCP1}
{\sc S. Tiwari, A. Klar, S. Hardt}, {\it A Particle–-Particle Hybrid
Method for Kinetic and Continuum Equations}, J. Comput. Phys., Vol.
228, (2009) pp. 7109-7124.



%




%
%

\end{thebibliography}
\end{document}